\documentclass[journal]{IEEEtran}

\RequirePackage[numbers]{natbib}
\usepackage{graphicx,colordvi,psfrag}
\usepackage{graphicx,epsfig,rotating,algorithm,algpseudocode}
\usepackage{amsmath,amssymb}
\usepackage{epsfig}
\usepackage{calc,pstricks, pgf, xcolor}
\usepackage{bbding}
\usepackage{bbm}
\usepackage{dsfont}
\usepackage{centernot}
\usepackage[caption=false]{subfig}
\usepackage{setspace}
\usepackage{ifthen}
\usepackage{xcolor}
\RequirePackage[colorlinks,citecolor=blue,urlcolor=blue]{hyperref}
\usepackage{color}
 
\usepackage{lineno}

 \newtheorem{theorem}{Theorem}
\newtheorem{corollary}{Corollary}
\newtheorem{proposition}{Proposition}

\newtheorem{lemma}{Lemma}
	\newtheorem{definition}{Definition}
	\newtheorem{example}{Example}
	\newtheorem{remark}{Remark}
\newcommand{\beq}{\begin{equation}}
	\newcommand{\eeq}{\end{equation}}
\newcommand{\beas}{\begin{align*}}
	\newcommand{\eeas}{\end{align*}}
\newcommand{\bea}{\begin{align}}
	\newcommand{\eea}{\end{align}}
\newcommand{\bei}{\begin{itemize}}
	\newcommand{\eei}{\end{itemize}}
\newcommand{\ben}{\begin{enumerate}}
	\newcommand{\een}{\end{enumerate}}
\newcommand{\bet}{\begin{theorem}}
	\newcommand{\eet}{\end{theorem}}
\newcommand{\bel}{\begin{lemma}}
	\newcommand{\eel}{\end{lemma}}
\newcommand{\bep}{\begin{proposition}}
	\newcommand{\eep}{\end{proposition}}
\newcommand{\bed}{\begin{definition}}
	\newcommand{\eed}{\end{definition}}
\newcommand{\bec}{\begin{corollary}}
	\newcommand{\eec}{\end{corollary}}
\newcommand{\bex}{\begin{example}}
	\newcommand{\eex}{\end{example}}

\newcommand{\bu}{\bold{u}}

\newcommand{\bv}{\bold{v}}

\newcommand{\bg}{\bold{g}}

\newcommand{\bh}{\bold{h}}

\newcommand{\bE}{\bold{E}}

\newcommand{\bJ}{\bold{J}}
\newcommand{\bB}{\bold{B}}
\newcommand{\bA}{\bold{A}}
\newcommand{\bR}{\bold{R}}

\newcommand{\bG}{\bold{G}}
\newcommand{\bQ}{\bold{Q}}
\newcommand{\bH}{\bold{H}}
\newcommand{\bZ}{\bold{Z}}
\newcommand{\bL}{\bold{L}}

\newcommand{\bY}{\bold{Y}}

\newcommand{\bD}{\bold{D}}

\newcommand{\bTheta}{\bold{\Theta}}

\newcommand{\bzeta}{\boldsymbol{\zeta}}

\newcommand{\xx}{\boldsymbol{x}}

\newcommand{\R}{\mathbb{R}}
\newcommand{\E}{\mathbb{E}}

\newcommand{\calR}{\mathcal{R}}

\newcommand{\calS}{\mathcal{S}}
\newcommand{\calT}{\mathcal{T}}

\newcommand{\argmin}{\mathop{\rm arg\min}}
\newcommand{\argmax}{\mathop{\rm arg\max}}

\newcommand{\vertiii}[1]{{\left\vert\kern-0.25ex\left\vert\kern-0.25ex\left\vert #1 
		\right\vert\kern-0.25ex\right\vert\kern-0.25ex\right\vert}}

\def\sf{{\cal F}}

\def\sc{{\cal C}}

\def\sp{{\cal P}}

\def\xx{\bold{x}}

\def\liminf{\mathop{\underline{\rm lim}}}

\makeatletter
\newcommand{\oset}[2]{%
	{\mathop{#2}\limits^{\vbox to -.5\ex@{\kern-\tw@\ex@
				\hbox{\scriptsize #1}\vss}}}}
\makeatother

\newenvironment{proof}[1][Proof]{\noindent\textbf{#1.} }{\ \rule{0.5em}{0.5em}}

\begin{document}
\title{Matrix Reordering for Noisy Disordered Matrices: Optimality and Computationally Efficient Algorithms}

\author{T. Tony Cai \ and \ Rong Ma
	\thanks{T. Tony Cai is with the Department of Statistics and Data Science at the University of Pennsylvania, Philadelphia, PA 19104 USA (email: tcai@wharton.upenn.edu).} 
	\thanks{Rong Ma is with the Department of Biostatistics at Harvard University, Boston, MA 02115 USA (email: rongma@hsph.harvard.edu).}
	\thanks{The research of Tony Cai was supported in part by NSF Grant DMS-2015259 and NIH grant R01-GM129781.}
		
		\thanks{Manuscript received September 19, 2022; revised August 11, 2023.}}


\date{}
 
\maketitle

\begin{abstract}
Motivated by applications in single-cell biology and metagenomics, we investigate the problem of matrix reordering based on a noisy disordered monotone Toeplitz matrix model. We establish the fundamental statistical limit for this problem in a decision-theoretic framework and demonstrate that a constrained least squares estimator achieves the optimal rate. However, due to its computational complexity, we analyze a popular polynomial-time algorithm, spectral seriation, and show that it is suboptimal. To address this, we propose a novel polynomial-time adaptive sorting algorithm with guaranteed performance improvement. Simulations and analyses of two real single-cell RNA sequencing datasets demonstrate the superiority of our algorithm over existing methods.
\end{abstract}

	\section{Introduction}

\begin{figure}[h!]
	\centering
	\includegraphics[angle=0,width=9cm]{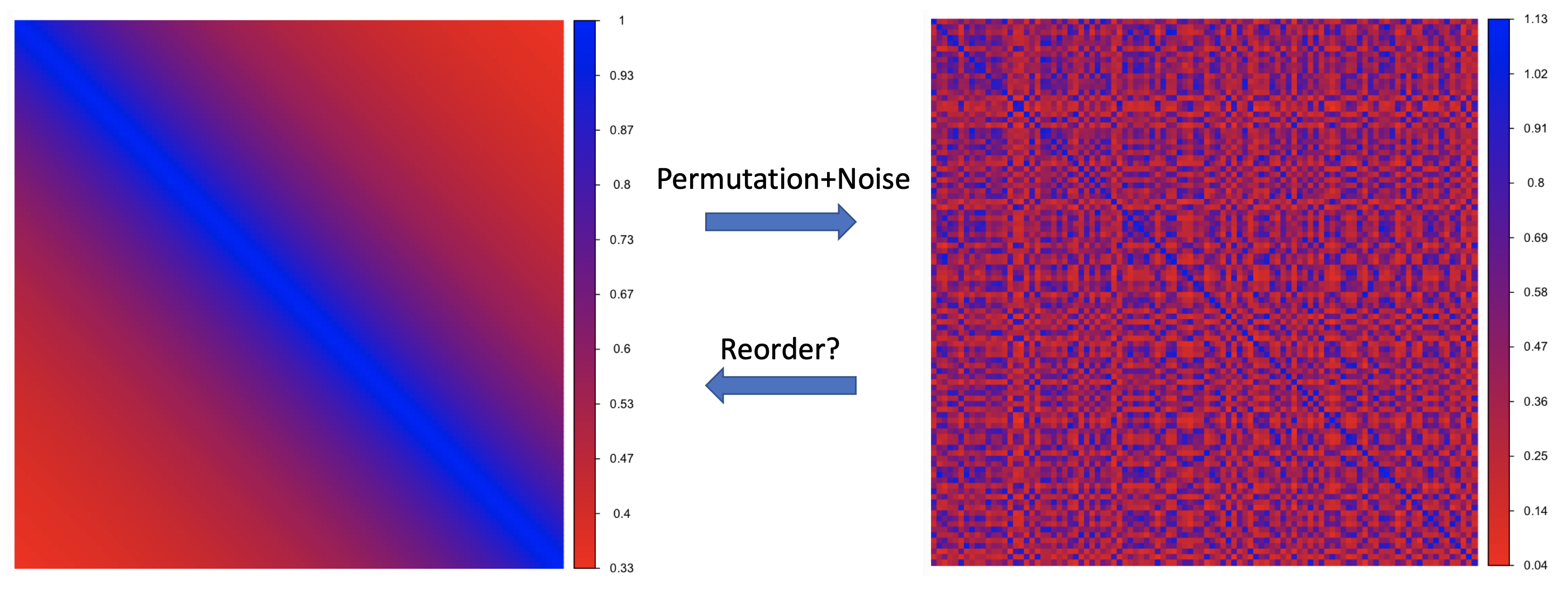}
	\caption{A graphical illustration of the matrix reordering problem. Left: a symmetric monotone Toeplitz matrix of dimension $100\times 100$ under Setting 6 of Section \ref{simu.sec}. Right: the observed noisy disordered matrix with an arbitrary permutation and the Gaussian noise.} 
	\label{Tmat}
\end{figure}

Consider the following noisy disordered   matrix model
\beq \label{R.model}
\bY =\Pi \bTheta \Pi^\top+\bZ,
\eeq
where $\bY,\bTheta,\bZ\in \R^{n\times n}$ are symmetric matrices, $\bY$ is observed, $\bZ$ is the noise matrix with independent (up to symmetry) sub-Gaussian entries with mean zero and variance $\sigma^2$,  $\bTheta$ is a deterministic signal matrix with certain structural patterns of interest, and $\Pi\in\R^{n\times n}$ is an unknown permutation matrix  that simultaneously permutes the columns and rows of the signal matrix $\bTheta$. This paper investigates the noisy matrix reordering problem, where the aim is to recover the underlying permutation $\Pi$ (Figure \ref{Tmat}) based on the observed noisy disordered matrix $\bY$.

This matrix reordering problem, also known as the matrix seriation problem, has a long history in data analysis and data mining \citep{liiv2010seriation,behrisch2016matrix}. It is often encountered when there is structural information contained in  the true signal matrix $\bTheta$, while only the corresponding noisy and disordered matrix $\bY=\Pi\bTheta\Pi^\top+\bZ$ is available. The importance of such a statistical problem lies in many applications (Section \ref{sec1.1}), where important features about the object under investigation are embedded in the structural pattern of  $\bTheta$, but are largely concealed from its noisy disordered observations.
For instance, when $\bTheta$ is a symmetric monotone Toeplitz matrix (Figure \ref{Tmat} left), after an arbitrary simultaneous permutation of its rows and columns, such a structural pattern is not easily discernible or recovered from a noisy   observation $\bY$ of the disordered matrix (Figure \ref{Tmat} right).

\subsection{Reordering Monotone Toeplitz Matrices and  Applications} \label{sec1.1}

The matrix reordering problem arises naturally in many applications \citep{liiv2010seriation,behrisch2016matrix,wu2008matrix}. 
The goal is to recover latent regularity and structural patterns contained in the noisy disordered data.  When viewed through the lens of model (\ref{R.model}), many of the applications 
involve reordering a noisy disordered matrix with some underlying monotone and Toeplitz structure, that is,
\beq  \label{T.layout}
\bTheta=\begin{bmatrix}
	\theta_0 & \theta_1 &\theta_2 &... &\theta_{n-1}\\
	\theta_1 & \theta_0 &\theta_1&...&\theta_{n-2}\\
	\theta_2 &\theta_1&\theta_0&...&\theta_{n-3}\\
	\vdots    &       &       & \ddots &\\
	\theta_{n-1} & \theta_{n-2} & \theta_{n-3}&...& \theta_0      
\end{bmatrix},
\eeq
where  $\theta_0\ge \theta_1\ge\theta_2\ge...\ge\theta_{n-1}$.
The following are two specific examples.

\begin{example}[Pseudotemporal ordering of single cells] \label{ss.ex}
	Single-cell analysis promises to revolutionize the treatment of common and rare diseases and provides insights into some of the most fundamental processes in biology. As an important problem in single-cell biology,
	pseudotemporal cell ordering aims to determine the pattern of a dynamic process experienced by cells and then arrange cells according to their progression through the process,
	based on single-cell RNA sequencing data collected at multiple time points. This problem can be formulated as a noisy matrix reordering problem (\ref{R.model}), where the entries of $\bTheta$ represent true cell-to-cell  similarities of $n$ cells in their transcriptomic profiles. For studies involving cells undergoing a dynamic process such as differentiation, monotonic patterns are often observed in the gene expression of the cells along the progression path \citep{trapnell2014dynamics,zeng2017pseudotemporal,saelens2019comparison}.
	In such cases,
	the true similarity matrix $\bTheta$ for the well-ordered cells can be modelled by a symmetric monotone Toeplitz matrix as in (\ref{T.layout}), and the goal is to recover the total order of the cells based on the noisy and disordered measurements  $\bY$ of  cell-to-cell similarities. For example, \cite{karin2022scprisma} considered a power decaying Toeplitz matrix with parameters $\theta_i=\alpha^i$ for some $\alpha\in(0,1)$, $i=1,2,...,n-1$,  to infer the underlying pseudotemporal ordering. See also Section \ref{data.sec} for more detailed discussion and the analyses of two real datasets.
\end{example}

\begin{example}[Genome assembly]
	In metagenomics and bioinformatics, genome assembly refers to the process of taking a large number of short DNA sequences and putting them back together in correct order to reconstruct the original sequence. In particular, the task of assembling a draft genome from shotgun metagenomic sequencing data can be treated as a noisy matrix reordering problem (\ref{R.model}), where each entry of $\bTheta$ characterizes the true genome distance between a pair of  contigs, or fragments of DNA sequence in the target genome.  Specifically, for a set of $n$ ordered contigs on one arm of the circular chromosome, their true pairwise genome distance matrix may be well approximated by a symmetric monotone Toeplitz matrix as in (\ref{T.layout}), and our goal is to recover the original genome order of these contigs based on the noisy and disordered measurements $\bY$ of their pairwise genome distances. Among existing works,  \citep{bagaria2020hidden} considered a graph-based model, which assumed a Hamiltonian path structure, or tridiagonal Toeplitz distance matrix (defined below) for the $n$ contigs on one arm; \citep{ma2021optimalb,ma2021optimala} considered a linear monotone model for the contigs on one arm, which  implies a linear decaying Toeplitz matrix (defined below) for the pairwise distance among these contigs. 
\end{example}


The matrix reordering problem also has important applications in combinatorial exploratory data analysis and data visualization \citep{wu2008matrix,behrisch2016matrix}. Specifically, for a given data matrix of interest, a proper reordering of its columns and rows may bring forth a more informative representation with structural patterns directly accessible or even visible to the analysts, thus providing critical guidance for downstream analysis \citep{friendly2002corrgrams,murdoch1996graphical}. { The importance of such a task in data analysis may be partially reflected by the current availability of a large variety of matrix visualization tools with an automatic matrix reordering option in standard statistical softwares, such as R. Among them, heatmap \citep{pheatmap} and corrplot \citep{corrplot2021} are probably the most commonly used tools.}

Motivated by these  interesting applications, this paper takes model (\ref{R.model}) as  a prototype underlying various matrix reordering problems  and focuses on reordering noisy symmetric matrices with latent monotone and Toeplitz structures. In particular, we consider the simple noise structure with independent sub-Gaussian entries up to symmetry.
Throughout, for any matrix $\bTheta$ admitting the expression (\ref{T.layout}), we call the entries corresponding to the value $\theta_0$ the main diagonals of $\bTheta$, and the entries corresponding to the value $\theta_i$ for $i\in\{1,...,n-1\}$ the $i$-th principal diagonals of $\bTheta$. 

In this paper, we consider the class of \emph{ridged monotone Toeplitz matrices} defined by
\beq \label{calT}
\calT_n = \left\{\bTheta\in\R^{n\times n}: \begin{aligned} &
	\text{$\bTheta$ admits the expression (\ref{T.layout})}\\
	& \theta_1\ge\theta_2\ge...\ge\theta_{n-1}\ge 0,\\
	& \theta_1-\theta_{\lceil n/2\rceil }\ge \theta_{\lceil n/2\rceil}-\theta_{n-1}
\end{aligned}\right\}.
\eeq
Note that the above definition puts no restriction on the main diagonal entries $\theta_0$ -- this is because for any $\bTheta\in\calT_n$, the main diagonal entries of $\Pi\bTheta\Pi^\top$ are invariant to the permutation $\Pi$, and therefore does not contain any information that helps for matrix reordering. The nonnegativity condition $\theta_i\ge 0$ and the direction of monotonicity are not essential here (Section \ref{dis.sec}).
The "ridge" condition $\theta_1-\theta_{\lceil n/2\rceil}\ge \theta_{\lceil n/2\rceil}-\theta_{n-1}$  ensures that the total amount of variations in the  first $\lceil n/2\rceil$ principal diagonals is no less than the variations in the rest of the diagonals. 
Such a characterization is required for technical reasons (Remark \ref{ridge}) but is in conformity with a wide range of applications. In particular, the class $\calT_n$  includes as special cases many interesting matrices that arise commonly in practice and have been  discussed in different contexts. 
As a few  examples, we note that $\calT_n$ includes, 
\begin{itemize}
	\item
	the tridiagonal Toeplitz matrices \citep{meurant1992review,noschese2013tridiagonal,da2020ninety} where $\theta_1>0$ and $\theta_j=0$ for all $j\ge 2$;
	\item
	the band monotone Toeplitz matrices \citep{bottcher2005spectral} where $\theta_1\ge\theta_2\ge...\ge\theta_k=...=\theta_{n-1}=0$ where $k<\lfloor n/2\rfloor$;
	\item
	the linear decaying Toeplitz matrices \citep{bunger2014inverses} where $\theta_j=\alpha+\beta (n-j)$ for all $1\le j\le n-1$ for some constants $\alpha,\beta>0$; and 
	\item
	the polynomial decaying Toeplitz matrices \citep{berenhaut2005monotone,cai2013optimal} where $\theta_j= Mj^{-\beta}$ for all $1\le j\le n-1$ for some constants $\beta,M\in(0,\infty)$. 
\end{itemize}

\subsection{Exact Matrix Reordering}

Throughout, we identify a permutation matrix $\Pi\in \R^{n\times n}$ with its corresponding permutation $\pi$, as an element in the symmetric group $\mathcal{S}_n$. Suppose $\Pi\in\calS_n$ is the underlying true permutation in $\Pi\bTheta\Pi^\top$ and let $\Pi'\in \calS_n$ be any given permutation. We quantify the distance between $\Pi'$ and $\Pi$ by the following 0-1 loss function 
\beq
\tau_{\bTheta}(\Pi,\Pi')=1\{\Pi \bTheta \Pi^\top\ne \Pi'\bTheta\Pi'^\top\}.
\eeq
Note that $\Pi'=\Pi$ is only a special case of $\tau_{\bTheta}(\Pi,\Pi')=0$. The loss $\tau_{\bTheta}(\Pi,\Pi')=0$ if and only if the two permutations produce the same disordered matrix, allowing for $\Pi\ne \Pi'$. The loss function $\tau_{\bTheta}(\Pi,\Pi')$ takes into account the possible equivalence classes among the permutations, caused by the specific structures of the signal matrix $\bTheta$. That is, we consider the exact permutation recovery modulo any inherent ambiguity caused by structure of the signal matrix. For example, when $\bTheta$ is Toeplitz with distinct diagonal elements, then the loss function will identify any two permutations up to a complete reversal. 

Let $\widehat\Pi$ be any estimator of $\Pi$ based on the observed matrix $\bY$. We define the estimation risk associated with the true parameters $\bTheta$ and $\Pi$ as 
\beq \label{risk}
\E_{\bTheta,\Pi} [\tau_{\bTheta}(\widehat\Pi,\Pi)]=P_{\bTheta,\Pi}(\Pi\bTheta\Pi^\top\ne \widehat\Pi \bTheta \widehat\Pi^\top),
\eeq
where the expectation on the left-hand side and the probability measure on the right-hand side are both with respect to the random observation $\bY$ for given $(\bTheta,\Pi)$. 

To evaluate the performance of an estimator, we consider the probability of exact matrix reordering  over a parameter space $\calT'_n\times \calS_n'=\{(\bTheta,\Pi):\bTheta\in\calT'_n, \Pi\in\calS'_n\}$ for some subsets  $\calT'_n\subseteq \calT_n$ and $\mathcal{S}'_n\subseteq \mathcal{S}_n$. In particular, we will identify sufficient and/or necessary conditions for parameter spaces of the form $\calT'_n\times \calS_n'$ such that, as $n\to\infty$, a given estimator could achieve exact matrix reordering uniformly over $\calT'_n\times \calS_n'$ with high probability. In this way, rigorous comparisons between various estimators can be made by comparing their respective conditions for exact matrix reordering. 
Specifically, for any $\calT'_n\times \calS_n'$, we define a hyper-parameter $\rho^*(\calT'_n,\mathcal{S}'_n)$ by
\beq
\rho^*(\calT'_n,\mathcal{S}'_n) 
=\inf_{\bTheta\in \calT'_n}\inf_{\substack{\Pi_1,\Pi_2\in\mathcal{S}'_n\\\Pi_1 \bTheta \Pi_1^\top\ne \Pi_2\bTheta\Pi_2^\top}}\|\Pi_1 \bTheta \Pi_1^\top-\Pi_2\bTheta\Pi_2^\top\|_F.
\eeq
Intuitively, the parameter $\rho^*(\calT'_n,\mathcal{S}'_n)$ quantifies how distinguishable two permuted versions of a matrix are within $\calT'_n\times \calS_n'$ -- a larger value of $\rho^*(\calT'_n,\mathcal{S}'_n)$ indicates a potentially bigger contrast between any two permuted matrices and vice versa. As will be shown shortly, this hyper-parameter reflects the overall signal strength, and therefore the fundamental difficulty of reordering  matrices over  a given parameter space.

\subsection{Main Results and Contributions}

The main results of this paper can be summarized as  follows.

\begin{enumerate}
	\item
	(Rate-optimal estimator)
	For any subsets $\calT'_n\subseteq \calT_n$ and $\mathcal{S}'_n\subseteq \mathcal{S}_n$ such that $\rho^*(\calT'_n,\mathcal{S}'_n)\gtrsim\sigma\sqrt{n\log n}$,  there exists a constrained least square estimator  (LSE)  that achieves exact matrix reordering for any $(\bTheta,\Pi)\in \calT'_n\times\mathcal{S}'_n$ with high probability.  \label{pp1}
	\item
	(Fundamental information threshold) 
	There exist some subsets $\calT'_n\subset \calT_n$ and $\mathcal{S}'_n\subset\mathcal{S}_n$ satisfying  $\rho^*(\calT'_n,\mathcal{S}'_n)\asymp\sigma\sqrt{n\log n}$, such that no permutation estimator could achieve exact matrix reordering for all $(\bTheta,\Pi)\in \calT'_n\times\mathcal{S}'_n$ with high probability. \label{pp2}
	\item(Suboptimality of the spectral seriation algorithm) 
	There exist some subsets $\calT'_n\subset \calT_n$ and $\mathcal{S}'_n\subset\mathcal{S}_n$ satisfying  $\rho^*(\calT'_n,\mathcal{S}'_n)\asymp\sigma n^3$, such that the commonly used spectral seriation algorithm cannot achieve exact matrix reordering for any $(\bTheta,\Pi)\in \calT'_n\times\mathcal{S}'_n$  with high probability. \label{pp4}
	\item (An improved polynomial-time algorithm) 
	For any subsets $\calT'_n\subseteq \calT_n$ and $\mathcal{S}'_n\subseteq \mathcal{S}_n$ such that $\rho^*(\calT'_n,\mathcal{S}'_n)\gtrsim \sigma n^2$,  there exists a polynomial-time algorithm, proposed in Section \ref{ipdd.sec}, that achieves exact matrix reordering for any $(\bTheta,\Pi)\in \calT'_n\times\mathcal{S}'_n$ with high probability.  \label{pp3}
\end{enumerate}

\begin{figure}
	\centering
	\includegraphics[angle=0,width=9cm]{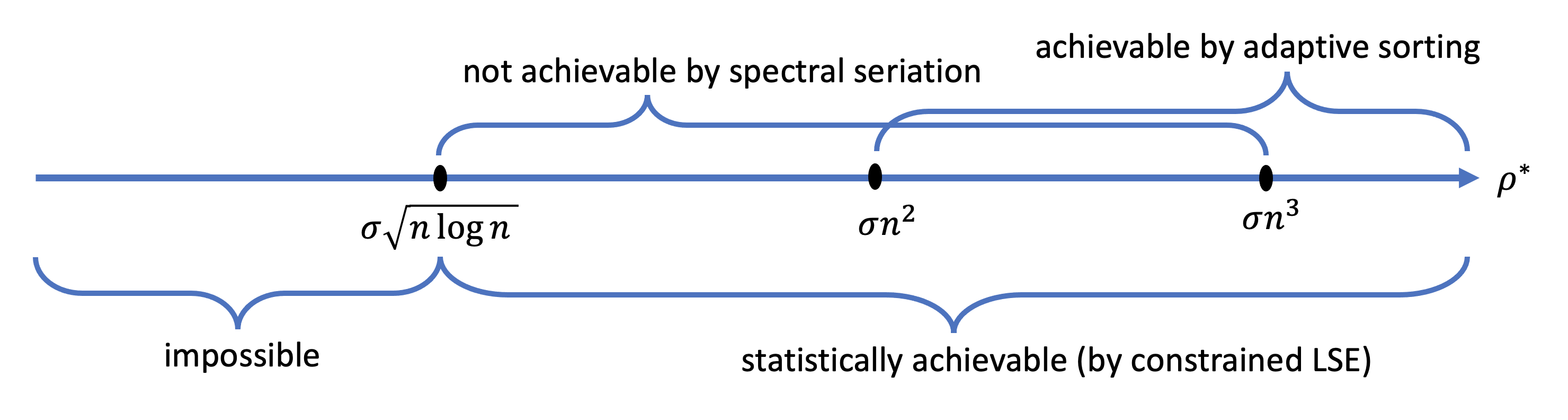}
	\caption{\small A graphical illustration of the main theoretical results. The rates under the horizontal line are thresholds for $\rho^*(\calT'_n,\mathcal{S}'_n)$ for three different matrix reordering algorithms, constrained LSE,  spectral seriation, and adaptive sorting, which are defined in Sections \ref{mle.sec},  \ref{ss.sec}, and \ref{ipdd.sec}, respectively.} 
	\label{rec}
\end{figure}

The results are  illustrated in Figure \ref{rec}.  Parts  \ref{pp1}) and \ref{pp2}) of the main results together show that a phase transition occurs at order $\sigma\sqrt{n\log n}$: no method can achieve the exact order recovery with high probability if the signal strength (as measured by $\rho^*$)  is at or below this level and the constrained LSE recovers the order exactly with high probability whenever the signal strength is above this level.

Although the constrained LSE is rate-optimal, it requires solving an optimization  over a potentially large permutation set $\calS'_n$ and is in general computationally infeasible for large $n$. Alternatively, polynomial-time algorithms have been used in practice. Among them, spectral seriation \citep{barnard1995spectral,atkins1998spectral,fogel2014serialrank} is arguably  the current state-of-art  polynomial-time matrix reordering algorithm\footnote{See Section \ref{simu.sec} for empirical evidences of the superiority of the spectral seriation over other existing methods.}.  We analyze the performance of the spectral seriation algorithm and show that it is suboptimal for reordering Toeplitz matrices as in $\mathcal{T}_n$.

We then develop a novel adaptive sorting algorithm, which runs in polynomial time, and show that it has strictly weaker signal strength requirement compared to the  spectral seriation algorithm. This result indicates the advantage of the proposed adaptive sorting algorithm over the current state-of-art matrix reordering method, and explains its overall superior empirical performance over various existing methods across a wide range of simulation settings and real data examples. On the other hand, despite its numerical advantages, the adaptive sorting algorithm is still statistically suboptimal. We conjecture in Section \ref{tradeoff.sec} that there is a fundamental gap between statistical optimality and computational efficiency.


\subsection{Related Works} \label{related.sec}

Many statistical seriation problems
 that in one way or another aim to find an element in the discrete permutation set optimizing certain objective function have been  studied from various aspects under different settings. These include the well-known consecutive one's problem \citep{fulkerson1964incidence,kendall1963statistical,kendall1970mathematical} that dates back to the 1960s;  the feature matching problem \citep{collier2016minimax,jeong2020recovering,galstyan2021optimal} and the noisy ranking problem \citep{braverman2009sorting,jiang2011statistical,negahban2012iterative,chatterjee2016estimation,shah2016stochastically,mao2017minimax,negahban2018learning,chen2019spectral}; the matrix seriation problem for various shape-constrained matrices including the monotone or bi-monotone matrices \citep{flammarion2019optimal,mao2020towards,ma2021optimalb,pananjady2020isotonic},  the Robinson matrices \citep{atkins1998spectral,fogel2014serialrank,recanati2018reconstructing,armstrong2021optimal},  and the Monge matrices \citep{hutter2020estimation};  and more recently, the seriation problem under the latent space models \citep{giraud2021localization,issartel2021optimal}.  

Many of the existing works have focused on recovering the underlying permutations, estimation of the (disordered) signal structures, or both. However, about the matrix seriation problems,  statistical limit and  optimal procedures for the permutation recovery problem are relatively less-understood, compared to the estimation of the  signal matrices \citep{flammarion2019optimal,mao2020towards,hutter2020estimation,ma2021optimalb}. 
In particular, Bagaria et al. \cite{bagaria2020hidden} considered the hidden Hamiltonian cycle recovery model, which is related to our model (\ref{R.model}) under a tridiagonal signal matrix, and Ding et al. \cite{ding2021consistent} considered a hidden $2k$-nearest neighbor graph recovery model, corresponding to our model (\ref{R.model}) under a band signal matrix. Both papers established the information-theoretical threshold for exact recovery and/or almost exact recovery, but under some related but in general different signal-to-noise ratio parameter. Moreover, these papers considered exact recovery over a special set of signal matrices  (say, tridiagonal Toeplitz matrices) and the full permutation set $\calS_n$ alone, whereas our work considers exact recovery over all possible combinations of subsets of signal matrices and permutations, which may reveal fundamentally more difficult scenarios. Besides, for the general monotone Toeplitz matrix reordering problem as formulated in the current paper,  most of the existing results  only concern permutation recovery in the noiseless setting \citep{atkins1998spectral,fogel2014serialrank,recanati2018reconstructing}, thus leaving the fundamental behavior of the problem under  the more realistic noisy observations largely unexplored.

\subsection{Organization}

The rest of the paper is organized as follows. We finish this section with notation that will be used throughout the paper.  Section \ref{rate.sec} establishes the fundamental statistical limit for exact matrix reordering, including the minimal signal strength required by the constrained LSE, and the matching fundamental information threshold. Section \ref{ss.sec} analyzes the spectral seriation algorithm and shows its fundamental suboptimality. In Section \ref{ipdd.sec}, the adaptive sorting algorithm is proposed and its theoretical properties are investigated. Section \ref{tradeoff.sec} discusses the potential tradeoff between computational efficiency and statistical optimality.  Section \ref{simu.sec} contains simulation studies that compare the empirical performances of several matrix reordering algorithms in various settings.  Section \ref{data.sec} presents the analyses of two real datasets, showing the advantage of the adaptive sorting in real-world applications.  Possible extensions  of the current work are discussed in Section \ref{dis.sec}. The proofs of main results are given in Section \ref{proof.sec} and the proofs of other technical results are given in \ref{supp.sec}.

\subsection{Notation}
For a vector $\bold{a} = (a_1,...,a_n)^\top \in \mathbb{R}^{n}$,  $\text{diag}(a_1,...,a_n)\in\R^{n\times n}$ denotes the diagonal matrix whose $i$-th diagonal entry is $a_i$, and define the $\ell_p$ norm $\| \bold{a} \|_p = \big(\sum_{i=1}^n |a_i|^p\big)^{1/p}$.  
For a matrix $ \bold{A}=(a_{ij})\in \R^{n\times n}$,  its Frobenius norm is $\| \bold{A}\|_F = \sqrt{ \sum_{i=1}^{n}\sum_{j=1}^{n} a^2_{ij}}$ 
;  its $i$-th column is denoted by $ \bold{A}_{.i}\in \R^{n}$ and its $i$-th row by $ \bold{A}_{i.}\in \R^{n}$. Moreover, we denote $\bA_{i,-j}$ as a subvector of $\bA_{i.}$ with its $j$-th component removed.
For any integer $n>0$, we denote the set $[n]=\{1,2,...,n\}$. 
For any $a>0$, $\lfloor a\rfloor$ denotes the largest integer no greater than $a$, and $\lceil a \rceil$ denotes the smallest integer no less than $a$.
For a finite set $S$,  its cardinality is denoted by $|S|$.  A random variable $X$ is sub-Gaussian if there are positive constants $C,v$ such that for every $t>0$, we have $P(|X|>t)\le Ce^{-vt^2}$.
For sequences $\{a_n\}$ and $\{b_n\}$, we write $a_n = o(b_n)$ or $a_n\ll b_n$ if $\lim_{n} a_n/b_n =0$, and write $a_n = O(b_n)$, $a_n\lesssim b_n$ or $b_n \gtrsim a_n$ if there exists a constant $C$ such that $a_n \le Cb_n$ for all $n$. We write $a_n\asymp b_n$ if $a_n \lesssim b_n$ and $a_n\gtrsim b_n$. Throughout, $C,C_1,C_2,...$ are universal constants independent of $n$, and can vary from place to place.

\section{Fundamental Statistical Limit for Matrix Reordering} \label{rate.sec}

Our main result on the statistical limit for matrix reordering consists of two parts: a fundamental information threshold that benchmarks all the matrix reordering algorithms, and a constrained least square permutation estimator, whose performance is rate-optimal among all the estimators. We start with the rate-optimal estimator.

\subsection{The Constrained Least Square Estimator} \label{mle.sec}

For any given parameter space $\calT'_n\times\mathcal{S}'_n$ where $\calT'_n\subseteq \calT_n$ and $\calS'_n\subseteq\calS_n$, suppose one observes $\bY$ from (\ref{R.model}) for some $(\bTheta^*,\Pi^*)\in \calT'_n\times \calS'_n$. A natural estimator for the unknown permutation is the constrained least square estimator (or the maximum likelihood estimator in the $i.i.d.$ Gaussian case) over $\calT'_n\times \calS'_n$ defined through
\beq \label{mle}
(\widehat{\bTheta}^{lse}, \widehat{\Pi}^{lse})=\argmin_{(\bTheta,\Pi)\in \calT'_n\times \mathcal{S}'_{n}} \|\bY-\Pi\bTheta\Pi^\top\|^2_F.
\eeq
The following result provides the theoretical guarantee of $\widehat\Pi^{lse}$ over $\calT'_n\times\mathcal{S}'_n$. 

\bet[Theoretical guarantee for constrained LSE] \label{mle.thm}
Under model (\ref{R.model}),  there exists some absolute constants $C,c>0$ such that,  for sufficiently large $n$, for any $\calT'_n\subseteq \calT_n$ and any $\mathcal{S}_n'\subseteq \mathcal{S}_n$ such that $\rho^*(\calT'_n,\mathcal{S}'_n)\ge C\sigma\sqrt{n\log n}$, the permutation estimator $\widehat{\Pi}^{lse}$ given by (\ref{mle}) satisfies
\beq
\sup_{(\bTheta,\Pi)\in \calT'_n\times \mathcal{S}'_n}P_{\bTheta,\Pi}(\widehat\Pi^{lse}\bTheta (\widehat\Pi^{lse})^\top\ne \Pi\bTheta\Pi^\top)\le n^{-c}.
\eeq
\eet

Theorem \ref{mle.thm} characterizes the explicit correspondence between the signal strength condition on $\rho^*$ and the final exact recovery error probability; it applies to any subsets of ridged monotone Toeplitz matrices, any subsets of permutations, and the general sub-Gaussian noises. The theorem identifies a sufficient minimal signal strength condition
\beq \label{suf.cond}
\rho^*(\calT'_n,\calS'_n)\gtrsim\sigma\sqrt{n\log n},
\eeq 
under which the constrained LSE $\widehat{\Pi}^{lse}$ achieves exact matrix reordering uniformly over $(\bTheta,\Pi)\in\calT'_n\times \calS'_n$ with high probability. Importantly,  combined with the information lower bound obtained in Section \ref{lbnd.sec}, Theorem \ref{mle.thm} essentially implies that $\widehat\Pi^{lse}$ is rate-optimal. 

The proof of Theorem \ref{mle.thm}, detailed in Section \ref{mle.thm.sec}, relies on analyzing the probability of exact matrix reordering for the constrained LSE at any given parameters $(\bTheta,\Pi)\in \calT'_n\times \mathcal{S}'_n$. To do so, we develop a general reduction scheme that connects the risk of matrix reordering to the risk of matrix denoising, i.e., estimating the permuted matrix $\Pi\bTheta\Pi^\top$ from $\bY$, under model (\ref{R.model}). We summarize our reduction scheme as the following proposition, proved in Section \ref{prop1.sec}.

\bep[Reduction scheme] \label{red.prop}
For any $\bTheta, \widehat\bTheta\in\calT_n$ and any permutations  $\Pi,\widehat\Pi\in \calS'_n\subseteq\calS_n$, it holds that
\beq \label{p.red}
P(\widehat\Pi \bTheta \widehat\Pi^\top\ne \Pi\bTheta\Pi^\top)\le e^{{2\E \|\widehat{\Pi}\widehat{\bTheta} \widehat{\Pi}^\top-\Pi\bTheta\Pi^\top\|_F}-{\rho(\bTheta,\mathcal{S}'_n)}},
\eeq
where 
\[
\rho(\bTheta,\mathcal{S}'_n)=\min_{\substack{\Pi_1,\Pi_2\in\mathcal{S}'_n\\\Pi_1 \bTheta \Pi_1^\top\ne \Pi_2\bTheta\Pi_2^\top}}\|\Pi_1 \bTheta \Pi_1^\top-\Pi_2\bTheta\Pi_2^\top\|_F.
\]
\eep

Proposition \ref{red.prop} provides a tool for bounding the probability $P_{\bTheta,\Pi}(\widehat\Pi^{lse}\bTheta (\widehat\Pi^{lse})^\top\ne \Pi\bTheta\Pi^\top)$ by analyzing the matrix denoising risk 
\beq \label{mle.risk}
\E \|\widehat{\Pi}^{lse}\widehat{\bTheta}^{lse}( \widehat{\Pi}^{lse})^\top-\Pi\bTheta\Pi^\top\|_F,
\eeq
which is easier to handle. This reduction step paves the way for an asymptotically sharp risk analysis of $(\widehat{\bTheta}^{lse}, \widehat{\Pi}^{lse})$ using powerful tools developed for general shape-constrained least square estimators. The key ingredients of our proof, which generalizes the idea for proving Theorem 3.1 of \citep{flammarion2019optimal}, include Chatterjee's variational formula (Lemma \ref{var.lem}), an improved Dudley's integral inequality (Lemma \ref{dudley.lem}), and a nontrivial calculation of the metric entropy of a set of permuted Toeplitz matrices (Lemma \ref{entropy.lem}).

\subsection{Fundamental Information Threshold and the Planted Path Reconstruction Problem} \label{lbnd.sec}

We investigate the necessity of the minimal signal strength condition (\ref{suf.cond}), and uncover the fundamental information threshold underlying the matrix reordering problem.
The following theorem provides a lower bound for the minimum signal  strength.

\bet[Fundamental information threshold] \label{low.bnd.thm}
Suppose $n\ge 48$. Then there exist subsets $\calT'_n\subset \calT_n$ and $\mathcal{S}'_n\subset\mathcal{S}_n$ satisfying $\rho^*(\calT'_n,\mathcal{S}'_n)= 0.02\sigma\sqrt{n\log n}$ such that 
\beq
\inf_{\widehat\Pi}\sup_{(\Theta,\Pi)\in \calT'_n\times \mathcal{S}'_n}P_{\bTheta,\Pi}(\widehat\Pi\bTheta \widehat\Pi^\top\ne \Pi\bTheta\Pi^\top) \ge 0.6.
\eeq
\eet

Theorem \ref{low.bnd.thm} shows that, under the Gaussian noise there exists a certain parameter space $\calT'_n\times \calS'_n$ with minimal signal strength $ \rho^*(\calT'_n,\calS'_n)\asymp \sigma\sqrt{n\log n}$ such that no method could achieve exact matrix reordering uniformly with high probability.  Theorems \ref{mle.thm} and \ref{low.bnd.thm} together show that the condition (\ref{suf.cond}) is asymptotically sharp and the rate
\beq \label{bnd}
\rho^*(\calT'_n,\calS'_n)\asymp \sigma\sqrt{n\log n},
\eeq
is the fundamental information threshold for the matrix reordering problem and the permutation estimator $\widehat\Pi^{lse}$ given by (\ref{mle}) is minimax rate-optimal.

To prove the information lower bound in Theorem \ref{low.bnd.thm}, we connect the matrix reordering problem with the following planted path reconstruction problem in graph-information theory, and obtain the information lower bound of the former problem by analyzing that of the latter problem.

\begin{definition}[Planted Hamiltonian path reconstruction]
	Consider a weighted undirected graph $G(V,E)$ with an adjacecy matrix $\bZ$, where $\bZ$ has $i.i.d.$ standard normal entries up to symmetry. Suppose an arbitrary Hamiltonian path $P(\theta)$ of constant edge weight $\theta$ and length $|V|-1$ connecting all the vertices in $V$, is added to $G(V,E)$, resulting to a new graph $G'(V,E')$. Then we refer problem of reconstructing $P$, or equivalently the recovery of the support of its adjacency matrix,  from $G'(V,E')$, as the planted Hamiltonian path reconstruction problem.
\end{definition}

To study the statistical limit of the above path reconstruction problem, we construct a least favorable class $\mathcal{H}$ of Hamiltonian paths, each with a constant edge weight $\theta>0$, such that (i) the class is sufficiently large in the sense that $\log |\mathcal{H}| \gtrsim n\log n$; and (ii) any two paths in $\mathcal{H}$ are sufficiently distinct from each other under the Hamming distance (defined in Lemma \ref{deza}). 
In particular, we show that whenever $\theta\lesssim \sqrt{\log n}$, there is no way to tell  with confidence from the new graph  $G'(V,E')$ which path in $\mathcal{H}$ is planted in $G'(V,E')$.
To construct the set $\mathcal{H}$, we introduce to the current context a useful result due to \cite{deza1976matrices}
concerning the Hamming packing in the permutation space \citep{cameron2005covering,quistorff2006survey,hendrey2020covering}. 

\bel[Deza's bound on permutation packing] \label{deza}
For any $\pi_1,\pi_2\in\mathcal{S}_n$, we define their Hamming distance  $d_H(\pi_1,\pi_2)=|\{i:\pi_1(i)\ne \pi_2(i)\}|$. A $d$-packing in the finite metric space $(\mathcal{S}_n,d_H)$ is a subset $M\subset \mathcal{S}_n$ such that its elements are at a distance of at least $d$ from each other. Then the largest cardinality of a $d$-packing $\beta_n(d)$ satisfies
\beq
\beta_n(d) \ge \frac{n!}{V_d},
\eeq
where 
\beq
V_d=\sum_{k=0}^d{n\choose k} k!\sum_{x=0}^k \frac{(-1)^x}{x!}.
\eeq
Consequently, for $2\le d\le n-1$, we have
\beq
\beta_n(d) \ge \frac{n!}{(n-d)!}\cdot \frac{(n-d+1)}{2(n-d)}.
\eeq
\eel

The detailed  proof of Theorem \ref{low.bnd.thm} is provided in Section \ref{low.bnd.thm.sec}. The proof techniques developed there may be applied to obtain information lower bounds for other permutation related problems.

In this study, we define the fundamental information threshold to be the minimum separation $\rho^*(\calT'_n, \calS'_n)$ required for {any possible} subset $\calT'_n$ of ridged monotone matrices, and {any subset} $\calS'_n$ of permutations, in order for the exact recovery to be achievable. This is different from the existing work such as \citep{bagaria2020hidden} and \citep{ding2021consistent} where the separation is considered for a specific subset of signal matrices such as tridiagonal Toeplitz matrices, and the {full} set $\calS_n$ of the permutations. In particular, from our proof of Theorem \ref{low.bnd.thm}, one can show that when restricted to the settings of \citep{bagaria2020hidden}, a smaller  separation condition $\rho^*\gtrsim\sigma\sqrt{\log n}$ may be obtained. 
The stronger requirement on the minimum separation obtained by Theorem \ref{low.bnd.thm}, is essentially due to the greater variety of scenarios allowed by our framework, some of which can be more difficult than those considered in the existing work. In this sense, our specific construction mainly reveals the fundamental difficulty caused by the unknown underlying permutation set. Finally, we remark that due to the nature of our minimax lower bound argument, it remains unclear if there exist many other hard scenarios, concerning possibly different signal matrices, under which a $\sigma\sqrt{n\log n}$-order separation is needed. Nevertheless, we believe this is an important problem that deserves further investigation.

\begin{remark} \label{remark}
Our proof of Theorems 1 and 2 does not involve the ``ridge" condition $\theta_1-\theta_{n/2}\ge \theta_{n/2}-\theta_{n-1}$. However, such a condition plays an important role in our analysis of the proposed adaptive sorting algorithm (Section \ref{ipdd.sec}), and is satisfied by our  suboptimality argument of the spectral seriation algorithm (Section \ref{ss.sec}). Therefore, we include it in our definition of the parameter space for the integrity of our theoretical statements.  
\end{remark}

\section{Suboptimality of Spectral Seriation} \label{ss.sec}

The constrained LSE introduced in Section \ref{mle.sec} is rate-optimal. However, obtaining such an estimator requires solving an optimization over the discrete permutation set $\calS'_n\subseteq \calS_n$, which could be computationally infeasible, either because the subset $\calS'_n$ is unknown, or $\calS'_n$ contains a large number of permutations that grows exponentially fast as the matrix size $n$ increases. 

Alternatively, spectral approaches have been widely used  for matrix reordering tasks.  
Among them, a spectral seriation method based on the Fiedler vector is particularly popular and has been extensively studied in the literature \citep{barnard1995spectral,atkins1998spectral,fogel2014serialrank,recanati2017spectral,recanati2018reconstructing,giraud2021localization}. In this section, we show that, despite the success of such a spectral seriation algorithm in many applications, it is nonetheless suboptimal for reordering Toeplitz matrices compared to the constrained LSE.

To formally introduce the spectral seriation estimator $\check\Pi$, 
we define the following ranking function.

\bed[Ranking function] \label{r.def}
The ranking operator $\frak{r}: \R^n \to \mathcal{S}_n$  is defined such that for any vector $\xx\in \R^n$, $\frak{r}(\xx)$ contains the ranks of the components of $\xx$ in increasing order.  Whenever there are ties, increasing orders are assigned from left to right.
\eed

As an example, for a vector $\xx = (2,5,1,6,2)^\top$, we have $\frak{r}(\xx) = (2,4,1,5,3)$. Following \cite{atkins1998spectral,fogel2014serialrank,recanati2018reconstructing}, the spectral seriation estimator $\check\Pi$ is then defined in Algorithm \ref{al1} below.

\begin{algorithm}
	\caption{Spectral seriation} \label{al1}
	\begin{algorithmic}
		\State {\bf Input:} Observed matrix $\bY=(Y_{ij})_{1\le i,j\le n}\in\R^{n\times n}$. 
		\State \hspace{4mm} 1. Compute the Laplacian matrix $\bL = \bD-\bY$ where $\bD=\text{diag}(d_1,...,d_n)$ and $d_i=\sum_{j=1}^n Y_{ij}$.
		\State \hspace{4mm}  2. Obtain the Fiedler eigenvector $\widehat\bv\in\R^n$ corresponding to the smallest nonzero eigenvalue of $\bL$.
		\State {\bf Output:} $\check{\Pi}=[\frak{r}(\widehat\bv)]^{-1}$, where 
		the inverse $[\cdot]^{-1}$ means the reversion of a permutation.
	\end{algorithmic}
\end{algorithm}

For matrix reordering, 	Algorithm \ref{al1} has been shown in \citep{atkins1998spectral} to achieve exact recovery in the noiseless case $\bZ=0$ for all Robinson matrices $\bR=(r_{ij})_{1\le i,j\le n}$ satisfying
\beq \label{Rob}
\bR=\bR^\top, \qquad r_{ij}\le \min\{r_{ik}, r_{kj}\}\text{ for all $1\le i<k<j\le n$},
\eeq
which contain Toeplitz matrices in $\mathcal{T}_n$ as a special case. This result is summarized in the following proposition.


\bep \label{ss.prop}
For any  $\bR\in\R^{n\times n}$ satisfying  (\ref{Rob}) and $\Pi\in\calS_n$, suppose the Fiedler eigenvector $\bv\in\R^n$ associated to the  smallest nonzero eigenvalue of $\Pi\bR\Pi^\top$ contains $n$ distinct components, then we have $\Pi=[\frak{r}(\bv)]^{-1}.$
\eep

Proposition \ref{ss.prop} essentially implies  that,  in the noiseless setting, the spectral seriation algorithm is able to achieve exact reorder recovery of any matrix $\bTheta\in\calT_n$ based on $\Pi\bTheta\Pi^\top$. 

However, the story is different in the noisy settings. Our analysis shows that the spectral seriation can be sensitive to the eigen-structure of the signal matrices and the noises. As a consequence, it may suffer from inconsistent estimation, and therefore significant suboptimality, due to insufficient separation between the Laplacian eigenvalues. Similar results concerning the suboptimality of spectral methods have been obtained in \cite{bagaria2020hidden} under a hidden Hamiltonian cycle recovery model, but with a slightly different formulation and signal strength measure.

\bet \label{ss.thm}
Suppose the noise matrix has $i.i.d.$ entries up to symmetry generated from $N(0,\sigma^2)$.  Then there exists some $\calT'_n\subseteq \calT_n$ and  $\mathcal{S}_n'\subseteq \mathcal{S}_n$ with $\log |\calS'_n|\gtrsim n\log n$ satisfying $\rho^*(\calT'_n,\mathcal{S}'_n)=C\sigma n^3$ for some absolute constant $C>0$, such that $$\liminf_{n\to \infty}\inf_{(\bTheta, \Pi)\in \calT'_n\times \mathcal{S}'_n}P_{\Theta,\Pi}(\check{\Pi}\bTheta\check\Pi\ne \Pi\bTheta\Pi^\top)\ge1/2.$$
\eet

Theorem \ref{ss.thm} shows that even when the minimal signal strength $\rho^*$ is of order $\sigma n^3$, 
there still exists nontrivial cases over which the exact matrix reordering using the spectral seriation is impossible. 
The proof of Theorem \ref{ss.thm}, given in Section \ref{ss.thm.sec} relies on a delicate eigenvector  analysis of  a  deformed random Laplacian matrix. To this end, we develop a novel triangulation argument inspired by \cite{lu2002} and \cite{johnstone2009sparse} that allows us to show  inconsistency of the sample Fiedler vector $\widehat\bv$ in relation to the underlying true Fiedler vector, in the so-called subcritical regime \citep{bloemendal2016principal,bao2021singular,cai2021optimal}. Our analytic framework can be useful for other lower bound problems in statistics and random matrix theory, especially when they involve  characterizing the asymptotic behavior of Laplacian eigenvectors  associated to the bulk eigenvalues.

\section{Efficient Matrix Reordering via Adaptive Sorting} \label{ipdd.sec}
The suboptimality of  the spectral seriation  motivates us to develop an alternative algorithm with improved performance. We propose in this section a novel polynomial-time matrix reordering algorithm. The method, referred as the adaptive sorting,  is summarized below in Algorithm \ref{al0}.

\begin{algorithm}
	\caption{Adaptive sorting} \label{al0}
	\begin{algorithmic}
		\State {\bf Input:} Observed matrix $\bY=(Y_{ij})\in \R^{n\times n}$. 
		\State 1. {\bf Locating the initial element.} For given matrix $\bY$,
		\State \hspace{4mm} (i) calculate  $S_i=\sum_{j\in[n]\setminus\{i\}}Y_{ij}$  for $i\in[n]$;
		\State \hspace{4mm}  (ii) set $\widetilde\pi(1)=\argmin_{i\in[n]}S_i$.
		\State 2. {\bf Iterative sorting.} For $i=1,2,...,n-1$, set
		$$\widetilde\pi(i+1) = \argmin_{j\in [n]\setminus \{\widetilde\pi(1),...,\widetilde\pi(i)\}} \|\bY_{\tilde{\pi}(i),-\tilde{\pi}(i)}-\bY_{j,-j}\|_1, $$ where $\bY_{i,-i}\in\R^{n-1}$ is the $i$-th row of $\bY$ with $i$-th component removed.
		\State {\bf Output:}  $\widetilde\pi=(\widetilde\pi(1),...,\widetilde\pi(n))$, and $\widetilde\Pi$ as the corresponding permutation matrix.
	\end{algorithmic}
\end{algorithm}	

Step 1 of Algorithm \ref{al0} identifies the location of the first (or equivalently, the last) row of the signal matrix after permutation. It uses the fact that the sum  of the first or the last row of the original signal matrix $\bTheta$ is the smallest among all the row sums. Once a beginning point of $n$ elements in their original order has been identified, Step 2 starts from there and builds up the complete permutation map $\tilde{\pi}$ iteratively based on the following rationale: after removing the main diagonals, each row of a monotone Toeplitz matrix is more similar to its nearby rows than the more distant rows. In particular, it can be shown that any reduced row (main diagonal removed) of the correctly ordered signal matrix, say $\bTheta_{i,-i}$, has its nearest $\ell_1$-neighbors the rows $\bTheta_{j,-j}$ for $j\in\{i-1,i+1\}$.
Thus, for given $\widetilde\pi(i)$, to identify $\widetilde\pi(i+1)$,  we look  among the remaining rows  $\{\bY_{j,-j}\}_{j\in[n]\setminus \{\widetilde\pi(1),...,\widetilde\pi(i)\}}$ for the one that minimizes  $\|\bY_{\tilde{\pi}(i),-\tilde{\pi}(i)}-\bY_{j,-j}\|_1$. This is applied iteratively until all the rows (or columns) are properly ordered.

	Related to the above reasoning,  the application of the $\ell_1$-norm for comparing the reduced rows $\{\bY_{i,-i}\}_{1\le i\le n}$ in Step 2 of Algorithm \ref{al0} is  rooted in the bias-variance tradeoff. For example, 
	although the nearest neighbor of a given reduced row $\bTheta_{i,-i}$ can be invariably determined based on either the $\ell_1$- or the $\ell_2$-norm,
	the variance of the observed $\ell_2$-distance $\|\bY_{\tilde{\pi}(i),-\tilde{\pi}(i)}-\bY_{j,-j}\|_2$ can be much larger than the variance of $\|\bY_{\tilde{\pi}(i),-\tilde{\pi}(i)}-\bY_{j,-j}\|_1$. From our analysis in Section \ref{eff.thm.sec}, it can be seen that an inflated variability may significantly deteriorate the performance of the sorting algorithm.

The following theorem provides the theoretical guarantee for the adaptive sorting algorithm, whose proof can be found in Section \ref{eff.thm.sec}.

\bet[Theoretical guarantee for adaptive sorting] \label{eff.thm}
Under model (\ref{R.model}), there exists some absolute constants $C,c>0$ such that,  for sufficiently large $n$, for any $\calT'_n\subseteq \calT_n$ and any $\mathcal{S}_n'\subseteq \mathcal{S}_n$ such that $\rho^*(\calT'_n,\mathcal{S}'_n)\ge C\sigma n^2$, we have $$\sup_{(\bTheta, \Pi)\in \calT'_n\times \mathcal{S}'_n}P_{\Theta,\Pi}(\widetilde\Pi\bTheta \widetilde\Pi^\top\ne \Pi\bTheta\Pi^\top)\le n^{-c}.$$
\eet

Compared with the results from the previous sections, although the adaptive sorting algorithm requires a condition
\beq \label{s.c}
\rho^*(\calT'_n,\mathcal{S}'_n)\gtrsim \sigma n^2,
\eeq
that is till stronger than the optimal condition (\ref{suf.cond}), it is nevertheless much weaker than that required by the spectral seriation, by  a factor at least of order $n$. In general, the adaptive sorting algorithm takes advantage of the Toeplitz structure to achieve better performance in reordering matrices in $\mathcal{T}_n$, whereas the spectral seriation may have wider applicability when signal-to-noise ratio is sufficiently large.
As for the fundamental limit of the adaptive sorting, in Section \ref{lb.sec} below, we show that there exist $\calT'\subset \calT_n$ and $\calS'_n\subset \calS_n$ satisfying  $\rho^*(\calT'_n,\mathcal{S}'_n)\asymp \sigma n^{3/2}$ over which the adaptive sorting algorithm does not always work, that is,
\beq
\sup_{(\bTheta, \Pi)\in \calT'_n\times \mathcal{S}'_n}P_{\Theta,\Pi}(\widetilde\Pi\bTheta \widetilde\Pi^\top\ne \Pi\bTheta\Pi^\top)\ge 0.2.
\eeq
This information lower bound suggests that a signal condition stronger than (\ref{suf.cond})  is also necessary for the adaptive sorting to perform well.

\begin{remark}
	From our theoretical analysis, it can be seen that the adaptive sorting algorithm may actually perform well in cases beyond the  ridged monotone Toeplitz class $\calT_n$ considered in this paper. For example, it is shown in Section \ref{eff.thm.sec} (Proposition \ref{eff.thm.2}) that exact matrix reordering can be achieved for all the monotone Toeplitz matrices of the form (\ref{T.layout}) satisfying $\theta_1-\theta_{\lceil n/2\rceil }\ge C\sigma n$ for some positive constant $C$, regardless what is the underlying permutation $\Pi$. See also Section \ref{dis.sec} for possible extensions.
\end{remark}

\section{Interplay Between Computational Efficiency and Statistical Accuracy} \label{tradeoff.sec}

Theorems \ref{mle.thm} to \ref{eff.thm} altogether suggest a fundamental information gap between the statistically optimal procedure and the computationally efficient algorithms. On the one hand, the constrained LSE  has asymptotically the weakest signal strength requirement, but can be computationally infeasible. On the other hand, the adaptive sorting or the spectral seriation algorithm has polynomial running time but requires a strictly stronger signal strength. Thus, it remains unclear whether the existing  gap between the statistical optimal procedure and the computationally efficient procedure is essential and unsurpassable, or it can be reduced, or even closed by devising a better polynomial-time algorithm. 

The  tradeoff between computational efficiency and statistical accuracy has been observed in other permutation-related statistical problems such as sparse/submatrix detection \citep{ma2015computational, cai2020statistical}, structured PCA \citep{cai2013sparse,wang2016statistical}, permuted isotonic regression \citep{mao2020towards, pananjady2020isotonic}, tensor spectral clustering \citep{luo2020tensor}, among many others. In particular, assuming the computational hardness  of the well-known planted clique problem, many of these problems \citep{ma2015computational, wang2016statistical,cai2020statistical,pananjady2020isotonic,luo2020tensor} have been shown to preserve a regime with fundamental computational barrier; that is,  any randomized polynomial-time algorithm must be statistically suboptimal.

In light of these existing work, it is of interest to prove or disprove the existence of any polynomial-time algorithm that succeeds over the region 
\beq
\sigma\sqrt{n\log n}\lesssim \rho^*(\calT'_n,\calS'_n)\lesssim \sigma n^2.
\eeq
Solving this problem requires a quite different set of tools and we leave this fundamental and challenging problem for future investigation.

\section{Numerical Studies} \label{simu.sec}

In this section, we evaluate the empirical performance of the proposed adaptive sorting algorithm and compare it with several existing matrix reordering methods.

\begin{figure}
	\centering
	\includegraphics[angle=0,width=4cm]{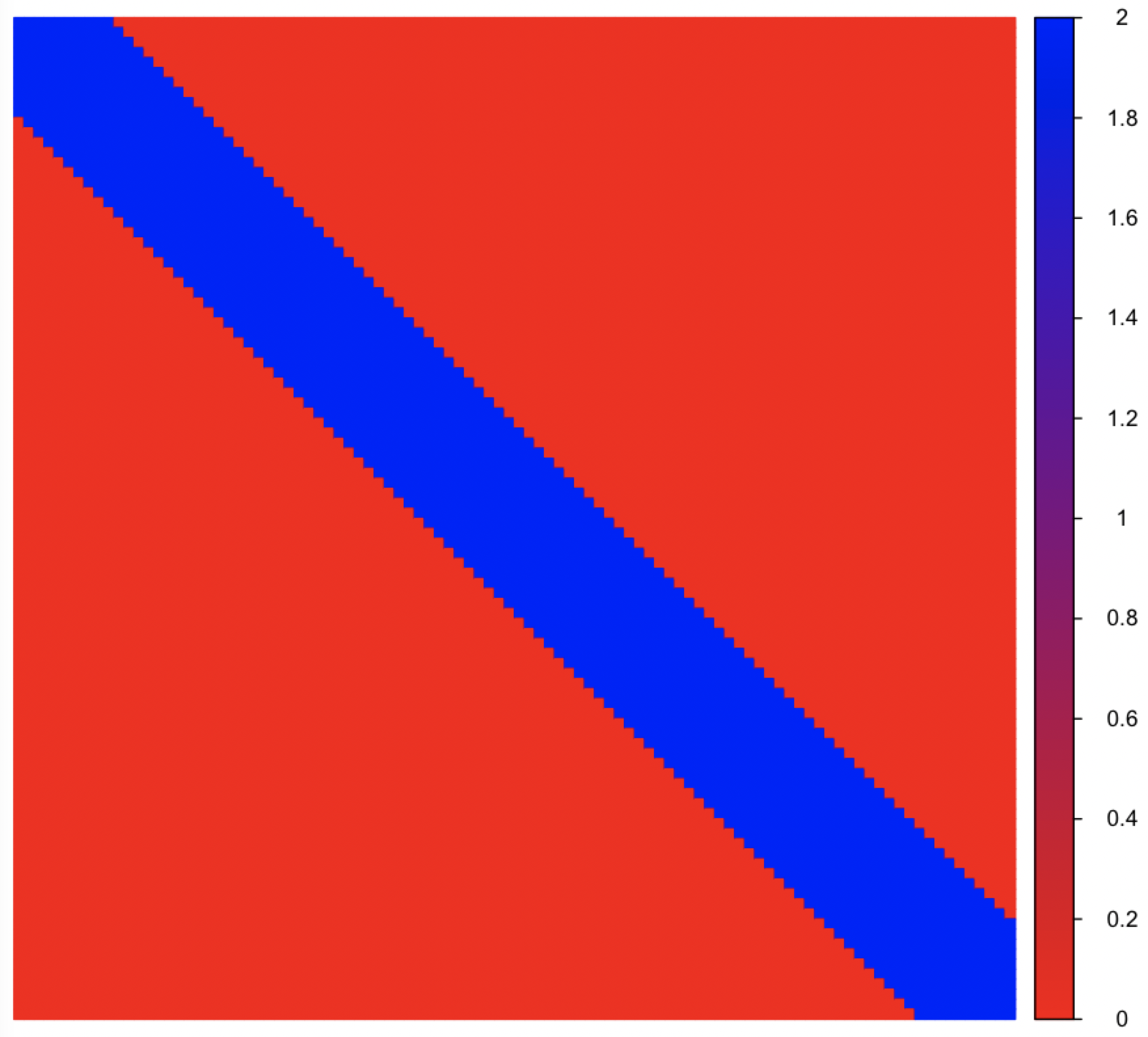}
	\includegraphics[angle=0,width=4cm]{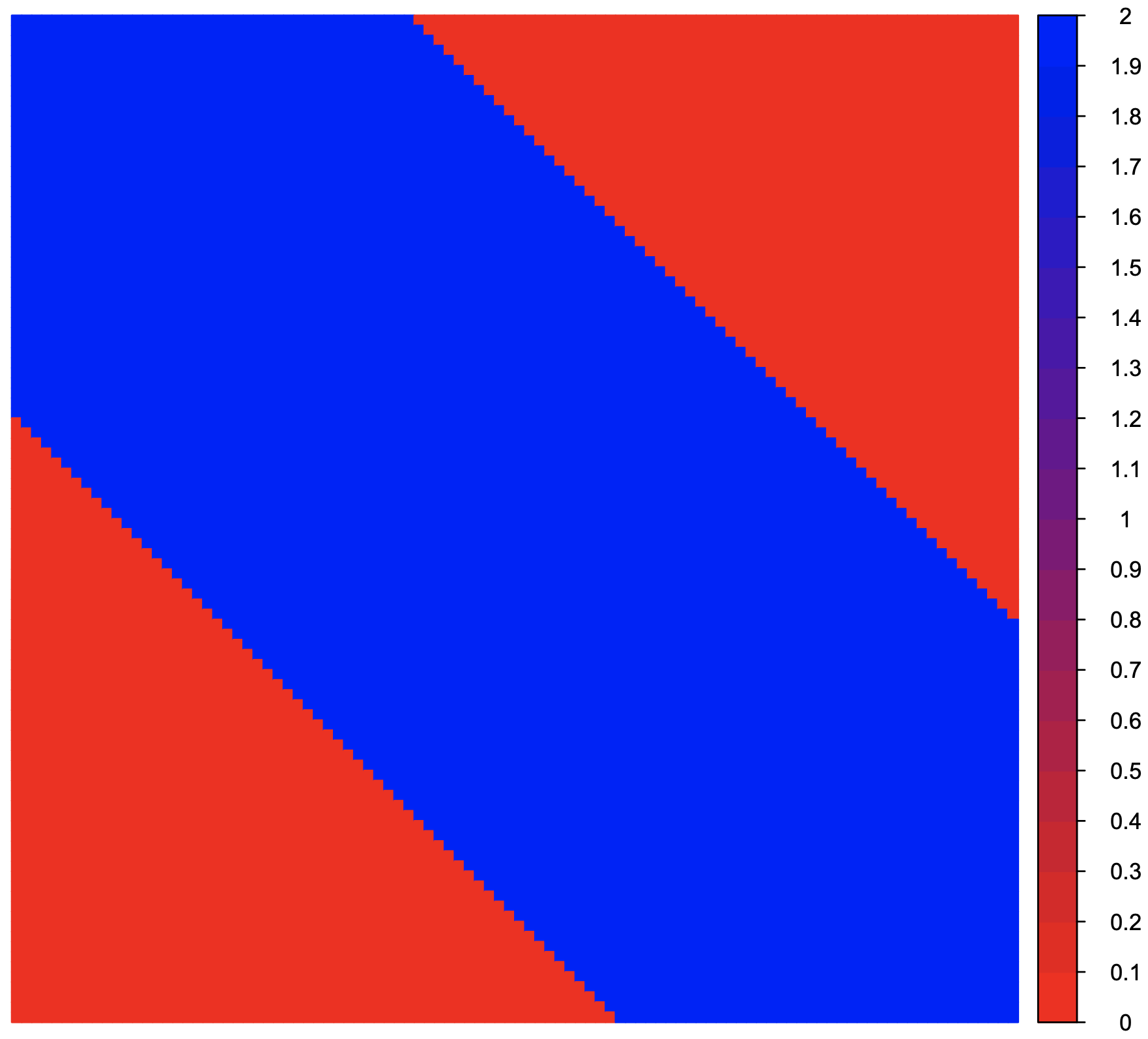}\\
	\includegraphics[angle=0,width=4cm]{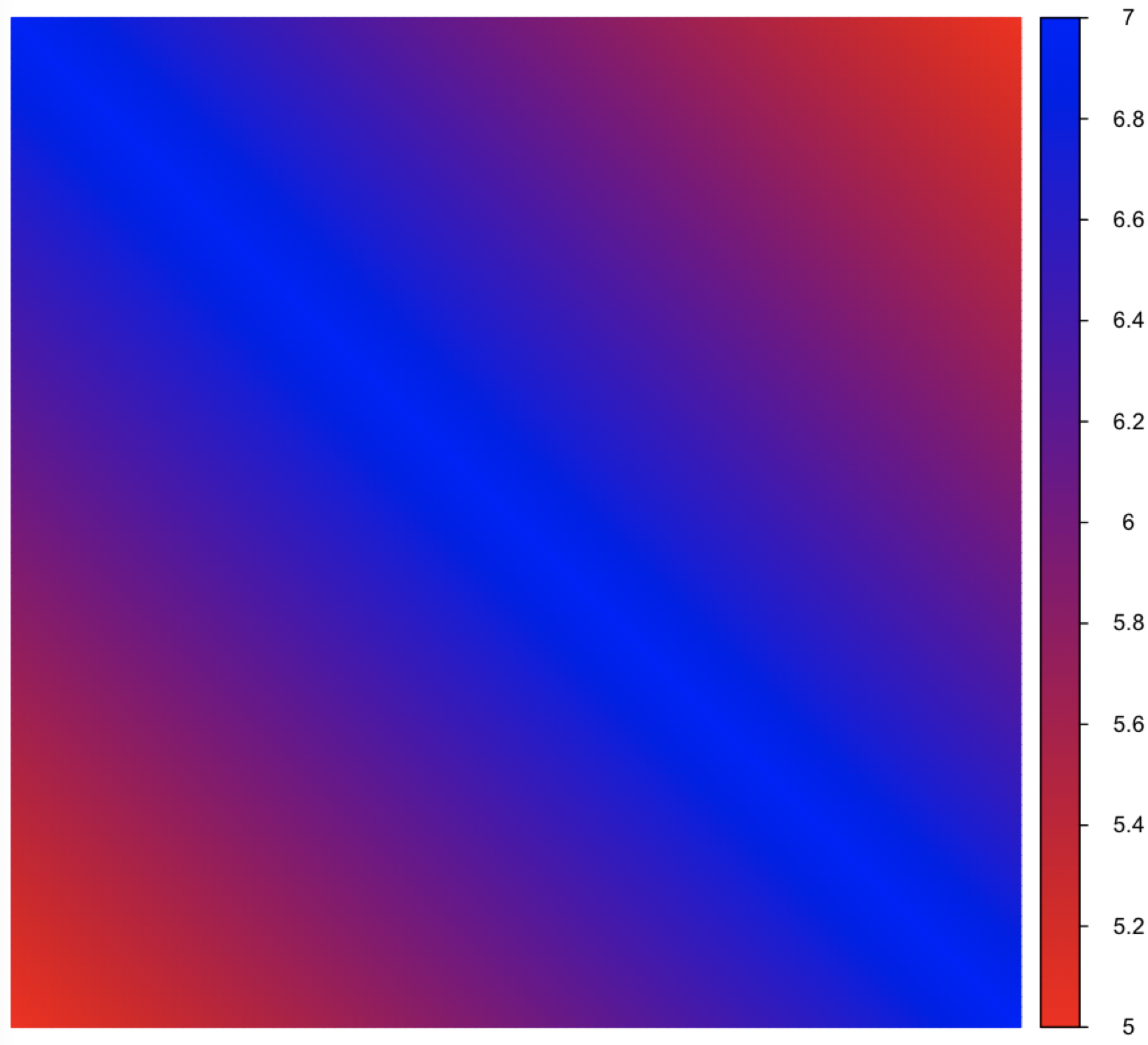}
	\includegraphics[angle=0,width=4cm]{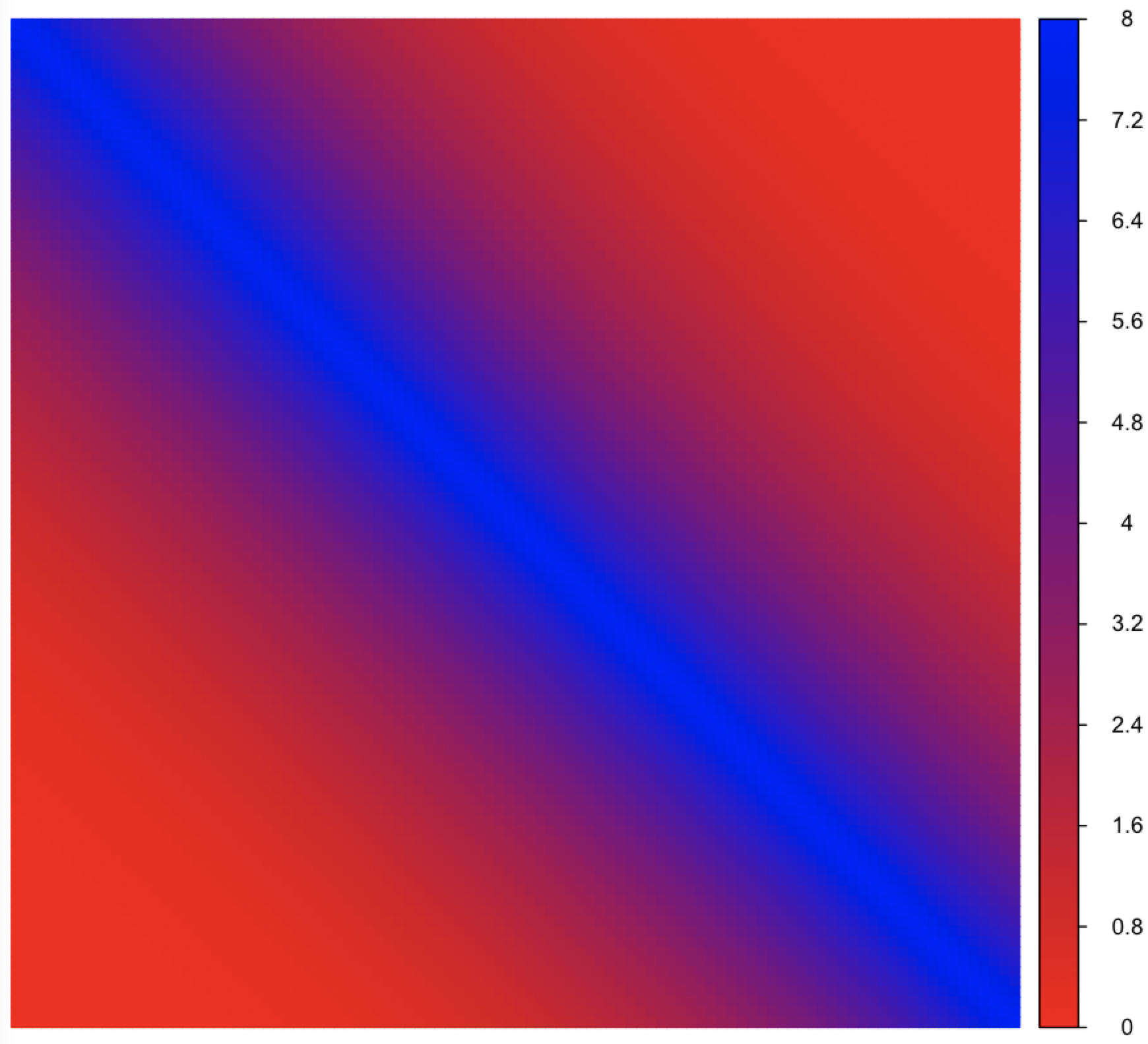}\\
	\includegraphics[angle=0,width=4cm]{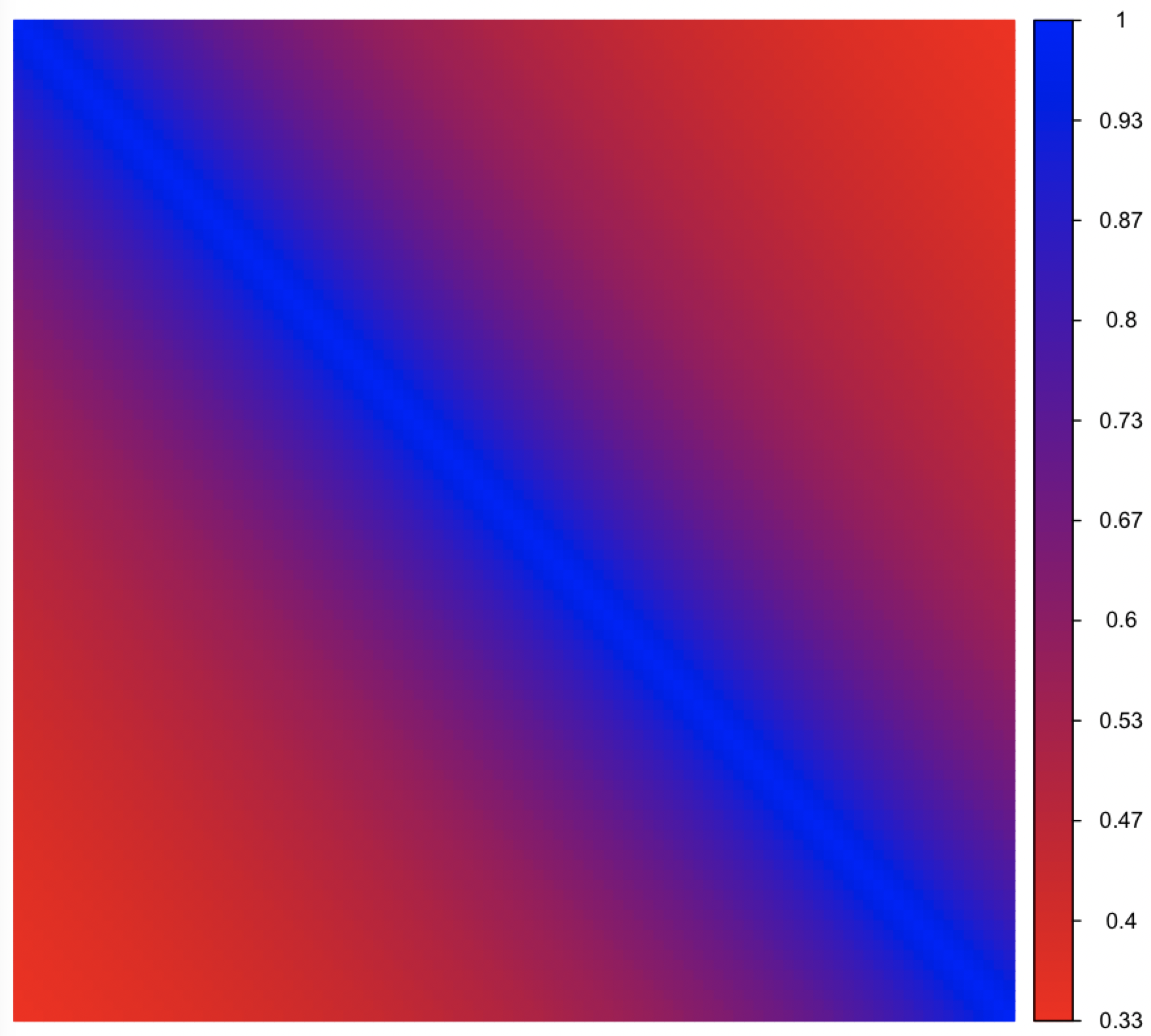}
	\includegraphics[angle=0,width=4cm]{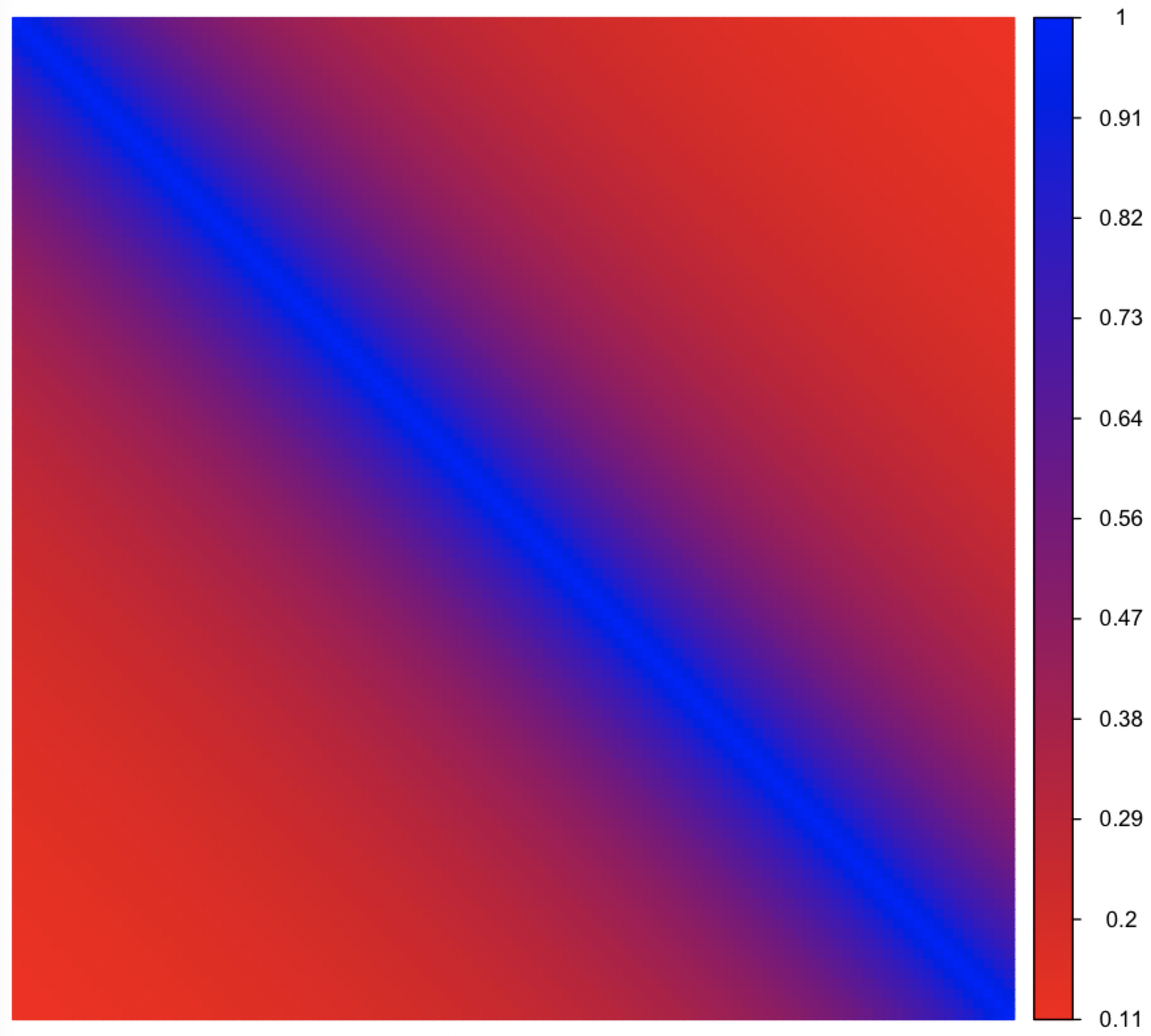}
	\caption{Visualizations of six true signal matrices. Settings 1 to 6 arranged by row.} 
	\label{simu.sig.fig}
\end{figure}

We first set the dimensionality $n=100$ for each matrix, and generate the noise matrix with $i.i.d.$ entries from either a Gaussian distribution $N(0,\sigma^2)$, or a heavier-tailed Laplace distribution $Lap(0,\sigma)$. To better assess the range of applicability of different methods, we consider the following six settings for the underlying true signal matrix (Figure \ref{simu.sig.fig}), containing band, linear decaying and nonlinear decaying monotone Toeplitz matrices: \\
(1) Narrow-band matrix: $\bTheta\in\R^{n\times n}$ is Toeplitz of the form (\ref{T.layout}), where $\theta_1=\theta_2=...=\theta_{10}=2$, and $\theta_i=0$ for all $i\ge11$; (2) Wide-band matrix: $\bTheta\in\R^{n\times n}$ is Toeplitz of the form (\ref{T.layout}), where $\theta_1=\theta_2=...=\theta_{40}=2$, and $\theta_i=0$ for all $i\ge 41$; (3) Linear decaying matrix: $\bTheta\in\R^{n\times n}$ is Toeplitz  of the form (\ref{T.layout}), where $\theta_i=5+0.02\cdot (n-i)$ for all $i\in[n-1]$; (4) Polynomial decaying matrix: $\bTheta\in\R^{n\times n}$ is Toeplitz  of the form (\ref{T.layout}), where $\theta_i=[(n-i)\cdot 0.02]^3$  for all $i\in[n-1]$; (5) Inverse linear decaying matrix: $\bTheta\in\R^{n\times n}$ is Toeplitz  of the form (\ref{T.layout}), where $\theta_i=[1+0.02\cdot i]^{-1}$  for all $i\in[n-1]$; (6) Inverse polynomial decaying matrix: $\bTheta\in\R^{n\times n}$ is Toeplitz  of the form (\ref{T.layout}), where $\theta_i=[1+0.02\cdot i]^{-2}$  for all $i\in[n-1]$.

For each setting, we evaluate the performance of various methods by comparing their ability to achieve exact order recovery across a range of noise levels. The underlying permutations are generated uniformly from the permutation group $\calS_n$. The performance of each algorithm is measured by the empirical proportion of failures in exact order recovery (as determined by the loss function $\tau_{\bTheta}(\Pi,\Pi')$, where we do not distinguish between the true permutation and its complete reversal) over 500 rounds of simulations at each noise level $\sigma$.

In addition to the adaptive sorting algorithm ("AS") proposed in Section \ref{ipdd.sec} and the spectral seriation algorithm ("SS")  defined in Section \ref{ss.sec}, we also evaluate the following five  existing matrix reordering methods:
\begin{itemize}
	\item Best permutation analysis ("BP") proposed by \cite{rajaratnam2013best}.
	\item The multidimensional scaling based method ("MDS") implemented by the function \texttt{seriate} in the R  package \texttt{seriation}, with the argument option \texttt{method="MDS"} \citep{seriation}.
	\item The rank-two ellipse seriation algorithm ("R2E") proposed by \cite{chen2002generalized}, and implemented by the function \texttt{seriate} in  the R  package \texttt{seriation}, with the option \texttt{method="R2E"}.
	\item The "VAT" (visual assessment of tendency) algorithm proposed by \cite{bezdek2002vat},  and implemented by the function \texttt{seriate} in  the R  package \texttt{seriation}, with the option \texttt{method="VAT"}.
	\item The normalized spectral seriation algorithm ("SS.n"), which  differs from the spectral seriation algorithm only in its definition of the Laplacian matrix $\bL={\bf I}_n-\bD^{-1}\bY$.
\end{itemize}

\begin{figure}[h!]
	\centering
	\includegraphics[angle=0,height=3.5cm]{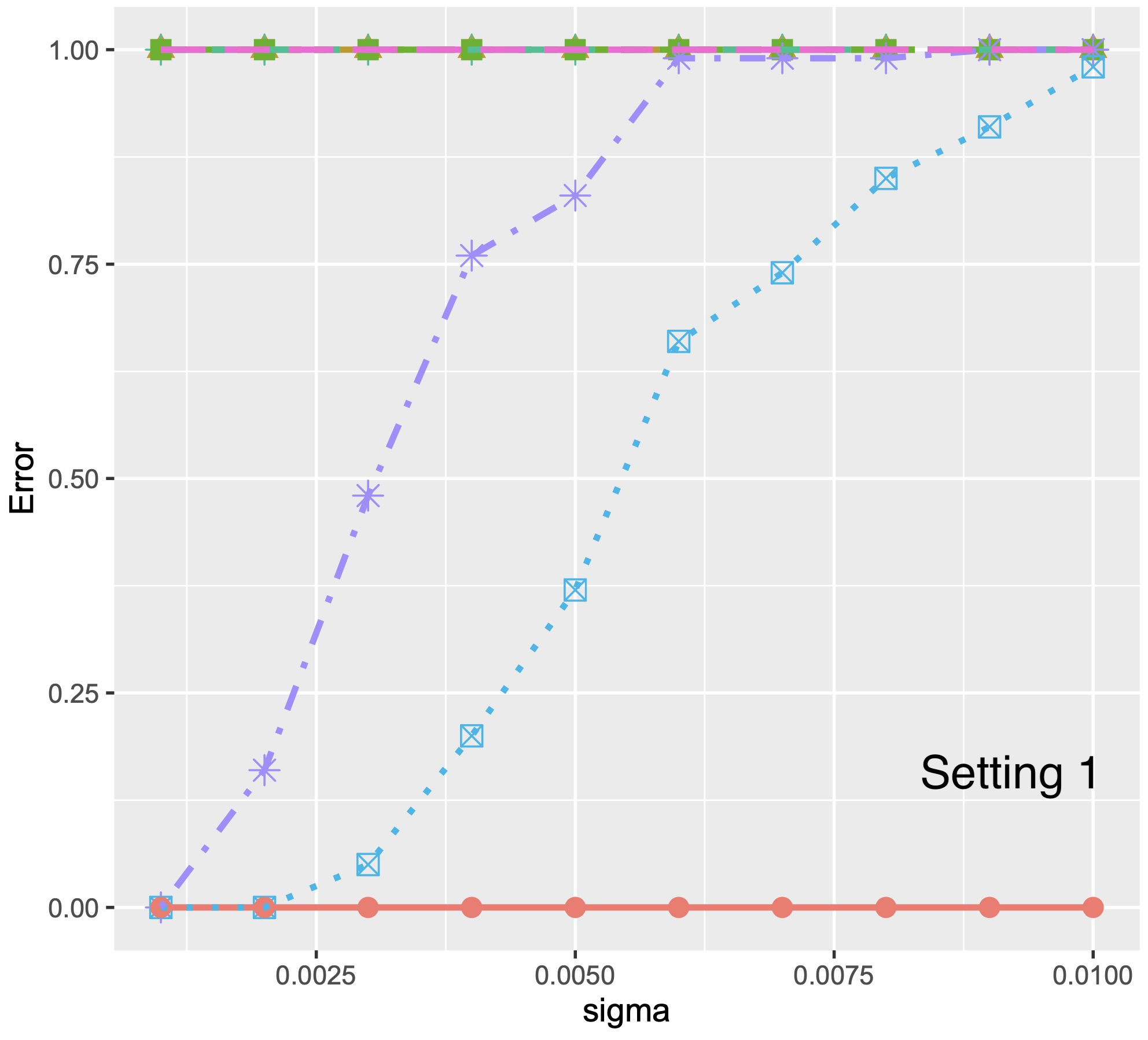}
	\includegraphics[angle=0,height=3.5cm]{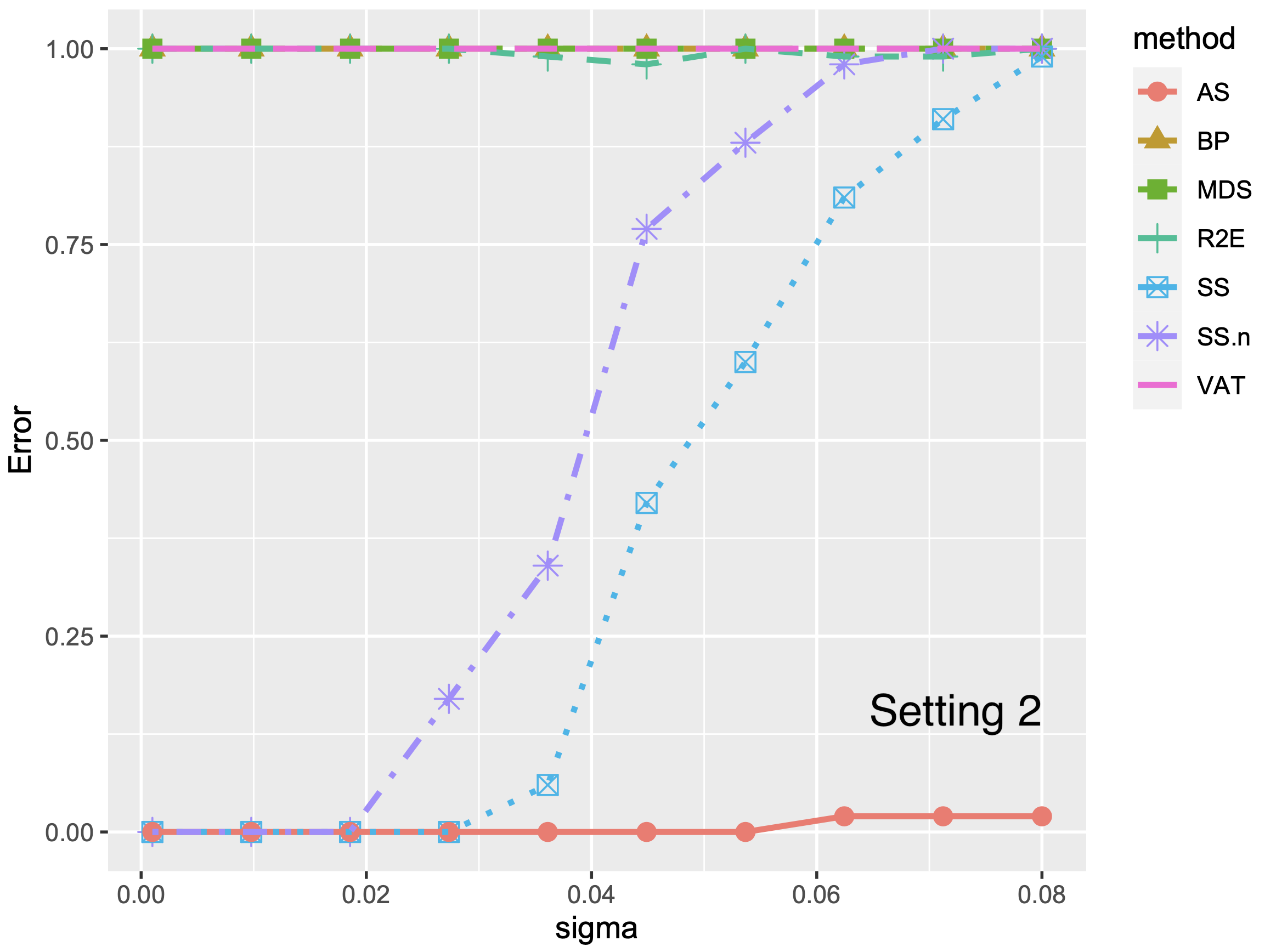}
	\includegraphics[angle=0,height=3.5cm]{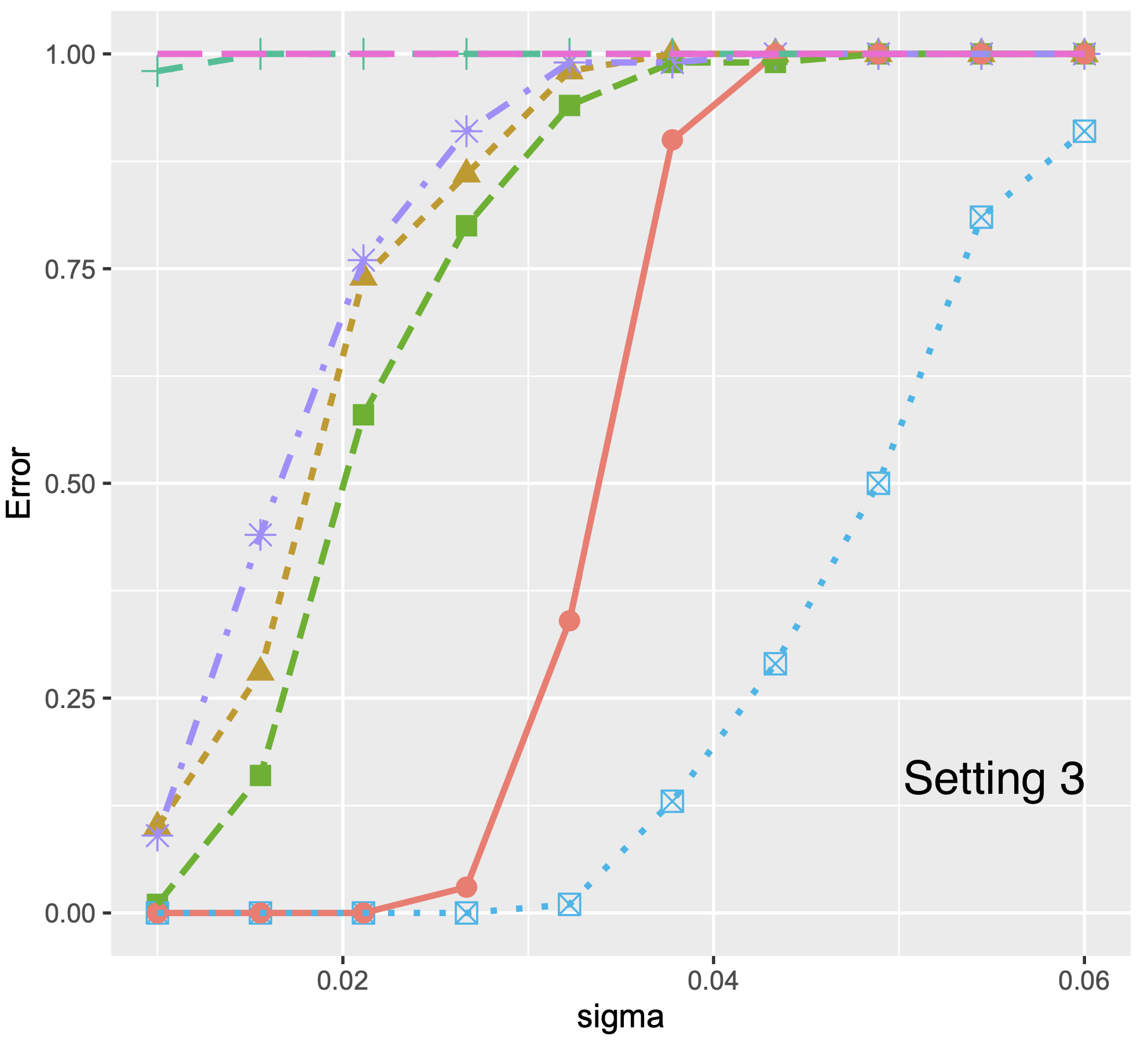}
	\includegraphics[angle=0,height=3.5cm]{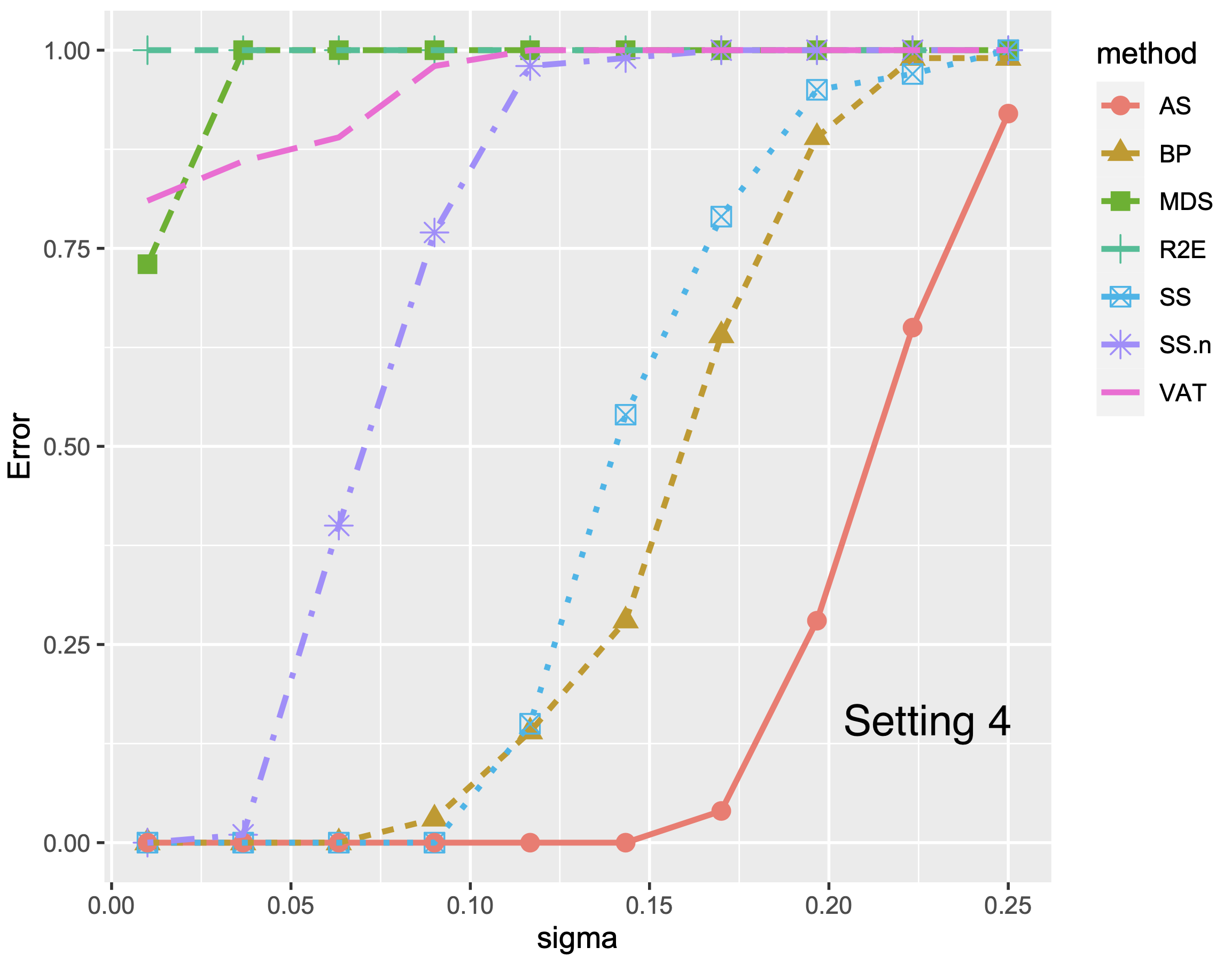}
	\includegraphics[angle=0,height=3.5cm]{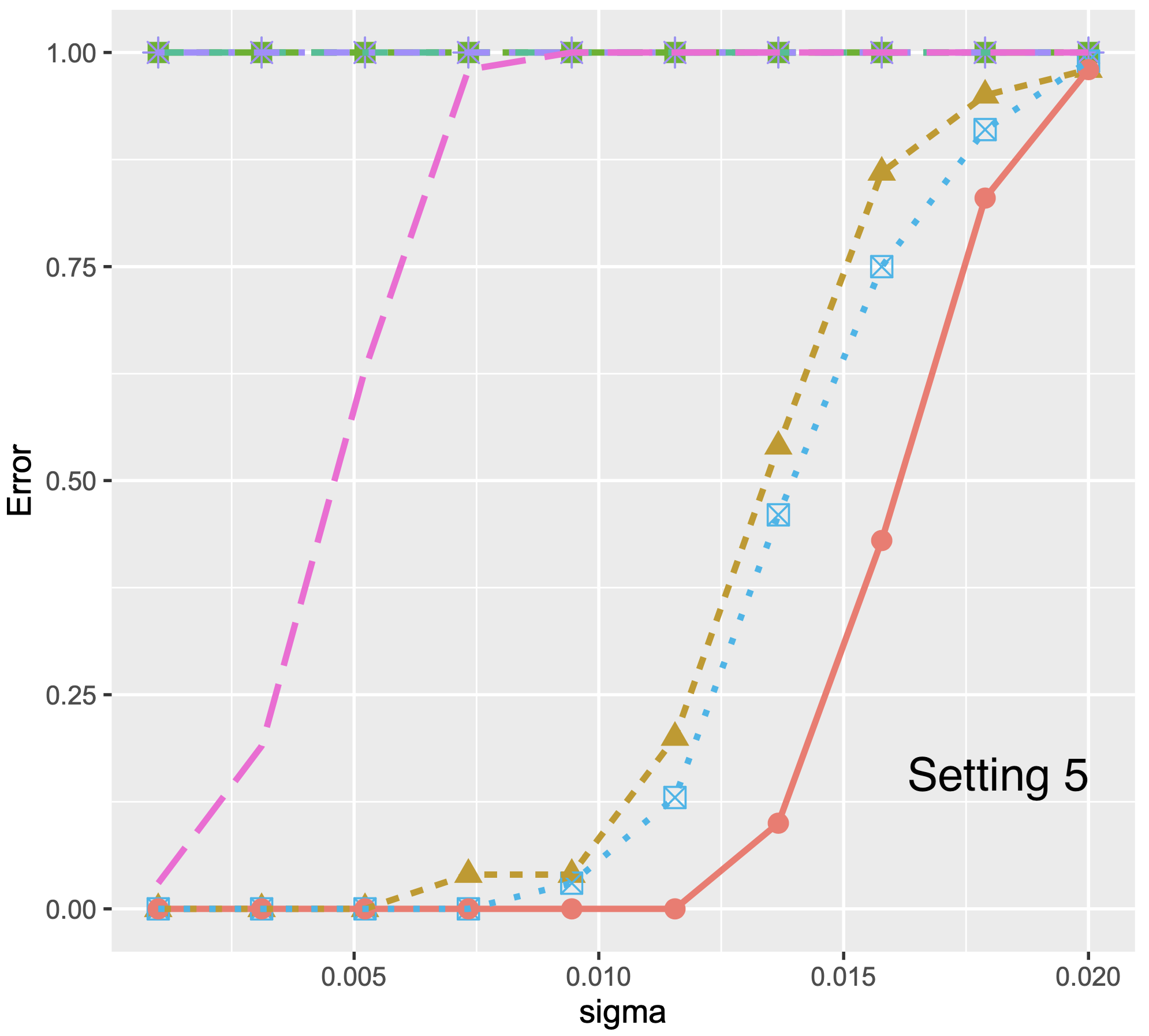}
	\includegraphics[angle=0,height=3.5cm]{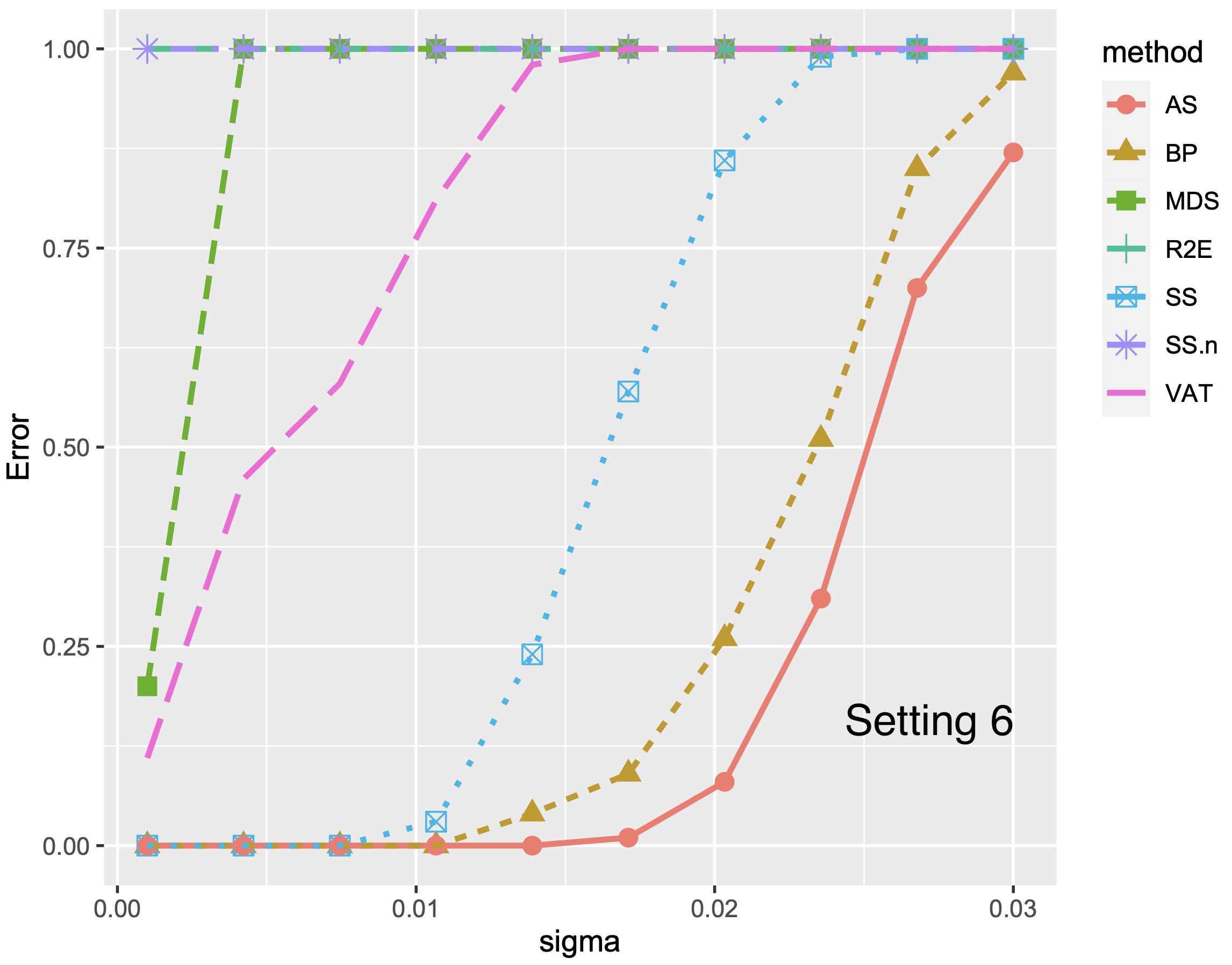}
	\caption{Comparison of seven matrix reordering methods under the Gaussian noise with $n=100$ and variance $\sigma^2$. The errors quantify the empirical probability of failures in exact order recovery over 500 rounds of simulations. AS: adaptive sorting, BP: best permutation analysis, MDS: multidimensional scaling based seriation, R2E: rank-two ellipse seriation, SS: spectral seriation, SS.n: normalized spectral seriation, VAT: visual assessment of tendency. } 
	\label{simu.fig}
\end{figure}

\begin{figure}[h!]
	\centering
	\includegraphics[angle=0,height=3.5cm]{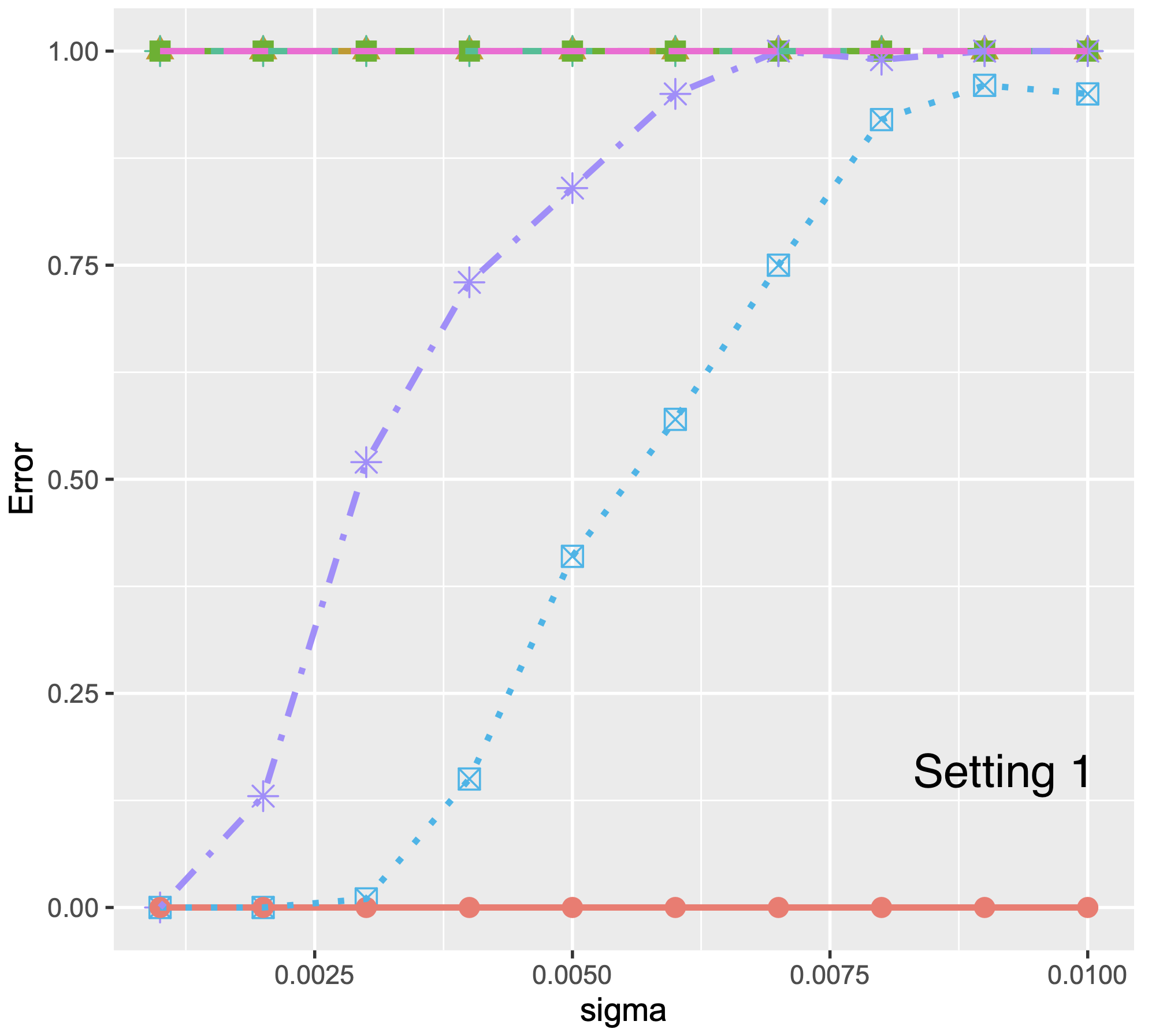}
	\includegraphics[angle=0,height=3.5cm]{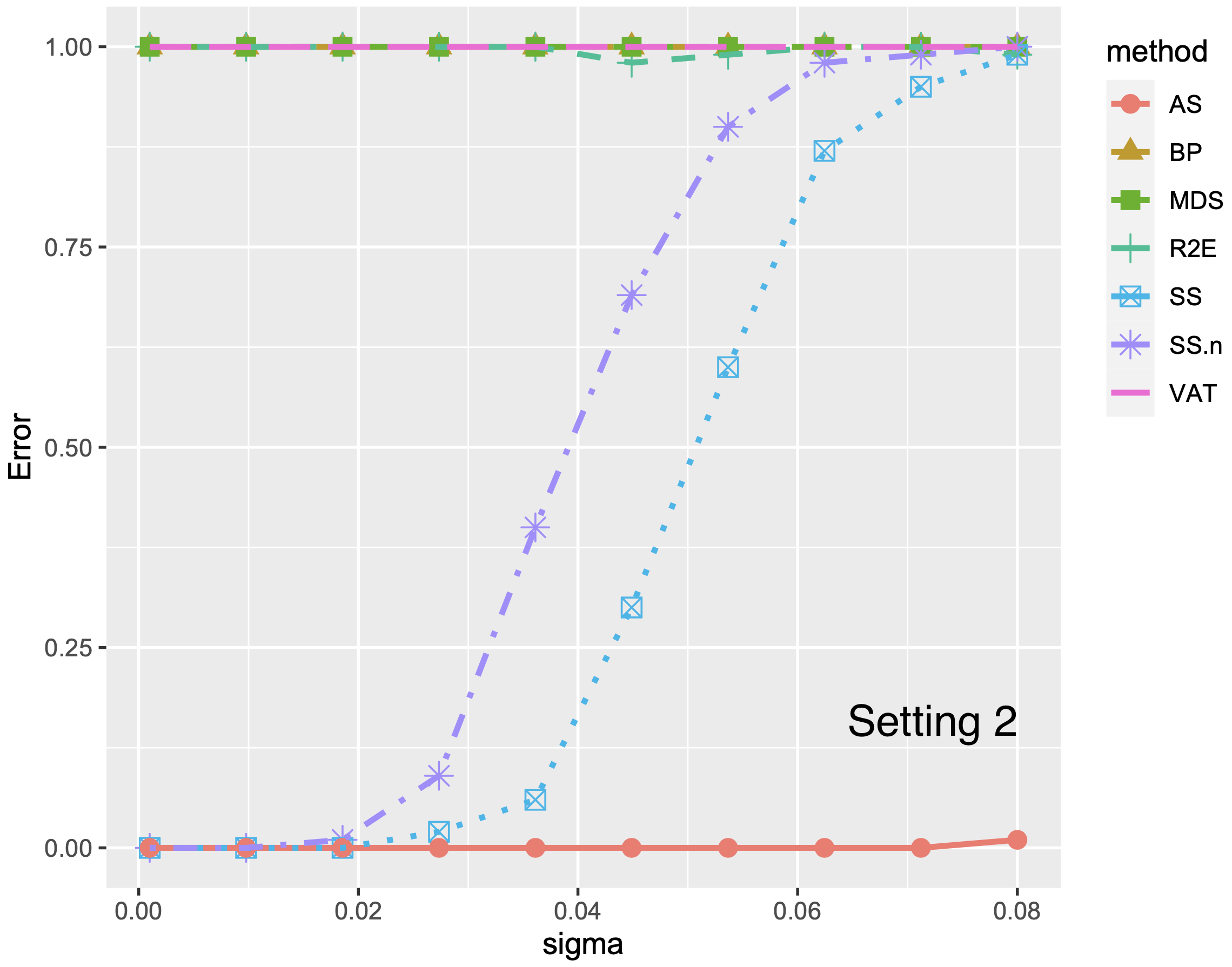}\\
	\includegraphics[angle=0,height=3.5cm]{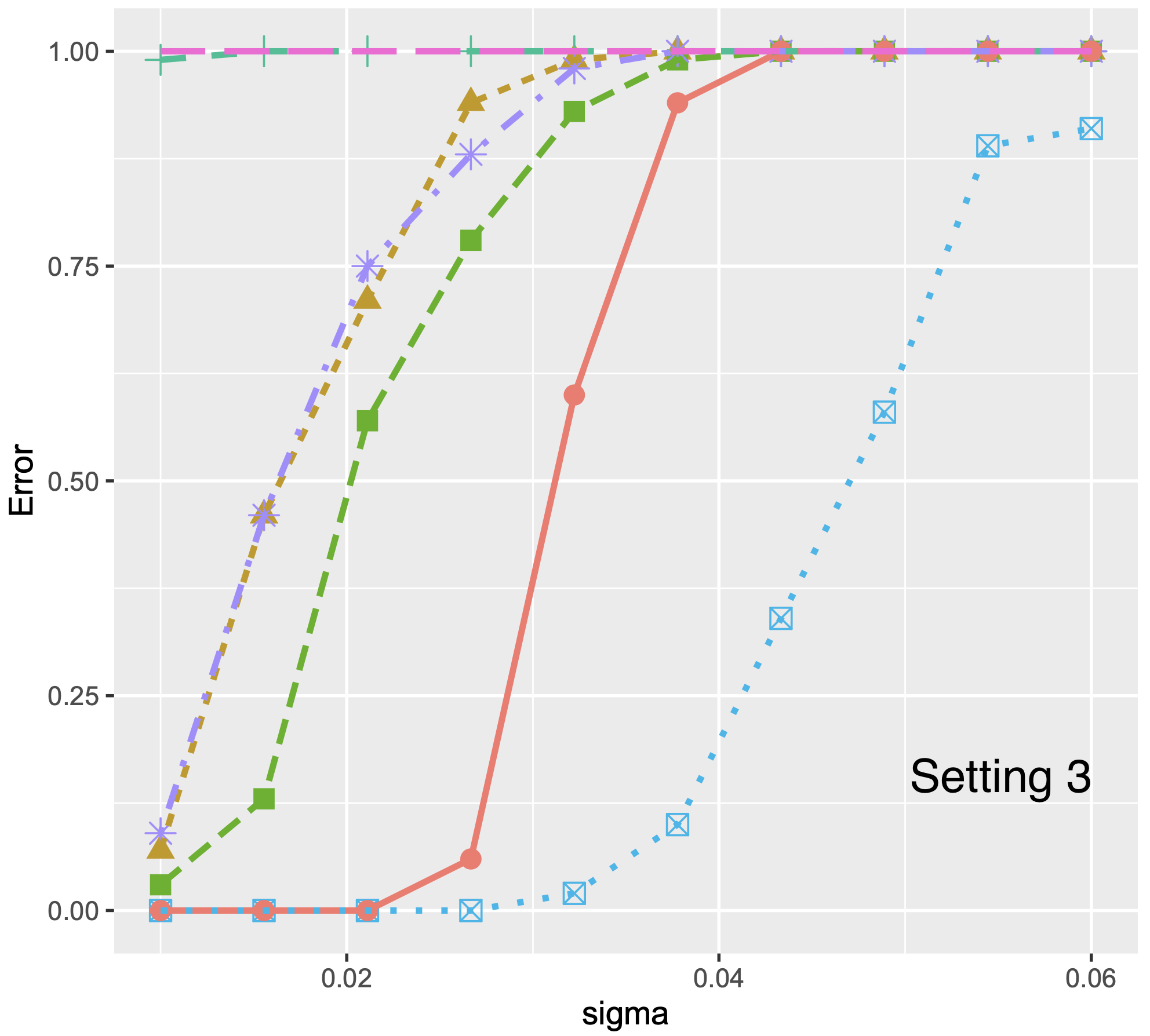}
	\includegraphics[angle=0,height=3.5cm]{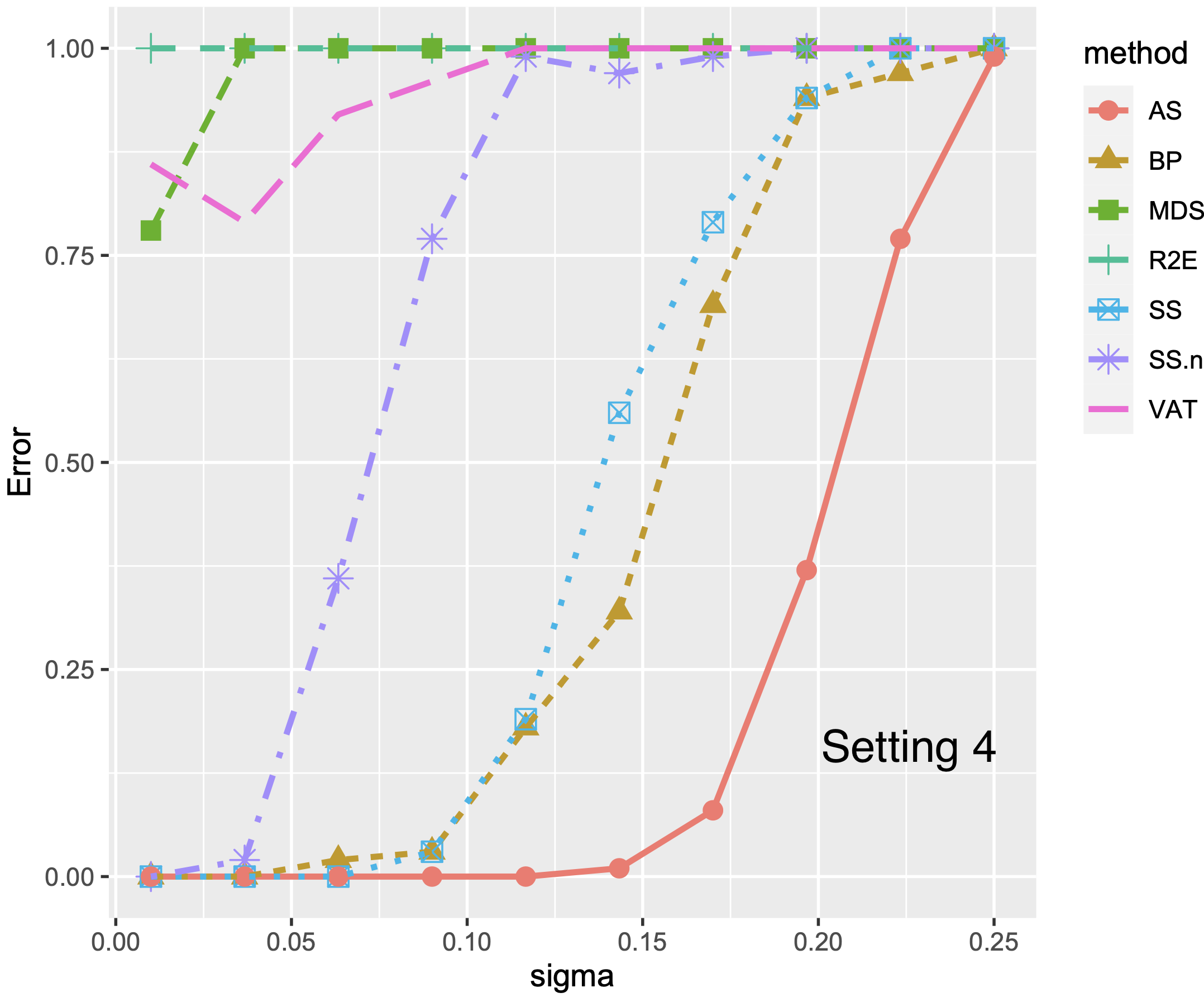}\\
	\includegraphics[angle=0,height=3.5cm]{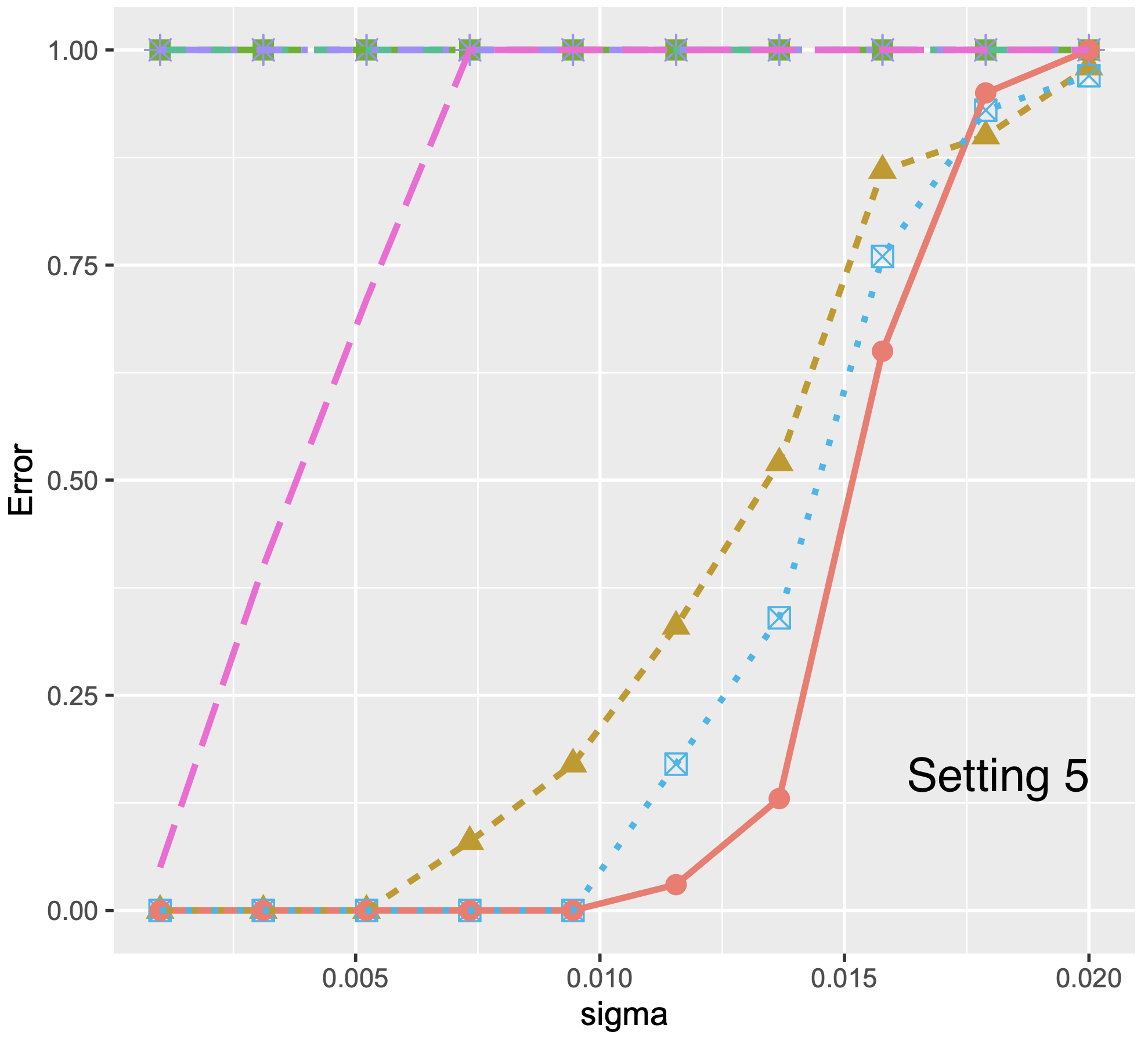}
	\includegraphics[angle=0,height=3.5cm]{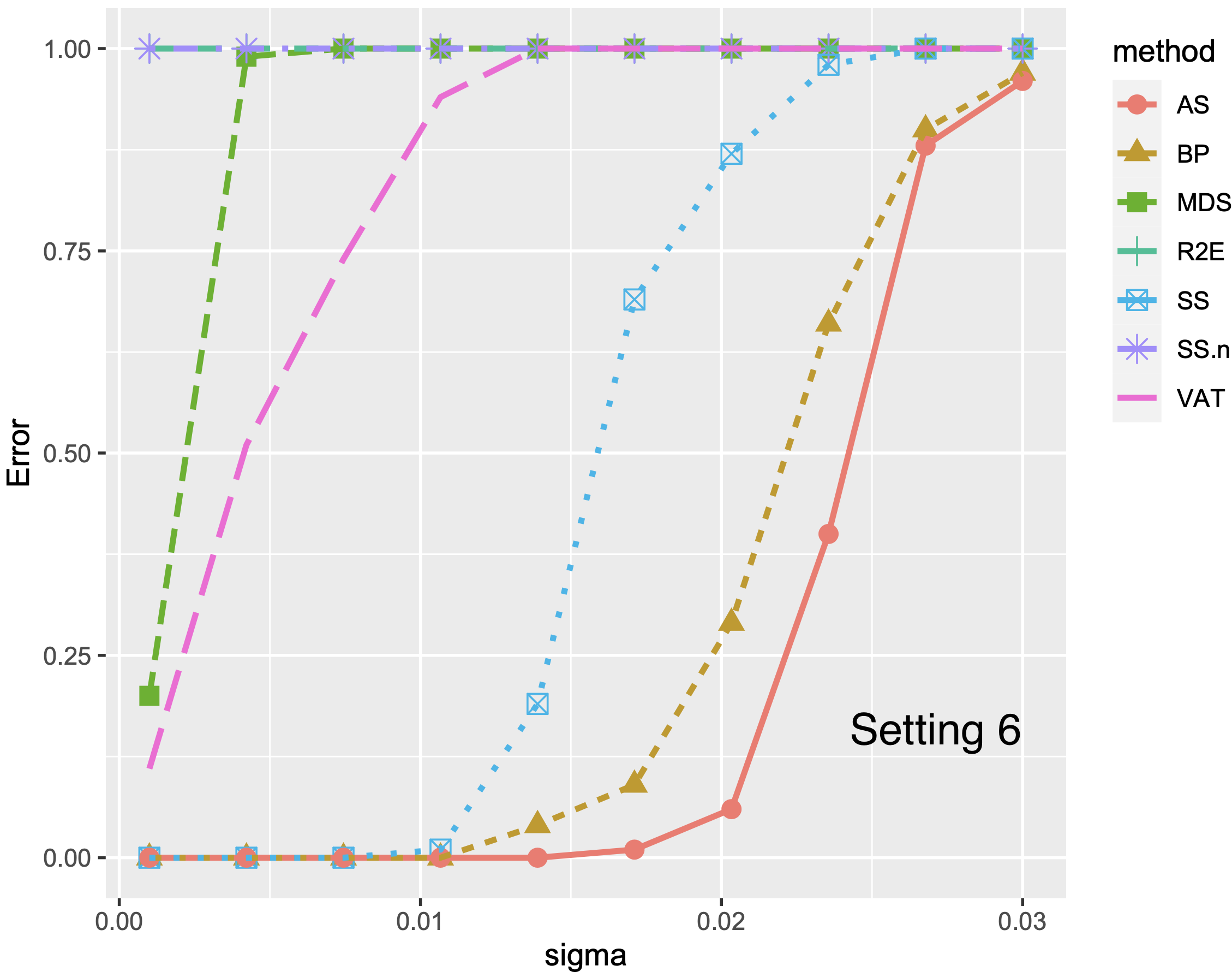}
	\caption{Comparison of seven matrix reordering methods under the Laplacian noise with $n=100$ and scale parameter $\sigma$ (i.e., with variance $2\sigma^2$).} 
	\label{simu.fig2}
\end{figure}

The numerical results for both noise settings are presented in Figures \ref{simu.fig} and \ref{simu.fig2}. Overall, the proposed AS algorithm demonstrates the best performance across all settings, followed by SS in most cases. In particular, AS outperforms other methods most significantly for the band Toeplitz matrices (Settings 1 and 2). Under the other four settings with linear or nonlinear decaying diagonals, AS, SS, and often BP exhibit relatively better performance compared to the other methods. Among these three methods, AS consistently performs better than SS and BP in the nonlinear decaying cases (Settings 4 to 6), while in the strict linear decaying case (Setting 3), SS demonstrates the best performance, followed by AS.

\begin{figure}[h!]
	\centering
	\includegraphics[angle=0,height=3.5cm]{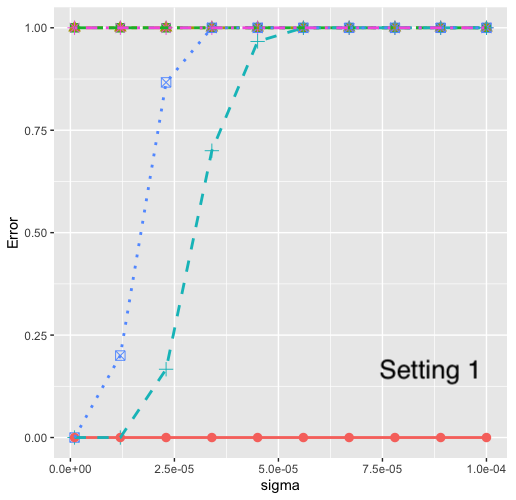}
	\includegraphics[angle=0,height=3.5cm]{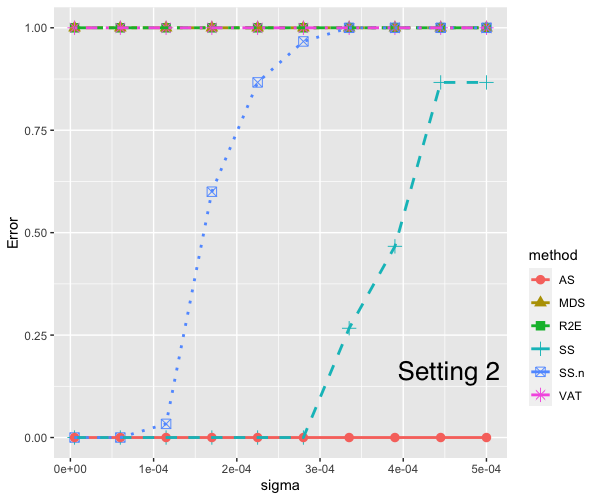}\\
	\includegraphics[angle=0,height=3.5cm]{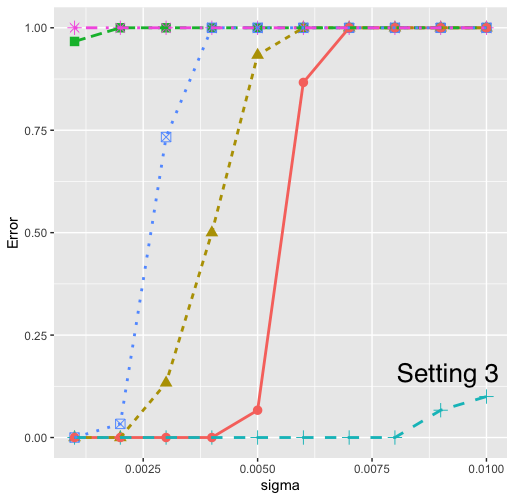}
	\includegraphics[angle=0,height=3.5cm]{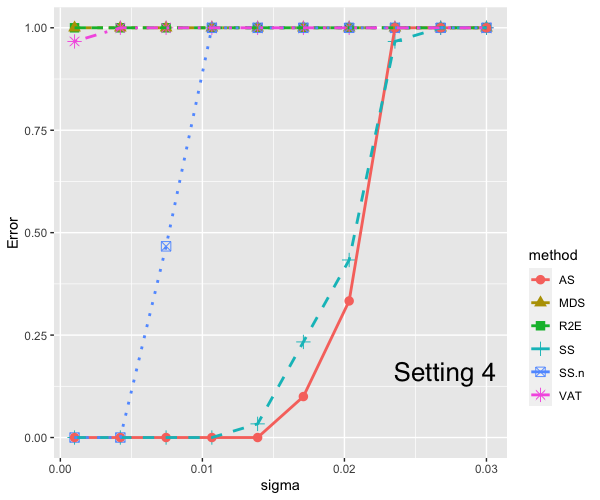}\\
	\includegraphics[angle=0,height=3.5cm]{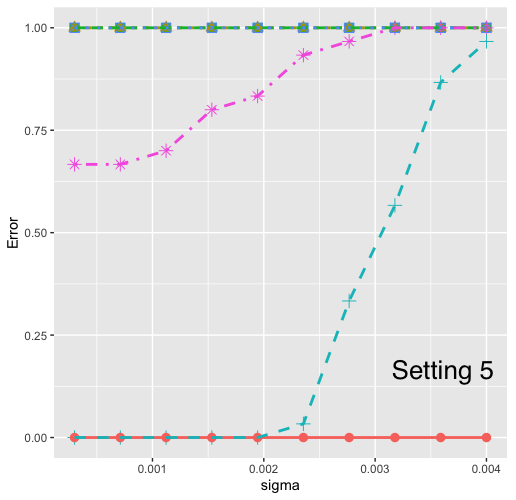}
	\includegraphics[angle=0,height=3.5cm]{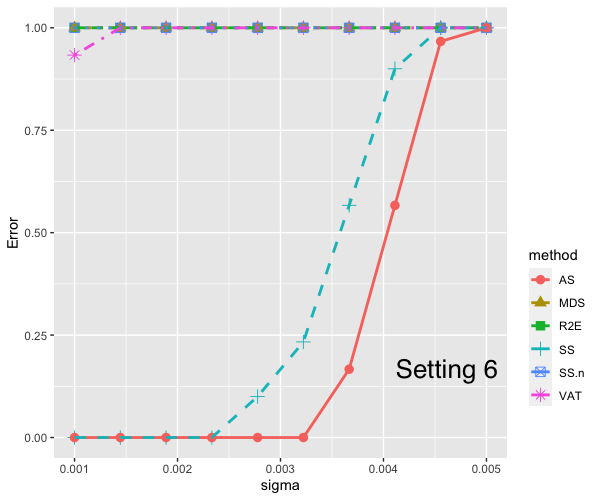}
	\caption{Comparison of six matrix reordering methods (``BP" removed due to lack of scalability) under the Gaussian noise with $n=1000$ and variance $\sigma^2$.} 
	\label{simu.fig3}
\end{figure}

To further assess the relative performance for large matrices, we repeated the experiments under Gaussian noise with $n=1000$. Notably, BP was excluded from this comparison due to its lack of scalability for large $n$, as the algorithm requires $n^2$ iterations, and in each iteration, the determinant of a matrix up to size $n\times n$ must be evaluated. Our results in Figure \ref{simu.fig3} suggest similar trends to those observed for $n=100$, highlighting the superior performance of AS in most settings and the best performance achieved by SS in Setting 3. Comparing Figures \ref{simu.fig} and \ref{simu.fig3}, particularly the values of $\sigma$ at which the phase transition occurs, we observe that, in line with our theory, the performance of AS and SS deteriorates as $n$ increases from 100 to 1000.

The simulation studies demonstrate the overall superiority and adaptivity of the proposed AS algorithm over the other six alternative methods. In particular, the numerical results indicate SS to be the overall best existing method, whose empirical performance is in turn dominated by the proposed AS algorithm in most cases. This phenomenon is consistent with our theoretical analysis of the two methods, showing the strength and practical relevance of the theoretical results developed in the preceding sections.

\section{Application to Two Real Datasets} \label{data.sec}

We analyze two real single-cell RNA sequencing datasets, and compare the performance of the adaptive sorting and the spectral seriation algorithms for inferring the latent pseudotemporal orders of  single cells.  

\begin{figure}[h!]
	\centering
	\includegraphics[angle=0,width=4cm]{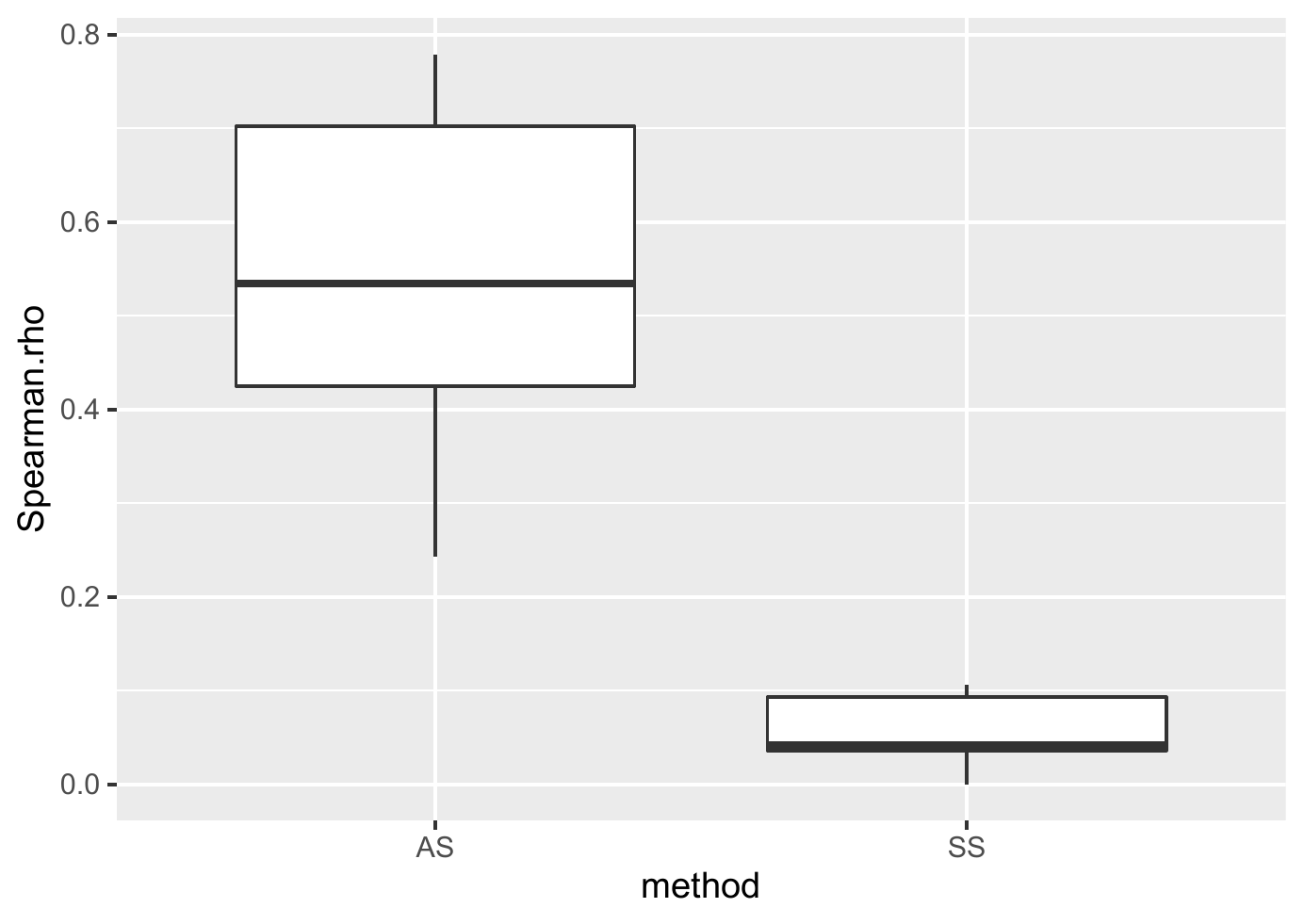}
		\includegraphics[angle=0,width=4cm]{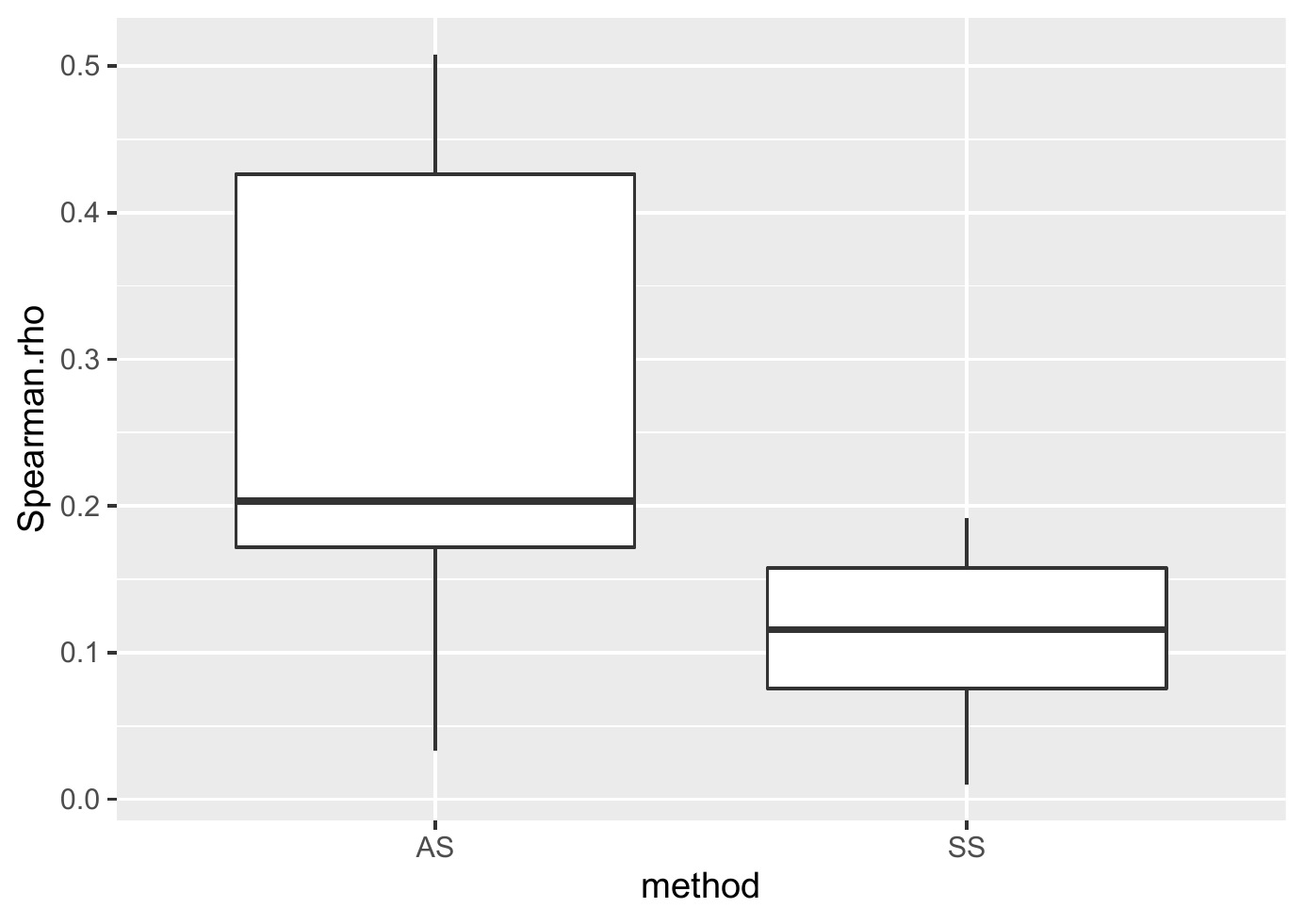}\\
				\includegraphics[angle=0,width=4cm]{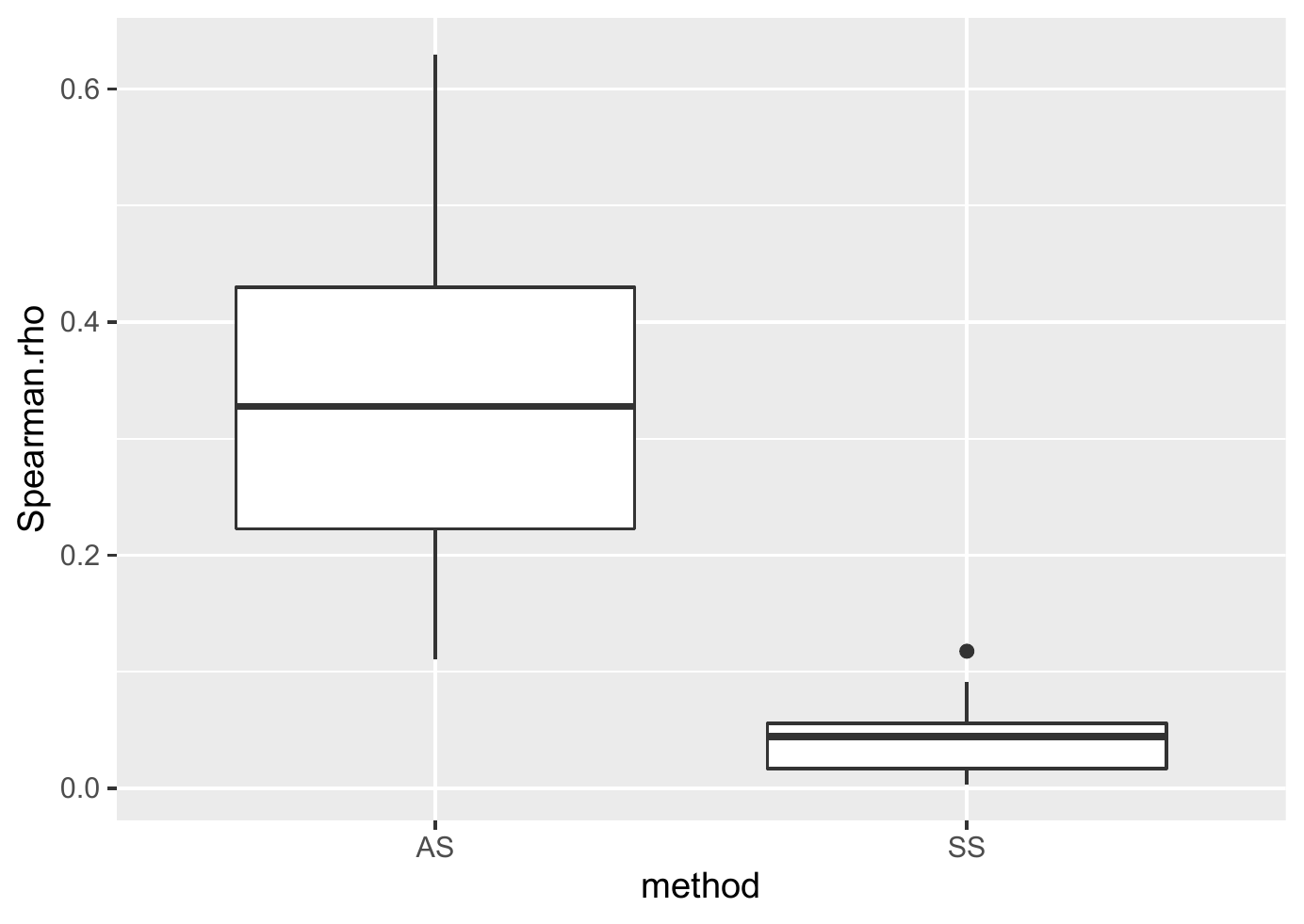}
		\includegraphics[angle=0,width=4cm]{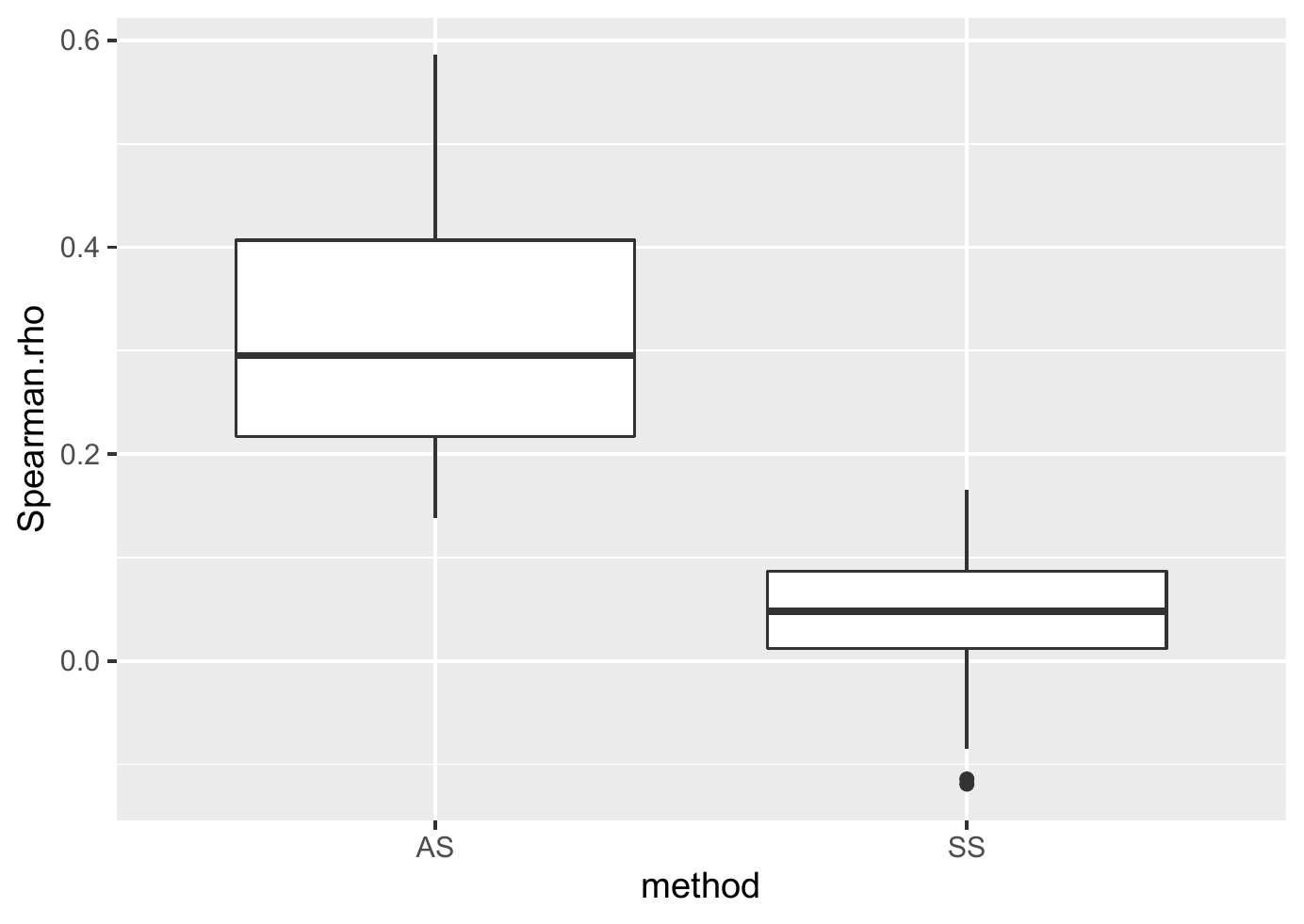}\\
	\includegraphics[angle=0,width=4cm]{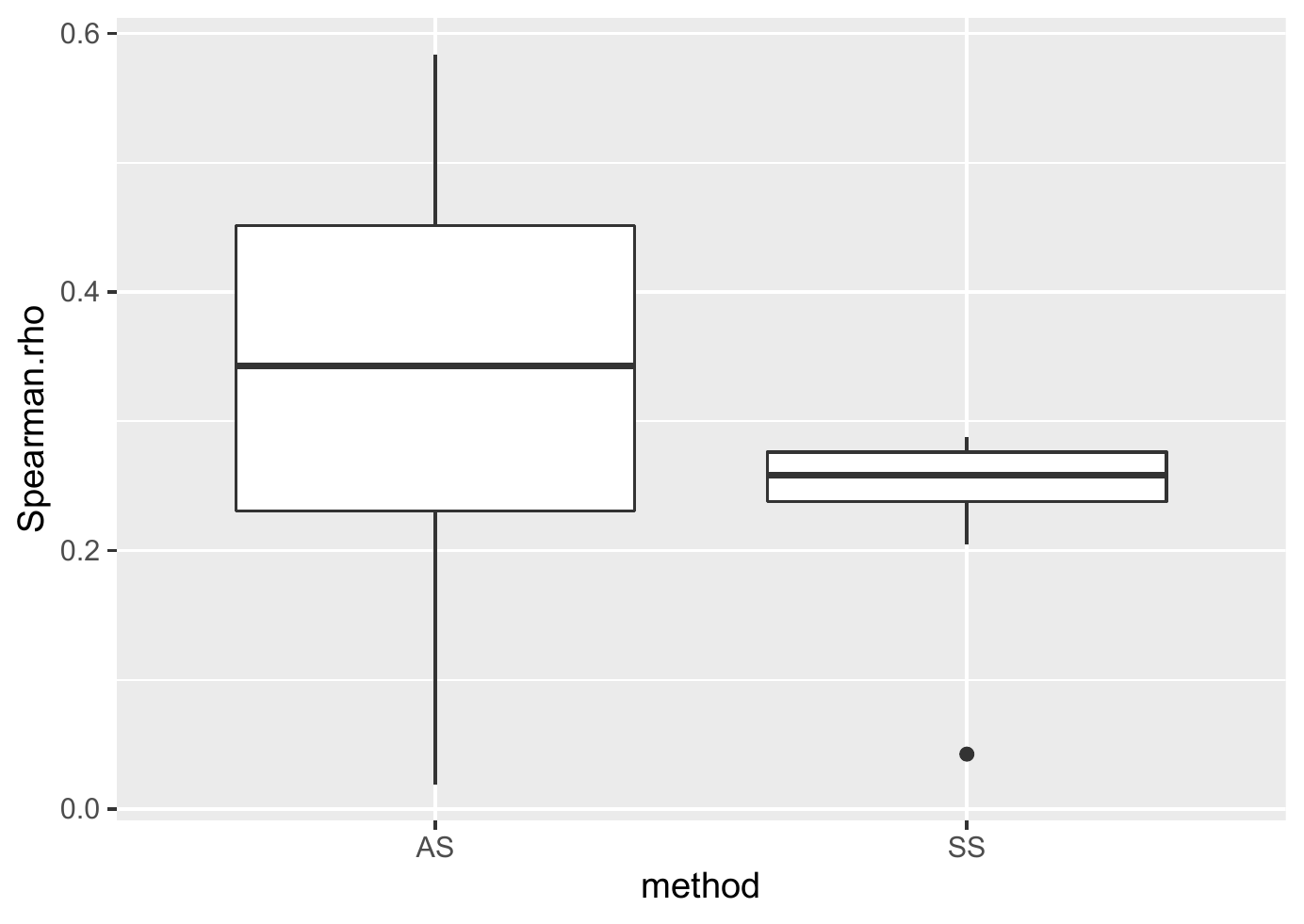}
	\includegraphics[angle=0,width=4cm]{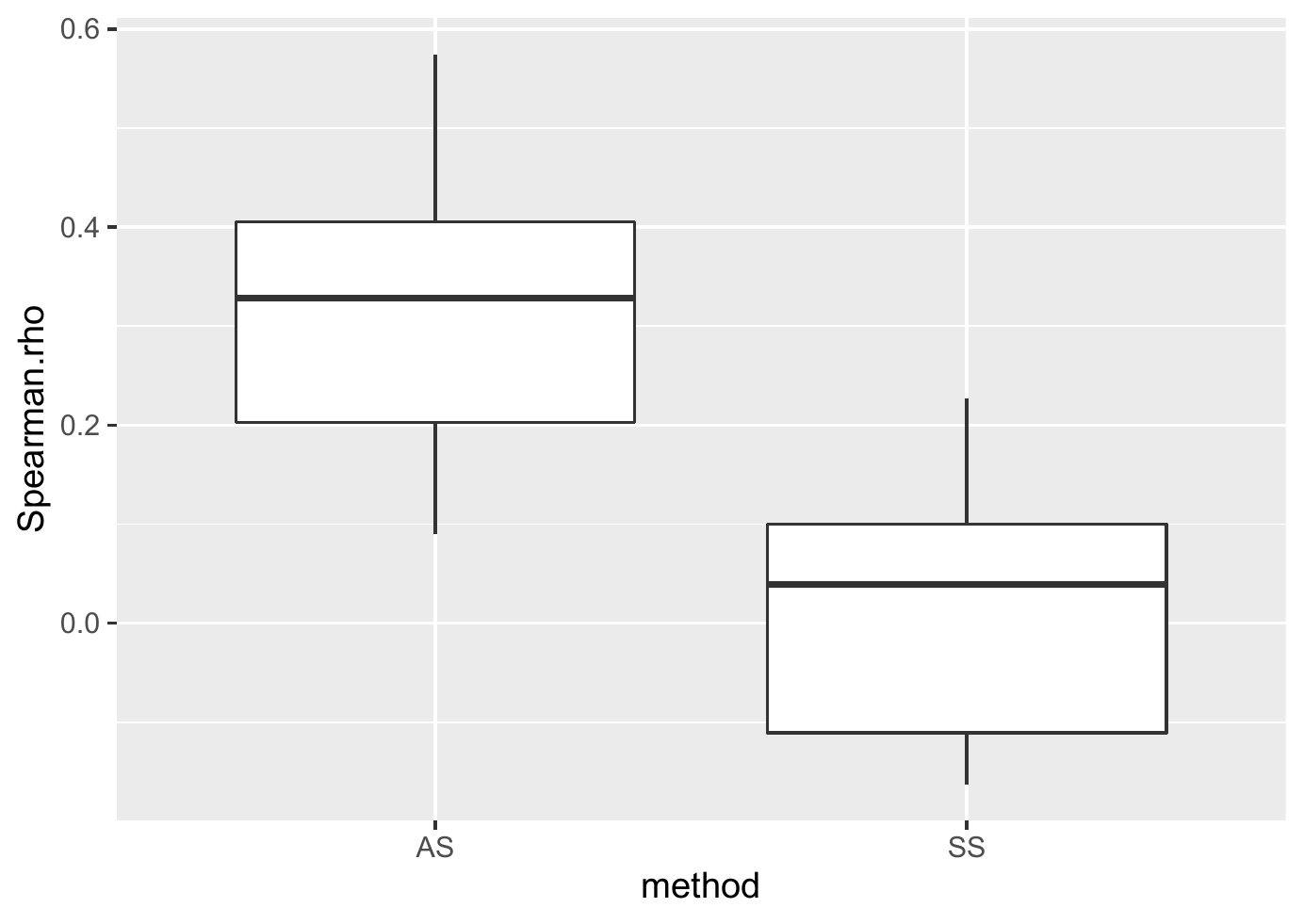}
	\caption{Comparison of adaptive sorting (AS) and spectral seriation (SS) for inferring temporal orders of single cells using Spearman's rho when $p=500$ (top), $p=1500$ (middle), and $p=2500$ (bottom). Left: study of primary human myoblasts \citep{trapnell2014dynamics}. Right: study of  human primordial germ cells \citep{guo2015transcriptome}.} 
	\label{comp.fig}
\end{figure}

The first dataset contains single-cell mRNA sequencing reads for 372 primary human skeletal muscle myoblasts undergoing differentiation \citep{trapnell2014dynamics}. Specifically, primary human myoblasts were cultured in high-serum medium; after switching to low-serum medium that induces differentiation, the cells were dissociated and  individually captured at 24-h intervals (0, 24, 48 and 72 h), and each cell was sequenced to obtain the final mRNA reads. As a result,  each of the four time points contains about 90 cells. Due to possible variations in the speed of differetiation across the cells, we expect the cells to be approximately uniform-distributed along the progression path, whose order may be recovered by the matrix reordering algorithm. The raw count data were preprocessed and normalized using the functions \texttt{CreateSeuratObject} and \texttt{NormalizeData} in the R package \texttt{Seurat}\footnote{https://cran.r-project.org/web/packages/Seurat/index.html} under default settings. We applied the functions \texttt{FindVariableFeatures} and \texttt{ScaleData} in \texttt{Seurat} to identify and standardize the levels of $p=1500$ most variable genes for subsequent analysis. Then  a pairwise similarity matrix  across the 372 cells was calculated as $\bY=c\bJ_n-\bD$, where  $\bJ_n$ is an all-one matrix,  $c$ is a sufficiently large constant making $\bY$ nonnegative (e.g., $c=\|\bD\|_\infty$), and $\bD$ is the pairwise Euclidean distance matrix of the $r$-dimensional spectral embeddings of the cells, based on the singular value decomposition of the standardized data. Specifically, after obtaining the standarized dataset, which contains  expression levels of $p$ genes for $n$ cells, we apply SVD to the data matrix, say $X\in\R^{n\times p}$, and define the $r$-dimensional embedding as the leading $r$ left singular vectors of $X$ weighted by their associated singular values.  Finally, we applied the adaptive sorting (AS) and the spectral seriation (SS) to reorder the similarity matrix $\bY$, which gives the inferred temporal order of cells, up to a possible reversion. To evaluate the performance of two methods, we compared the inferred temporal orders with  the true order using Spearman's rho statistic. On the left of Figure \ref{comp.fig}, we show a boxplot of Spearman's rho statistics evaluated over various values for $r\in\{2,3,...,20\}$. Our AS algorithm shows clear advantages over SS in terms of the preciseness of the inferred temporal orders, even though  the thus constructed similarity matrices have possibly dependent entries. In particular, the evaluation results were consistent for different choices of $p\in\{500,1500,2500\}$. 

The second dataset consists of single-cell RNA sequencing  reads for 149 human primordial germ cells ranging from 4 weeks  to 19 weeks old \citep{guo2015transcriptome}. Specifically, there were 6 cells of 4 weeks old, 37 cells of 7 weeks old, 20 cells of 10 weeks old, 27 cells of 11 weeks old, and 57 cells of 19 weeks old. The RNA-seq data were preprocessed and normalized using the same procedure,  leading to a similarity matrix across the 149 cells calculated from the $r$-dimensional spectral embeddings of the cells. On the right of Figure \ref{comp.fig}, we have a boxplot for the Spearman's rho statistics evaluated for different $r$ values ($r\in\{2,3,...,20\}$), again indicating AS to be overall much better than SS for inferring temporal orders of single cells. Like the previous example, we also found the results to be consistent over different choices of $p\in\{500,1500,2500\}$.

\section{Discussion} \label{dis.sec}

This paper studied the matrix reordering problem for a large class of monotone Toeplitz matrices characterized by $\calT_n$, where the diagonals decay as they deviate from the main diagonal. The analysis can be easily extended to the class of monotone Toeplitz matrices where the diagonals increase as they deviate from the main diagonal (i.e., $\theta_1\le \theta_2\le ...\le \theta_{n-1}$). Specifically, on the one hand, since the direction of monotonicity is not essential in both the lower and the upper bound arguments, the analysis of the fundamental limit of matrix reordering as presented in Section \ref{rate.sec} may be adopted with minor changes to obtain  similar results. On the other hand, a computationally efficient algorithm can be constructed by slightly modifying the adaptive sorting algorithm in Section \ref{ipdd.sec} with $\widetilde\pi(1)$ in the initialization step given by $\widetilde\pi(1)=\argmax_{i\in[n]}S_i$. 

In Section \ref{mle.sec}, we essentially reduced the problem of matrix reordering to an estimation problem, about which we obtained sharp upper bound for the estimation risk. In particular, our proofs of Theorem \ref{mle.thm} and Proposition \ref{red.prop} together yielded the risk upper bound
\begin{align} \label{risk1}
\sup_{\calT'_n\times \calS'_n \subseteq\calT_n\times \calS_n}\sup_{(\bTheta^*,\Pi^*)\in \calT'_n\times \calS'_n}&\E \|\widehat{\Pi}^{lse}{\bTheta}^* (\widehat{\Pi}^{lse})^\top-\Pi^*\bTheta^*(\Pi^*)^\top\|_F\nonumber\\
&\le C\sigma\sqrt{n\log n}.
\end{align}
The above rate of convergence is in fact minimax optimal as one can show using a similar argument as in the proof of Theorem \ref{low.bnd.thm} that
\begin{align} \label{risk2}
\inf_{\widehat\Pi}\sup_{\calT'_n\times \calS'_n\subseteq \calT_n\times \calS_n}\sup_{(\bTheta^*,\Pi^*)\in \calT'_n\times \calS'_n}&\E \|\widehat{\Pi}{\bTheta}^*\widehat{\Pi}^\top-\Pi^*\bTheta^*(\Pi^*)^\top\|_F\nonumber\\
&\ge C'\sigma\sqrt{n\log n}.
\end{align}
In particular, to prove (\ref{risk1}) there is no requirement on the minimum separation, unlike Theorem \ref{mle.thm} which requires $\rho^*\gtrsim \sigma\sqrt{n\log n}$. Nevertheless, Equations (\ref{risk1}) and (\ref{risk2}) imply a minimax rate of convergence of order $\sigma\sqrt{n\log n}$, indicating a close connection between the exact recovery problem under the loss $\tau_\bTheta$ and the estimation problem under the loss $|\hat\Pi\bTheta\hat\Pi^\top-\Pi\bTheta\Pi^\top|_F$. However, our analyses of the spectral seriation algorithm and the adaptive sorting algorithm cannot be directly adapted to obtain an estimation result.

In addition to monotone Toeplitz matrices, it is important to extend the current theoretical framework to other matrix classes that are relevant in practice, such as bimonotone matrices \citep{mao2020towards}, circulant matrices \citep{issartel2021optimal}, Robinson matrices \citep{atkins1998spectral}, and Monge matrices \citep{hutter2020estimation}. Furthermore, compared to the exact recovery problem considered in this study, evaluating performance under partial recovery criteria \citep{ma2021optimala} can be less restrictive and potentially more interesting, albeit more challenging. We plan to investigate this problem systematically in a future study.

This paper considered the prototypical setting where the noise matrix has independent, homoscedastic, and sub-Gaussian entries. However, noise structures in real-world applications can be more complex. It is intriguing to consider matrix reordering in settings with dependent, heteroscedastic, and/or heavy-tailed observations. Furthermore, in terms of computationally efficient algorithms, in the absence of a computational lower bound, there may exist other algorithms with weaker separation requirements than the proposed AS algorithm. Although we are currently unaware of such algorithms, we recognize this as an interesting follow-up question and plan to explore it in future work. Some of the technical tools and theoretical results developed in this paper may be useful for solving the problem in more intricate settings.

\section{Proof of Main Results} \label{proof.sec}

In this section, we present  the proofs of Proposition \ref{red.prop} and Theorems \ref{mle.thm}, \ref{low.bnd.thm}, \ref{ss.thm} and \ref{eff.thm}. The proofs of other theorems, propositions and technical lemmas are given in  Section \ref{supp.sec}.

\subsection{Proof of Proposition \ref{red.prop}} \label{prop1.sec}

Note that for any $\bTheta\in\calT_n$, the simple inequality
\beq \label{simple.ineq}
1\{\widehat\Pi \bTheta \widehat\Pi^\top\ne \Pi\bTheta\Pi^\top\} \le e^{{\|\widehat\Pi \bTheta \widehat\Pi^\top-\Pi\bTheta\Pi^\top\|_F}-{\rho(\bTheta,\mathcal{S}'_n)}},
\eeq
holds trivially for
\[
\rho(\bTheta,\mathcal{S}'_n)=\min_{\substack{\Pi_1,\Pi_2\in\mathcal{S}'_n\\\Pi_1 \bTheta \Pi_1^\top\ne \Pi_2\bTheta\Pi_2^\top}}\|\Pi_1 \bTheta \Pi_1^\top-\Pi_2\bTheta\Pi_2^\top\|_F.
\]
On the other hand, we note that for any  $\bTheta_1,\bTheta_2 \in \calT_n$, and any permutation matrices $\Pi_1,\Pi_2\in\mathcal{S}'_n$, it holds that
\beq \label{P}
\|\bTheta_1-\bTheta_2\|_F\le \|\Pi_1\bTheta_1\Pi_1^\top-\Pi_2\bTheta_2\Pi_2^\top\|_F.
\eeq
To see this, note that by the Toeplitz structure of $\bTheta_1$ and $\bTheta_2$, one can identify $\bTheta_1$ and $\bTheta_2$ with two monotonic nondecreasing vectors of dimension $n^2$, denoted as $\text{vec}(\bTheta_1)$ and $\text{vec}(\bTheta_2)$, respectively, by arranging the matrix entries in the following order: the main diagonals, the first principal diagonals, the second principal diagonals, etc. As a consequence, for any $\Pi_1,\Pi_2\in\calS'_n$,
\begin{align*}
&	\|\bTheta_1-\bTheta_2\|_F^2\\
	&=\|\text{vec}(\bTheta_1)-\text{vec}(\bTheta_2)\|_2^2\\
	&=\|\text{vec}(\bTheta_1)\|_2^2+\|\text{vec}(\bTheta_2)\|_2^2-2\text{vec}(\bTheta_1)^\top\text{vec}(\bTheta_2)\\
	&\le \|\text{vec}(\Pi_1\bTheta_1\Pi_1^\top)\|_2^2+\|\text{vec}(\Pi_2\bTheta_2\Pi_2^\top)\|_2^2\\
	&\quad-2\text{vec}(\Pi_1\bTheta_1\Pi_1^\top)^\top\text{vec}(\Pi_2\bTheta_2\Pi_2^\top)\\
	&=\|\text{vec}(\Pi_1\bTheta_1\Pi_1^\top)-\text{vec}(\Pi_2\bTheta_2\Pi_2^\top)\|_F^2\\
	&=\|\Pi_1\bTheta_1\Pi_1^\top-\Pi_2\bTheta_2\Pi_2^\top\|^2_F,
\end{align*}
where the third line follows from the permutation invariance of the $\ell_2$ norm and the following elementary rearragement inequality.

\bel [Rearrangement Inequality] \label{re.ineq}
If $a_1\ge a_2\ge ...\ge a_n$ and $b_1\ge b_2\ge ... \ge b_n$, then 
$a_nb_1+...+a_1b_n \le a_{\sigma(1)}b_1+...+a_{\sigma(n)}b_n \le a_1b_1+...+a_nb_n,$
where $\sigma$ is any permutation in $\calS_n$.
\eel

Thus for any $\bTheta\in\calT_n$  and $\Pi, \widehat\Pi\in \mathcal{S}'_n$, we have
\begin{align*}
&	\|\widehat{\Pi}\bTheta \widehat{\Pi}^\top-\Pi\bTheta\Pi^\top\|_F\\
	&\le \|\widehat{\Pi}\widehat{\bTheta} \widehat{\Pi}^\top-\Pi\bTheta\Pi^\top\|_F+\|\widehat{\Pi}\widehat{\bTheta} \widehat{\Pi}^\top-\widehat{\Pi}\bTheta \widehat{\Pi}^\top\|_F\\
	&= \|\widehat{\Pi}\widehat{\bTheta} \widehat{\Pi}^\top-\Pi\bTheta\Pi^\top\|_F+\|\widehat{\bTheta}-\bTheta\|_F\\
	&\le 2\|\widehat{\Pi}\widehat{\bTheta} \widehat{\Pi}^\top-\Pi\bTheta\Pi^\top\|_F,
\end{align*}
where the last inequality follows from the property (\ref{P}). Combining the above inequalities and taking expectations for both sides, we obtain the final inequality 
\beq 
P(\widehat\Pi \bTheta \widehat\Pi^\top\ne \Pi\bTheta\Pi^\top)\le e^{{2\E \|\widehat{\Pi}\widehat{\bTheta} \widehat{\Pi}^\top-\Pi\bTheta\Pi^\top\|_F}-{\rho(\bTheta,\mathcal{S}'_n)}}.
\eeq

\subsection{Proof of Theorem \ref{mle.thm}} \label{mle.thm.sec}

In the following, to simplify notation, we write $(\widehat{\bTheta}^{lse}, \widehat{\Pi}^{lse})$ as $(\widehat{\bTheta}, \widehat{\Pi})$, and denote the permutation maps associated to the permutation matrices $\widehat\Pi$ and $\Pi^*$ as $\widehat{\pi}:[n]\to[n]$ and $\pi^*:[n]\to[n]$, respectively. We also identify a permutation matrix with the associated permutation map when there is no confusion.
We will prove the following pointwise result, from which the uniform statement in Theorem \ref{mle.thm} follows directly.

\bet[Pointwise guarantee] \label{mle.thm.T}
Under model (\ref{R.model}), there exist some absolute constants $C,c>0$ such that, for sufficiently large $n$, for any $\bTheta^*\in\calT'_n\subseteq\calT_n$ and any $\mathcal{S}'_n\subseteq \mathcal{S}_n$  satisfying $\rho(\bTheta^*;\mathcal{S}'_n)\ge  C\sigma\sqrt{n\log n}$,  the constrained LSE $(\widehat{\Pi},\widehat{\bTheta})$ defined over  $(\calT'_n,\mathcal{S}'_n)$ satisfies
\beq
P_{\bTheta^*,\Pi^*}(\widehat\Pi\bTheta^* \widehat\Pi^\top\ne \Pi^*\bTheta^*(\Pi^*)^\top)\le \exp\{-c\sigma\sqrt{n\log n}\},
\eeq
for each $\Pi^*\in \calS'_n$.
\eet

To prove Theorem \ref{mle.thm.T}, by the general reduction scheme (Proposition \ref{red.prop}), it suffices  to obtain the upper bound for the matrix denoising risk $\E \|\widehat{\Pi}\widehat{\bTheta} \widehat{\Pi}^\top-\Pi^*\bTheta^*(\Pi^*)^\top\|_F\le C\sigma\sqrt{n\log n}$ under fixed parameters $(\bTheta^*,\Pi^*)$. The key ingredient is  Chatterjee's variational formula, originally developed in \cite{chatterjee2014new}. 
The following version is proved as Lemma A.1 in \cite{flammarion2019optimal}.

\bel[Chatterjee's variational formula] \label{var.lem}
Let $\mathcal{C}$ be a closed subset of $\R^d$. Suppose $y=a^*+z$ where $a^*\in\mathcal{C}$ and $z\in\R^d$. Let $\hat{a}\in \argmin_{a\in \mathcal{C}}\|y-a\|_2^2$ be a projection of $y$ onto $\mathcal{C}$. Define the function $f_{a^*}:\R_+\to \R$ by
\[
f_{a^*}(t)=\sup_{a\in \mathcal{C}\cap \mathbb{B}_d(a^*,t)}\langle a-a^*,z\rangle -\frac{t^2}{2}.
\]
Then we have
\[
\|\hat{a}-a^*\|_2 \in \argmin_{t\ge 0}f_{a^*}(t).
\]
Moreover, if there exists $t^*>0$ such that $f_{a^*}(t)<0$ for all $t\ge t^*$, then $\|\hat{a}-a^*\|_2\le t^*$.
\eel

We set $a^*=\Pi^*\bTheta^*(\Pi^*)^\top$, $z=\bZ$, $y=\bY$, $\mathcal{C}=\sp(\calT'_n,\mathcal{S}'_n)\equiv\{\bQ\in\R^{n\times n}: \bQ=\Pi\bTheta\Pi^\top, \bTheta\in\calT'_n,\Pi\in\mathcal{S}'_n\}$, $\mathbb{B}(\bTheta_0,t)=\{\bTheta\in\R^{n\times n}: \|\bTheta-\bTheta_0\|_F\le t\}$, and let $\langle \cdot,\cdot\rangle$ be the Hilbert-Schmidt inner product over symmetric matrices in $\R^{n\times n}$. Hence, by Lemma \ref{var.lem}, we have
\begin{align*}
&\|\widehat{\Pi}\widehat{\bTheta}\widehat{\Pi}^\top - \Pi^*\bTheta^*(\Pi^*)^\top\|_F\in \\
&\argmin_{t\ge 0} \bigg\{\sup_{\substack{ \bQ\in\sp(\calT'_n,\mathcal{S}'_n)\cap\\ \mathbb{B}(\Pi^*\bTheta^*(\Pi^*)^\top,t) }}\langle \bQ-\Pi^*\bTheta^*(\Pi^*)^\top,\bZ\rangle -\frac{t^2}{2}\bigg\}.
\end{align*}
Now since $(\Pi^*)^\top\Pi^*={\textup{Id}}_n$, we have
\begin{align*}
	&\quad\sup_{\bQ\in\sp(\calT'_n,\mathcal{S}'_n)\cap \mathbb{B}(\Pi^*\bTheta^*(\Pi^*)^\top,t)}\langle \bQ-\Pi^*\bTheta^*(\Pi^*)^\top,\bZ\rangle\\ &=t\cdot\sup_{\substack{\bQ\in\sp(\calT'_n,\mathcal{S}'_n)\cap\\ \mathbb{B}(t^{-1}\Pi^*\bTheta^*(\Pi^*)^\top,1)}}\langle \bQ-t^{-1}\Pi^*\bTheta^*(\Pi^*)^\top,\bZ\rangle \\
	&=t\cdot\sup_{\substack{\bQ\in\sp(\calT'_n,\mathcal{S}'_n)\cap\\ \mathbb{B}(t^{-1}\Pi^*\bTheta^*(\Pi^*)^\top,1)}}\langle (\Pi^*)^\top\bQ\Pi^*-t^{-1}\bTheta^*,(\Pi^*)^\top\bZ\Pi^*\rangle\\
	&=t\cdot\sup_{\substack{\bQ\in\sp(\calT'_n,(\pi^*)^{-1}\circ\mathcal{S}'_n)\cap\\ \mathbb{B}(t^{-1}\bTheta^*,1)}}\langle\bQ-t^{-1}\bTheta^*,(\Pi^*)^\top\bZ\Pi^*\rangle\\
	&\le t\cdot\sup_{\bH\in \mathcal{M}_n}\langle \bH,\bZ'\rangle,
\end{align*}
where $\mathcal{M}_n\equiv\mathcal{M}(\calT'_n,\mathcal{S}'_n)=\{\bQ-t^{-1}\bTheta^*: \bQ\in \sp(\calT'_n,(\pi^*)^{-1}\circ\mathcal{S}'_n)\cap \mathbb{B}(t^{-1}\bTheta^*,1)\}$, $(\pi^*)^{-1}\circ \mathcal{S}'_n=\{(\Pi^*)^\top \Pi: \Pi\in \mathcal{S}'_n\}$ and $\bZ'=(\Pi^*)^\top\bZ\Pi^*$.  In particular, we have $\mathcal{M}_n\subseteq \mathbb{B}({\bf 0},1)$ and ${\bf0}\in \mathcal{M}_n$. 

By the second statement of Lemma \ref{var.lem}, we can choose $t^*=2\sup_{\bH\in \mathcal{M}_n}|\langle \bH,\bZ'\rangle|+s$ for any constant $s>0$. Then, it can be easily checked that, for any $t\ge t^*$, we have 
\begin{align}
&\sup_{\substack{\bQ\in\sp(\calT'_n,\mathcal{S}'_n)\cap\\ \mathbb{B}(\Pi^*\bTheta^*(\Pi^*)^\top,t)}}\langle \bQ-\Pi^*\bTheta^*(\Pi^*)^\top,\bZ\rangle -\frac{t^2}{2}\nonumber\\
&\le  t\sup_{\bH\in \mathcal{M}_n}\langle \bH,\bZ'\rangle-\frac{t^2}{2}<0.
\end{align}
Thus, we have
\beq \label{re1}
\E \|\widehat{\Pi}\widehat{\bTheta}\widehat{\Pi}^\top-\Pi^*\bTheta^*(\Pi^*)^\top\|_F\le 2\E\sup_{\bH\in \mathcal{M}_n}|\langle \bH,\bZ'\rangle|+s.
\eeq
To control the right-hand side of the above inequality, we use the following version of Dudley's integral inequality, whose proof is given in Section \ref{supp.sec} below.

\bel[Dudley's integral inequality] \label{dudley.lem}
Let $\{X_t\}_{t\in T}$ be a mean zero random process on a metric space $(T,d)$ with sub-Gaussian increments, in the sense that there exists some constant $K>0$ such that $\inf\{t>0: \E \exp((X_t-X_s)^2/t^2)\le 2\} \le Kd(t,s)$ for all $t,s\in T$. Then, we have
\beq
\E \sup_{t\in T} |X_t |\le CK\int_0^{\textup{diam}(T)}\sqrt{\log \mathcal{N}(T, d, \epsilon)}d\epsilon+\inf_{t\in T}\E |X_{t}|,
\eeq
where $\mathcal{N}(T, d, \epsilon)$ is the $\epsilon$-covering number of $T$, that is, the smallest number of closed balls with centers in $T$ and radii $\epsilon$ whose union covers $T$.
\eel 

Since $\big\{\langle \bH,\bZ'\rangle\big\}_{\bH\in \mathcal{M}_n}$ is a mean-zero randon process with sub-Gaussian increment, and one can check that for any $\bH_1,\bH_2\in \mathcal{M}_n$, there exists some constant $C'>0$ such that
\[
\inf\{t>0: \E \exp(\langle \bH_1-\bH_2,\bZ'\rangle^2/t^2)\le 2\} \le \sigma C'\|\bH_1-\bH_2\|_F.
\]
From Lemma \ref{dudley.lem}, it follows that
\begin{align} \label{re2}
&\E \sup_{\bH\in \mathcal{M}_n}|\langle \bH,\bZ'\rangle|\nonumber\\
&\le C\sigma\int_0^{\text{diam}(\mathcal{M}_n)}\sqrt{\log \mathcal{N}(\mathcal{M}_n, d_2, \epsilon)}d\epsilon+\E |\langle \bH_0,\bZ'\rangle|,
\end{align}
for any $\bH_0\in \mathcal{M}_n$,
where the metric $d_2$ is defined as  $d_2(\bG_1,\bG_2)=\|\bG_1-\bG_2\|_F$. In particular, we can take $\bH_0={\bf 0}$ to get 
\beq \label{re3}
\E |\langle \bH_0,\bZ'\rangle|=0.
\eeq 
In this way, we further reduced the calculation of the denoising risk $\E \|\widehat{\Pi}\widehat{\bTheta} \widehat{\Pi}^\top-\Pi\bTheta\Pi^\top\|_F$ to that of the metric entropy of the set $\mathcal{M}_n$ under the Frobenius norm. The following lemma, proved in Section \ref{supp.sec} below,  provides an  estimate of such a entropy measure.

\bel \label{entropy.lem}
Under the conditions of Theorem \ref{mle.thm}, using the above notations, we have
\beq
\int_0^{\textup{diam}(\mathcal{M}_n)}\sqrt{\log \mathcal{N}(\mathcal{M}_n, d_2, \epsilon)}d\epsilon\le C\sqrt{n\log n}.
\eeq
\eel

The above upper bound estimate of the entropy integral along with inequalities (\ref{re1})  (\ref{re2}) and (\ref{re3}) leads to the upper bound in Theorem \ref{mle.thm.T}, as along as we choose $s$ sufficiently small, for example, $s=\sigma$.
This completes the proof of the theorem.

\subsection{Proof of Theorem \ref{low.bnd.thm}} \label{low.bnd.thm.sec}

From Lemma \ref{deza}, we can further simplify the lower bound for $\beta_n(d)$ using the inequality
\[
\bigg(\frac{n}{e}\bigg)^n e \le n!\le    \bigg(\frac{n+1}{e}\bigg)^{n+1}e,
\]
which implies
\begin{align*}
&\frac{n!}{(n-d)!}\ge \frac{(n/e)^n}{[(n-d+1)/e]^{n-d+1}}\ge \\
&\bigg( \frac{n}{n-d+1}\bigg)^{n-d+1} (n/e)^{d-1}\ge  (n/e)^{d-1}.
\end{align*}
Consequently, for $2\le d\le n-2$, we have
\beq
\log \beta_n(d)\ge  (d-1)\log (n/e)-\log (n-d).
\eeq
Set $d=(1-e^{-1})n$. We have
\[
\log \beta_n(d)\ge  (d-2)\log (n/e).
\]
Thus, for $n\ge 16$, we have
\beq \label{packing.lb}
\log \beta_n(d)\ge  \frac{n}{2}\log (n/e).
\eeq
This is the key packing number inequality that we will use for the proof of the fundamental limit.

The proof of the fundamental limit starts with a careful construction of a set of least favourable scenarios, over which exact matrix reorderinig is most difficult to achieve.
To this end, we consider the following tridiagonal signal matrix, which is also the adjacency matrix of a Hamiltonian path,
\beq \label{trid}
\bTheta_0 =  \begin{bmatrix}
	0& \delta &0 &... &0\\
	\delta & 0 &\delta&...&0\\
	0 &\delta&0&...&0\\
	\vdots    &       &       & \ddots &\\
	0 & 0 & 0&...& 0      
\end{bmatrix}.
\eeq
We also  consider the subset of $\mathcal{S}_n$ that only permutes the columns and rows in $\bTheta_0$ satisfying $i=3k-1$ for some $k\in\{1,..., \lfloor n/3 \rfloor \}$. The reason we consider such a class of permutations is that, by treating $\bTheta_0$ as a concatenation of small blocks of size $3\times 3$, it will be seen that, for any two permutations $\pi_1,\pi_2$ in the above subset, the distance $\|\Pi_1\bTheta_0\Pi_1^\top-\Pi_2\bTheta_0\Pi_2^\top\|_F^2$ is completely determined by the Hamming distance $d_H(\pi_1,\pi_2)$ between the two permutations (equation (\ref{dist.S}) below).

Let $n'=\lfloor n/3 \rfloor$. Suppose $n\ge 48$, or $n'\ge 16$, and set $d=(1-e^{-1})n'$. By  inequality (\ref{packing.lb}), there exists a subset $\mathcal{S}_{n'}^*$ of $\mathcal{S}_{n'}$ that attains the maximal  $d$-packing number $\beta_{n'}(d)$, where $\log \beta_{n'}(d)\ge \frac{n'}{2}\log (n'/e)\ge\frac{5n}{32}\log (5n/16e)$, as $n'\ge 5n/16$ for $n\ge 48$.

Now we identify the permutations in $\mathcal{S}_{n'}^*$ with the permutations in $\mathcal{S}_n$ that only involve the $i$-th element for $i=3k-1, k\in\{1,...,\lfloor n/3 \rfloor \}$. 
Thus, we have constructed a subset  $\mathcal{S}_n^*$ of  $\mathcal{S}_n$ where $\log | \mathcal{S}_n^*|\ge \frac{5n}{32}\log (5n/16e)$ elements with mutual distance at least $d\ge (1-e^{-1})5n/16$. 

Now note that for any two permutations $\Pi_1,\Pi_2\in\mathcal{S}_n^*$ such that $d_H(\pi_1,\pi_2)=d$, we have
\beq \label{dist.S}
2d\delta^2\le\|\Pi_1\bTheta_0\Pi_1^\top -\Pi_2\bTheta_0\Pi_2^\top\|_F^2\le 4n\delta^2
\eeq
To see this, we notice that
\[
\|\Pi_1\bTheta_0\Pi_1^\top -\Pi_2\bTheta_0\Pi_2^\top\|_F^2=\|\Pi\bTheta_0\Pi^\top-\bTheta_0\|_F^2,
\]
where $\Pi=\Pi_2^\top\Pi_1\in \mathcal{S}_n^*$.
We denote $\bTheta_0=(\theta_{ij})_{1\le i,j\le n}$ and calculate that, for any $i$ such that $\pi(i)>i$, we have
\begin{align*}
	&\sum_{j=1}^n(\theta_{\pi(i),\pi(j)}-\theta_{ij})^2 \\
	&\ge (\theta_{\pi(i),\pi(i)+1}-\theta_{i,\pi(i)+1})^2+(\theta_{\pi(i),i-1}-\theta_{i,i-1})^2 \\
	&\ge 2\delta^2,
\end{align*}
since the $(\pi(i)+1)$-th and the $(i-1)$-th components are fixed in the permutation $\Pi$. A similar result can be obtained for $i$ such that $\pi(i)<i$. Therefore, we obtain the lower bound
\begin{align}\label{lower}
&\|\Pi_1\bTheta_0\Pi_1^\top -\Pi_2\bTheta_0\Pi_2^\top\|_F^2\ge \sum_{i: i\ne \pi(i)}\sum_{j=1}^n(\theta_{\pi(i)\pi(j)}-\theta_{ij})^2\nonumber\\
&\ge 2d\delta^2.
\end{align}
The upper bound follows from the simple inequality
\beq \label{upper}
\|\Pi_1\bTheta_0\Pi_1^\top -\Pi_2\bTheta_0\Pi_2^\top\|_F^2\le2\|\bTheta_0\|_F^2\le 4(n-1)\delta^2
\eeq
In the following, we consider the fundamental limit for $\rho^*(\cdot,\cdot)$ over the subspace $\{\bTheta_0\}\times \mathcal{S}^*_n$.
The proof relies on the following lemma  from \cite{tsybakov2009introduction}.

\bel \label{low.bnd.lem}
Assume that for some integer $M\ge 2$ there exist distinct parameters $\theta_0,...,\theta_M$ from the parameter space $\Theta$ and mutually absolutely continuous probability measures $P_0,...,P_M$ with $P_j=P_{\theta_j}$ for $j=0,1,...,M$, defined on a common probability space $(\Omega, \mathcal{F})$ such that the averaged KL divergence $\frac{1}{M}\sum_{j=1}^M D(P_j,P_0) \le \frac{1}{8}\log M.$ Then, for every measurable mapping $\hat{\theta}: \Omega\to \Theta$,
\[
\max_{j=0,...,M}P_j(\hat{\theta}\ne \theta_j)\ge \frac{\sqrt{M}}{\sqrt{M}+1}\bigg( \frac{3}{4}-\frac{1}{2\sqrt{\log M}}\bigg).
\]
\eel

Applying the above lemma to the parameter subspace $\{\bTheta_0\}\times \mathcal{S}^*_n$, we could check that for any $\Pi_1,\Pi_2\in \mathcal{S}^*_n$, the KL divergence between the probability measures of $\bY_1=\Pi_1\bTheta_0\Pi_1^\top+\bZ$ and that of $\bY_2=\Pi_2\bTheta_0\Pi_2^\top+\bZ$ can be bounded by
\[
D(P_{\bY_1},P_{\bY_2})=\frac{ \|\Pi_1\bTheta_0\Pi_1^\top-\Pi_2\bTheta_0\Pi_2^\top\|_F^2}{2\sigma^2}\le \frac{2n\delta^2}{\sigma^2},
\]
where the last inequality follows from (\ref{dist.S}).
Set $\delta = 0.06\sigma\sqrt{\log (5n/16e)}$. By the lower bound on $|\mathcal{S}_n^*|$, we have
\[
D(P_{\bY_1},P_{\bY_2})\le 0.019n\log (5n/16e) \le\frac{1}{8}\log |\mathcal{S}^*_n|.
\]
Thus, by Lemma \ref{low.bnd.lem}, it follows that, for any $n\ge 48$,
\begin{align}
&\inf_{\hat{\Pi}}\max_{(\Theta,\Pi)\in \{\bTheta_0\}\times \mathcal{S}^*_n}P(\hat\Pi\bTheta\hat\Pi^\top\ne\Pi\bTheta\Pi^\top)\nonumber\\
&=\inf_{\hat{\pi}}\max_{(\Theta,\pi)\in \{\bTheta_0\}\times \mathcal{S}^*_n}P(\hat{\pi}\ne \pi)\ge 0.6.
\end{align}
Lastly,  by the lower bound (\ref{dist.S}) on the mutual distance between the elements in the set $\mathcal{S}_n^*$, we have
\begin{align*}
	&\quad\rho^*(\{\bTheta_0\}, \mathcal{S}_n^*) = \min_{\substack{\Pi_1,\Pi_2\in \mathcal{S}_n^*\\\Pi_1\ne \Pi_2}}\|\Pi_1\bTheta_0\Pi_1-\Pi_2\bTheta_0\Pi_2\|_F\\
	&\ge \sqrt{2d\delta^2}>0.037\sigma\sqrt{n\log (5n/16e)}>0.02\sigma \sqrt{n\log n},
\end{align*}
where the last inequality follows from $\log(5n/16e)\ge 0.44\log n$.

In other words, we have found a subset $\{\bTheta_0\}\times \mathcal{S}_n^*$ of $\calT_n\times \mathcal{S}_n$ with $\rho^*(\{\bTheta_0\}, \mathcal{S}_n^*) \ge 0.02\sigma \sqrt{n\log n}$ such that uniform exact matrix ordering over $\{\bTheta_0\}\times \mathcal{S}_n^*$ with high probability is not possible. 
In particular, from the above argument, we can see that for any $\rho\le 0.02\sigma \sqrt{n\log n}$, we can always construct similar subsets by choosing $\delta$ smaller such that $\rho^*(\{\bTheta_0\}, \mathcal{S}_n^*)=\rho$, and show that uniform exact matrix ordering over $\{\bTheta_0\}\times \mathcal{S}_n^*$ with high probability is not possible. This completes the proof of the theorem.

\subsection{Proof of Theorem \ref{ss.thm}} \label{ss.thm.sec}

We first define the degree and the Laplacian operators as follows.

\begin{definition}[Degree \& Laplacian Operators] \label{lap.def}
	For a symmetric matrix ${\bf A}=(a_{ij})_{1\le i,j\le n}\in \R^{n\times n}$, define the degree operator $\bD: \R^{n\times n}\to \R^{n\times n}$ by $\bD({\bf A})=\textup{diag}(\sum_{i=1}^n a_{i1}, ..., \sum_{i=1}^n a_{in})$, and the Laplacian operator $\bL: \R^{n\times n}\to \R^{n\times n}$  by $\bL({\bf A})=\bD({\bf A} )-{\bf A}$.
\end{definition}

Without loss of generality, we set $\Pi=\text{Id}_n$. We define $\calT'_n=\{\bTheta_0\}$ where $\bTheta_0$ is the tridiagonal matrix defined in (\ref{trid}), with $\delta=C_0\sigma n^2\sqrt{n}$, and let $\calS'_n$ in the theorem be the permutation set $\calS_n^*$ constructed in the proof of Theorem 2. By inequalities (\ref{lower}) and (\ref{upper}), we have $\rho^*(\calT'_n,\calS'_n)=C\sigma n^3$ for some absolute constant $C>0$.
Note that
$
\bL \equiv\bL(\bY)= \bL(\bTheta)+\bL(\bZ).
$
We denote $\bB=\bL(\bTheta)$ and $\bE=\bL(\bZ)$. Let $\bB=\sum_{i=1}^n\lambda_i\bu_i\bu_i^\top$ be the eigendecomposition of $\bB$, with $0=\lambda_1\le \lambda_2\le ...\le \lambda_n$. In particular, by Lemma \ref{trid.spe.lem} below,  $\bB$ has simple eigenvalues, which implies $\bu_2=\bv$, the Fiedler vector, up to a change of sign. The following lemma is well-known, and can be found, for example, on page 3234 of \cite{nakatsukasa2013mysteries}.

\bel \label{trid.spe.lem}
The Laplacian matrix of any tridiagonal matrix  in the form (\ref{trid}) has eigenvalues
$
\lambda_k=4\delta\sin^2\big( \frac{\pi (k-1)}{2n}\big),$ for $k=1,...,n,$
and eigenvectors $\bu_k=(u_{k1},u_{k2},...,u_{kn})^\top, k=1,...,n$ where
$
u_{kj}=\frac{1}{\sqrt{n/2}}\cos \big(\frac{\pi(k-1)(j-1/2)}{n} \big),
$
for $j=1,...,n.$
In particular, the Fiedler vector $\bu_2$ is
$
\bu_2=\frac{1}{\sqrt{n/2}}\big(\cos \big(\frac{\pi}{2n} \big),\cos \big(\frac{3\pi}{2n} \big),..., \cos \big(\frac{(2n-1)\pi}{2n} \big) \big)^\top.
$
\eel

In order to show $\Pi\bTheta_0\Pi^\top\ne \check{\Pi}\bTheta_0\check\Pi$, by the definition of $\bTheta_0$, it is equivalent to showing $\Pi\ne \check{\Pi}$, or $\frak{r}(\widehat\bv)\ne \frak{r}(\bv)$. The rest of the proof is devoted to
\beq
\liminf_{n\to \infty}\inf_{\Pi\in  \mathcal{S}'_n}P_{\bTheta_0,\Pi}(\frak{r}(\widehat\bv)\ne \frak{r}(\bv))\ge 1/2.
\eeq
Let $\bB+\bE=\sum_{i=1}^n\hat\lambda_i\hat\bu_i\hat\bu_i^\top$ be the eigendecomposition of  $\bB+\bE$ with $\hat\lambda_1\le \hat\lambda_2\le ...\le \hat\lambda_n$. In other words, we have $\widehat\bv=\hat\bu_2$.
We also define $\widehat\bv_-$ as the eigenvector associated to the second smallest eigenvalue of $\bB-\bE$.

In order to show $\frak{r}(\widehat\bv)\ne \frak{r}(\bv)$, we note that by Lemma \ref{trid.spe.lem}, the minimal  distance between any two consecutive components in $\bv$ is bounded by
\[
\min_{1\le i\ne j\le n} |u_{2i}-u_{2j}| \le  \sqrt{\frac{2}{n}} \cdot\frac{\pi}{n}\cdot\sin \bigg(  \frac{1}{2n}\bigg)\le \frac{\pi}{\sqrt{2}n^{5/2}}.
\]
If we are able to show that  $\|\widehat\bv-\bv\|_\infty> \frac{\pi}{\sqrt{2}n^{5/2}}$, then  it follows that  $\frak{r}(\widehat\bv)\ne \frak{r}(\bv)$. Since  $\|\widehat\bv-\bv\|_\infty\ge \frac{1}{\sqrt{n}}\|\widehat\bv-\bv\|_2$, it suffices to show that $\|\widehat\bv-\bv\|_2>  \frac{\pi}{\sqrt{2}n^{2}}$. To this end, we obtain the following proposition, proved in Section \ref{supp.sec} below.

\bep \label{tri.lem}
Define the  probability event $\mathcal{A}_n=\{\|\bE\bv\|_2\ge C_1\sigma\sqrt{n}, \|\bE\|\le C_2\sigma\sqrt{n\log n}, |\bv^\top\bE\bv|\le C_3\sigma\sqrt{\log n}\}.$
Then we have $\lim_{n\to\infty}P(\mathcal{A}_n)=1$, where $C_1, C_2, C_3>0$ are some universal constants.
In addition, under the event $\mathcal{A}_n$, for any constant $c>0$ we have $\|\widehat\bv-\bv\|_2+\|\widehat\bv_--\bv\|_2\ge c/n^{2}$  for  sufficiently large $n$. 
\eep

Since the distribution of $\bE$ is the same as that of $-\bE$, the distributions of the corresponding eigenvectors $\widehat\bv$ and $\widehat\bv_-$ should also be identical. In other words, we have $P( \|\widehat\bv-\bv\|_2\ge z)= P(\|\widehat\bv_--\bv\|_2\ge z)$ for all $z\in\R$. It then follows that, for any constant $c>0$,
\begin{align*}
	P(\mathcal{A}_n) &= P(\mathcal{A}_n, \|\widehat\bv-\bv\|_2+\|\widehat\bv_--\bv\|_2\ge  c/n^{2})\\
	&\le P\big((\mathcal{A}_n\cap\{\|\widehat\bv-\bv\|_2\ge \frac{c}{2n^{2}}\})\cup\\
	&\quad\qquad(\mathcal{A}_n\cap\{\|\widehat\bv_--\bv\|_2\ge \frac{c}{2n^{2}}\})\big)\\
	&\le P(\|\widehat\bv-\bv\|_2\ge cn^{-2}/2)\\
	&\quad+P(\|\widehat\bv_--\bv\|_2\ge cn^{-2}/2)\\
	&= 2P(\|\widehat\bv-\bv\|_2\ge cn^{-2}/2),
\end{align*}
which implies 
$
\liminf_{n\to\infty}P(\|\widehat\bv-\bv\|_2\ge c'n^{-2})\ge 1/2
$
This completes the proof.

\subsection{Proof of Theorem \ref{eff.thm}} \label{eff.thm.sec}

We define the class of $\lambda$-ridged Topelitz matrices as
\beq
\calT_n^{R}(\lambda)=\bigg\{\bTheta\in \calT_n:\theta_1-\theta_{\lceil n/2\rceil }\ge\lambda \bigg\}.
\eeq
Then we have the following propositions hold.

\bep[Theoretical guarantee over $\lambda$-ridged Topelitz matrices] \label{eff.thm.2}
There exists some absolute constants $C,c>0$ such that, whenever $\lambda\ge C\sigma n$,  we have $\sup_{(\bTheta, \Pi)\in \calT_n^{R}(\lambda)\times \mathcal{S}_n}P_{\Theta,\Pi}(\widetilde\Pi\bTheta \widetilde\Pi^\top\ne \Pi\bTheta\Pi^\top)\le n^{-c}$.
\eep

\bep[Sharp $\lambda$-$\rho^*$ correspondence] \label{p.space}
For any $\calT'_n\subseteq \calT_n$ and any $\mathcal{S}_n'\subseteq \mathcal{S}_n$ such that $\rho^*(\calT'_n,\mathcal{S}'_n)\ge C_1 n\lambda^*$ for some absolute constant $C_1>0$, there exists some absolute constant $C_2>0$ such that  $\calT'_n \subseteq \calT_n^{R}(\lambda)$ for  $\lambda= C_2\lambda^*$. On the other hand, the above characterization is asymptotically sharp in the sense that, for any $\lambda^*$, there exists some $(\calT'_n, \calS'_n)\subseteq \calT_n\times \calS_n$ such that $\calT'_n\subseteq \calT^R(\lambda^*)$ and $\rho^*(\calT'_n, \calS'_n)=Cn\lambda^*$.
\eep

The proofs of these propositions can be found in Section \ref{supp.sec} below.
With the above results, Theorem \ref{eff.thm} is then proved by combining Propositions \ref{eff.thm.2} and  \ref{p.space} with $\lambda^*=\sigma n$. In particular, the sharp $\lambda$-$\rho^*$ correspondence in Proposition \ref{p.space} suggests that the argument for proving Theorem \ref{eff.thm}, and therefore the obtained minimal signal strength condition (\ref{s.c}) for the adaptive sorting algorithm, are asymptotically tight.

\section{Proof of Technical Results} \label{supp.sec}

\subsection{Dudley's Integral Inequality: Proof of Lemma \ref{dudley.lem}} \label{dudley.lem.sec}

The integral inequality stated in the lemma generalizes the results in \cite{vershynin2018high}. Specifically, by Remark 8.1.5 (supremum of increments) of \cite{vershynin2018high}, we have
\[
\E \sup_{t\in T}|X_t-X_{t_0}|\le CK\int_{0}^{\textup{diam}(T)}\sqrt{\log \mathcal{N}(T,d,\epsilon)}d\epsilon,
\]
holds for any fixed $t_0\in T$. Then our lemma follows immediately from the simple inequality $
\E\sup_{t\in T}|X_t|\le \E\sup_{t\in T}|X_t-X_{t_0}|+\E|X_{t_0}|.$

\subsection{Metric Entropy Calculation: Proof of Lemma \ref{entropy.lem}} \label{entropy.lem.sec}

In order to study the metric entropy of the set $\mathcal{M}_n$, we consider the permutation set $\calS'_n$ and denote $\calS'_n=\{\Pi_1,...,\Pi_{|\calS'_n|}\}$.  For any given element $\bH\in\mathcal{M}_n$, there exists an element $\Pi_i\in \calS'_n$ and a Toeplitz matrix $\bQ^*\in \calT'_n$ such that 
\[
\bH = (\Pi^*)^\top \Pi_i\bQ^* \Pi_i^\top\Pi^*- t^{-1}\bTheta^*.
\]
In other words, if we define $\bar\sf_i \equiv \{(\Pi^*)^\top \Pi_i\bQ \Pi_i^\top\Pi^*- t^{-1}\bTheta^*: \bQ\in\calT^*_n\}$ where $\calT^*_n$ is the set of all the $n\times n$ symmetric Toeplitz matrices, then there exists an $i\in \{1,...,|\calS'_n|\}$ such that $\bH\in \bar\sf_i$. Hence, we have
\beq
\mathcal{M}_n\subseteq \bigcup_{1\le i\le |\calS'_n|}\bar\sf_i.
\eeq
Moreover, by definition we also have $	\mathcal{M}_n\subseteq\mathbb{B}({\bf 0}, 1)$, so that
\beq
\mathcal{M}_n\subseteq \bigcup_{1\le i\le |\calS'_n|}\bar\sf_i\cap \mathbb{B}({\bf 0}, 1).
\eeq
By the union bound, we have
\begin{align*}
	&\log \mathcal{N}(\mathcal{M}_n, d_2, \epsilon)\\
	&\le \log\bigg[ \sum_{1\le i\le |\calS'_n|}\mathcal{N}(\bar\sf_i\cap \mathbb{B}({\bf 0},1), d_2, \epsilon)\bigg]\\
	&\le \log |\calS'_n|+\log \bigg[\max_{1\le i\le |\calS'_n|}\mathcal{N}(\bar\sf_i\cap \mathbb{B}({\bf 0},1), d_2, \epsilon)\bigg]\\
	&\le \log (n!)+\log \bigg[\max_{1\le i\le m}\mathcal{N}(\bar\sf_i\cap \mathbb{B}({\bf 0},1), d_2, \epsilon)\bigg].
\end{align*}
For the first term in the last inequality, by Stirling's formula, we have
\beq \label{p1}
\log (n!)\le (n+1)\log \frac{n+1}{e}+1\le 3n\log n,
\eeq
where the last inequality holds for all $n\ge 2$. 
In the following, we control the second term in the last inequality, to obtain an upper bound for $\log \mathcal{N}(\mathcal{M}_n, d_2, \epsilon)$.

We control the metric entropy $\log \mathcal{N}(\bar\sf_i\cap \mathbb{B}({\bf 0},1), d_2, \epsilon)$ for any $i\in\{1,2,...,|\calS'_n|\}$. Suppose without loss of generality we consider the set $\bar\sf_{i_0}$, associated to the permutation $\Pi_{i_0}\in \calS'_n$. Recall that we denote $\pi_{i_0}:[n]\to[n]$ as the permutation map associated to the permutation matrix $\Pi_{i_0}$. Now we consider the map $\Phi_{i_0}: (\R^{n}, d_2)\to (\bar\sf_i, d_2)$, where for each $\gamma=(\gamma_1,...,\gamma_{n})\in \R^{n}$, we define a matrix $\bH\in\bar\sf_i$ as follows:
\begin{enumerate}
	\item Define a symmetric Toeplitz matrix $\bQ^*(\gamma)$ such that its first row is $\gamma$.
	\item Set $\Phi_{i_0}(\gamma)= (\Pi^*)^\top \Pi_{i_0}\bQ^*(\gamma) \Pi_{i_0}^\top\Pi^*- t^{-1}\bTheta^*$.
\end{enumerate}
The map $\Phi_{i_0}$ plays a key role in translating the metric entropy of $\bar\sf_{i_0}\cap \mathbb{B}({\bf 0},1)$ in $(\R^{n\times n},d_2)$ to that of the Euclidean ball $\mathbb{B}_{n}({\bf 0},1)$ in $(\R^n,d_2)$. 
Specifically, we will need the following lemma concerning the property of covering numbers with respect to Lipschitz maps.

\bel[\cite{szarek1998metric}] \label{szarek.lem}
Let $(M,d)$ and $(M_1,d_1)$ be metric spaces, $K\subset M$, $\Phi:M\to M_1$, and let $L>0$. If $\Phi$ satisfies $d_1(\Phi(x),\Phi(y)) \le Ld(x,y)$ for $x,y,\in M$, then, for every $\epsilon>0$, we have $\mathcal{N}(\Phi(K),d_1,L\epsilon)\le \mathcal{N}(K,d,\epsilon).$
\eel

Now to use Lemma \ref{szarek.lem}, for any $x,y\in \R^{n}$, we have, on the one hand,
\begin{align}
	&\|\Phi_{i_0}(x)-\Phi_{i_0}(y)\|_F^2\nonumber\\
	&=\|(\Pi^*)^\top \Pi_{i_0}\bQ^*(x) \Pi_{i_0}^\top\Pi^*-(\Pi^*)^\top \Pi_{i_0}\bQ^*(y) \Pi_{i_0}^\top\Pi^*\|_F^2\nonumber \\
	&= \|\bQ^*(x) -\bQ^*(y) \|_F^2\nonumber \\
	&\le 2n\|x-y\|_2^2, \label{ub.Phi}
\end{align}
where the last inequality follows from the property of symmetric Toeplitz matrices, and on the other hand
\begin{align}
	&\|\Phi_{i_0}(x)-\Phi_{i_0}(y)\|_F^2\nonumber\\
	&=\|(\Pi^*)^\top \Pi_{i_0}\bQ^*(x) \Pi_{i_0}^\top\Pi^*-(\Pi^*)^\top \Pi_{i_0}\bQ^*(y) \Pi_{i_0}^\top\Pi^*\|_F^2\nonumber \\
	&= \|\bQ^*(x) -\bQ^*(y) \|_F^2\nonumber \\
	&\ge \|x-y\|_2^2, \label{lb.Phi}
\end{align}
where the last inequality follows by only considering the first rows of $\bQ^*(x) -\bQ^*(y) $.
By the above lower bound (\ref{lb.Phi}), it follows that
\beq  \label{sub.sf}
\big\{ \gamma\in\R^{n}: \|\Phi_{i_0}(\gamma)\|_F\le 1 \big\}\subseteq \big\{  \gamma\in \R^{n}: \|\gamma\|_2\le  1 \big\}.
\eeq
In addition, for any $\bH\in \bar\sf_{i_0}\cap \mathbb{B}({\bf 0},1)$, there exists a vector $\gamma\in \big\{ \gamma\in\R^{n}: \|\Phi_{i_0}(\gamma)\|_F\le 1 \big\}$ such that $\Phi_{i_0}(\gamma)=\bH$. Then we also have
\beq \label{inc.sf}
\bar\sf_{i_0}\cap \mathbb{B}({\bf 0},1)\subseteq \Phi_{i_0}(\{ \gamma\in\R^{n}: \|\Phi_{i_0}(\gamma)\|_F\le 1 \}).
\eeq
By taking the images of both sides of (\ref{sub.sf}) under the map $\Phi_{i_0}$ and applying (\ref{inc.sf}), we have
\[
\bar\sf_{i_0}\cap \mathbb{B}({\bf 0},1)\subseteq \Phi_{i_0}(\mathbb{B}_{n}({\bf 0},1)),
\]
which implies
\[
\log \mathcal{N}(\bar\sf_i\cap \mathbb{B}({\bf 0},1), d_2, \epsilon)\le 	\log \mathcal{N}(\Phi_{i_0}(\mathbb{B}_{n}({\bf 0},1)), d_2, \epsilon).
\]
The right-hand side of the inequality can be further bounded by Lemma \ref{szarek.lem} and the upper bound (\ref{ub.Phi}) as
\begin{align} \label{p3}
&\log \mathcal{N}(\Phi_{i_0}(\mathbb{B}_{n}({\bf 0},1)), d_2, \epsilon)\le \log \mathcal{N}(\mathbb{B}_{n}({\bf 0},1), d_2, \epsilon/\sqrt{2n})\nonumber \\
&\le n\log\bigg(  \frac{2\sqrt{2n}}{\epsilon}+1\bigg),
\end{align}
where the last inequality is a direct consequence of the following entropy bound for the Euclidean unit ball 
\[
\mathcal{N}(\mathbb{B}_{n}(0,1), d_2, \epsilon)\le (2/\epsilon+1)^n,
\]
which can be found, for example, in Corollary 4.2.13 (covering numbers of the Euclidean ball) of \cite{vershynin2018high}.
Hence, combining inequalities (\ref{p1})  and (\ref{p3}), we have
\beq
\log \mathcal{N}(\mathcal{M}_n, d_2, \epsilon)\lesssim n\log n+n\log \bigg(  \frac{2\sqrt{2n}}{\epsilon}+1\bigg).
\eeq
Finally, since $\mathcal{M}_n\subseteq \mathbb{B}(0,1)$, we have for $n\ge 2$,
\begin{align*}
	&\int_0^{\text{diam}(\mathcal{M}_n)}\sqrt{\log \mathcal{N}(\mathcal{M}_n, d_2, \epsilon)}d\epsilon\\
	&\le 	C\int_0^{2}\sqrt{n\log n+n\log \bigg(  \frac{2\sqrt{2n}}{\epsilon}+1\bigg)}d\epsilon\\
	&\le  C\sqrt{n\log n}+C\sqrt{n}\int_0^{{2}}\sqrt{\log \bigg(  \frac{2\sqrt{2n}}{\epsilon}+1\bigg)}d\epsilon\\
	&\le C\sqrt{n\log n},
\end{align*}
where the last inequality follows from
\begin{align*}
	&\int_0^{{2}}\sqrt{\log \bigg(  \frac{2\sqrt{2n}}{\epsilon}+1\bigg)}d\epsilon\le \int_0^{{2}}\sqrt{\log \bigg(  \frac{4\sqrt{n}}{\epsilon}\bigg)}d\epsilon\\
	&\le  C\sqrt{\log n}+C\int_0^{{2}}\sqrt{\log \bigg(  \frac{1}{\epsilon}\bigg)}d\epsilon\\
	&\le C\sqrt{\log n}+C\int_{0}^\infty t^2e^{-t^2}dt\le C\sqrt{\log n}.
\end{align*}
This completes the proof of the lemma.

\subsection{Inconsistency of Sample Eigenvectors: Proof of Proposition \ref{tri.lem}}   \label{prop.4}

We define  $\bB-\bE=\sum_{i=1}^n\tilde\lambda_i\tilde\bu_i\tilde\bu_i^\top$ where $\tilde\lambda_1\le \tilde\lambda_2\le ...\le \tilde\lambda_n$ are the eigenvalues, and $\{\tilde\bu_i\}$ are the eigenvectors. 
We denote $\angle (\bg,\bh)$ as the angle between two vectors $\bh$ and $\bg$. Because the eigenvectors are identifiable up to a change of sign, we assume without loss of generality that all the angles are between 0 and $\pi/2$. 
The proof of this proposition is separated into the following steps:
\begin{enumerate}
	\item  We show that under event $\mathcal{A}_n$, if instead $\|\widehat\bv-\bv\|_2\le c/n^{2}$ and $\|\widehat\bv_--\bv\|_2\le c/n^{2}$ for some  constant $c>0$, then the angle $\angle (\widehat\bv,(\bB-\bE)\widehat\bv)\ge \frac{\pi}{2}-\epsilon_1$ for some sufficiently small constant $\epsilon_1\in(0,\pi/2)$.
	\item We show that the above statement implies that $\|\widehat\bv- \widehat\bv_-\|_2\ge \epsilon_2$ for some constant $\epsilon_2>0.$ A contradiction to the first statement by triangle inequality. Thus, we conclude that under event $\mathcal{A}_n$, we must have $\|\widehat\bv-\bv\|_2+\|\widehat\bv_--\bv\|_2\ge c/n^{2}$. 
	\item We finish the proof by showing that  event $\mathcal{A}_n$ holds in probability.
\end{enumerate}
The detailed proofs are presented in order.

\paragraph{Step I} Note that $\hat\lambda_2\widehat\bv=(\bB+\bE)\widehat\bv$.
Then 
\begin{align}
\angle (\widehat\bv,(\bB-\bE)\widehat\bv)&=	\angle ((\bB+\bE)\widehat\bv,(\bB-\bE)\widehat\bv)\nonumber\\
&\ge 	\angle ((\bB+\bE)\widehat\bv,\bB\widehat\bv).
\end{align}
Now consider the triangle with sides $
\bB\widehat\bv, 
\bE\widehat\bv$ and $
(\bB+\bE)\widehat\bv$. 
By the sine rule, we have
\beq \label{sine}
\sin\angle (
\bB\widehat\bv, 
(\bB+\bE)\widehat\bv)=\frac{\|
	\bE\widehat\bv\|_2}{\|
	\bB\widehat\bv\|_2}\cdot\sin \angle (
\bE\widehat\bv,
\widehat\bv)
\eeq
Note that 
\begin{align}
	&\|\bB\widehat\bv\|_2 \le\|\bB\bv\|_2+\|\bB(\bv-\widehat\bv)\|_2\nonumber\\
	&\le \lambda_2+\lambda_{n}\cdot\|\bv-\widehat\bv\|_2\lesssim \frac{\delta}{n^2}+\frac{\delta}{n^{2}}\lesssim \sigma\sqrt{n}, \label{Bv}
\end{align}
where the second last inequality follows from  the assumption that $\|\bv-\widehat\bv\|_2\lesssim n^{-2}$, and Lemma 7. On the other hand, we have
\begin{align}
	&\|\bE\widehat\bv\|_2\ge \|\bE\bv\|_2-\|\bE(\bv-\widehat{\bv})\|_2\ge \|\bE\bv\|_2-\|\bE\|\cdot\|\bv-\widehat{\bv}\|_2\nonumber	\\
	&\ge C_1\sigma\sqrt{n}-C_2\frac{\sigma\sqrt{n\log n}}{n^2}\gtrsim \sigma\sqrt{n}, \label{Ev}
\end{align}
where the second last inequality follows from the definition of $\mathcal{A}_n$ (i.e., $E_1$ and $E_2$) and the assumption that $\|\bv-\widehat\bv\|_2\lesssim n^{-2}$. Finally, note that
$
	\cos \angle 
	(\bE\widehat\bv,\widehat\bv)=\frac{\widehat\bv^\top\bE\widehat\bv}{\|\bE\widehat\bv\|_2}.
$
For the numerator, we have
\begin{align*}
	&|\widehat\bv^\top\bE\widehat\bv|\le |\bv^\top\bE\bv|+|(\bv-\widehat\bv)^\top\bE\bv|+|\widehat\bv^\top\bE(\bv-\widehat\bv)|
	\\
	&\lesssim  \sigma\sqrt{\log n}+\sigma\sqrt{\frac{\log n}{n^3}}\lesssim \sigma\sqrt{\log n},
\end{align*}
where the second last inequality follows from the definition of $\mathcal{A}_n$ (i.e., $E_2$ and $E_3$) and the assumption $\|\bv-\widehat\bv\|_2\lesssim n^{-2}$. For the denominator, we already have $\|\bE\widehat\bv\|_2\gtrsim \sigma\sqrt{n}$
for sufficiently large $n$. Hence, it follows that
\beq
\cos \angle (\bE\widehat\bv,\widehat\bv)\le C \sqrt{\frac{\log n}{n}},
\eeq
which implies that, for any small constant $c>0$, we have
\beq
\sin \angle (\bE\widehat\bv,\widehat\bv)\ge 1-c
\eeq
for all sufficiently large $n$. Plugging in the above results back to (\ref{sine}), by the above arguments (\ref{Bv}) and (\ref{Ev}), we can choose $\delta=C_0\sigma n^2\sqrt{n}$ for a sufficiently small $C_0>0$, such that $\|\bE\widehat\bv\|_2/\|\bB\widehat\bv\|_2>1$. Therefore,
\beq
\sin\angle (\bB\widehat\bv, (\bB+\bE)\widehat\bv)\ge 1-c
\eeq
for some small constant $c>0$. This proves the first statement.

\paragraph{Step II} Suppose the first statement hold. Then for any small constant $c>0$, we have
\beq \label{cos}
\cos(\widehat\bv, (\bB-\bE)\widehat\bv) = \frac{\widehat\bv^\top (\bB-\bE)\widehat\bv}{\| (\bB-\bE)\widehat\bv\|_2}<c.
\eeq
For the numerator, we have
\beq 
\widehat\bv^\top (\bB-\bE)\widehat\bv=\sum_{i=1}^n \tilde\lambda_i( \tilde\bu_i^\top\widehat\bv)^2.
\eeq
Under the assumption $\|\bv-\widehat\bv_-\|_2\le cn^{-2}$ and  event $\mathcal{A}_n$ (i.e., $E_2$ and $E_3$), we have
\beq \label{l2.err}
|\tilde\lambda_2-\lambda_2|\le C\sigma\sqrt{\log n}.
\eeq
To see this, note that $\tilde\lambda_2=\widehat\bv_-^\top(\bB-\bE)\widehat\bv_-=\lambda_2(\bv^\top\widehat\bv_-)^2+\sum_{i\ne 2}\lambda_i(\bu_i^\top\widehat\bv_-)^2-\widehat\bv_-^\top\bE\widehat\bv_-$. Then we have
\begin{align*}
	|\tilde\lambda_2-\lambda_2|&\le \lambda_2|1-(\bv^\top\widehat\bv_-)^2|+\lambda_n\sum_{i\ne 2}(\bu_i^\top\widehat\bv_-)^2\\
	&\quad+|\bv^\top\bE\bv|+2|(\bv-\widehat\bv_-)^\top\bE\bv|\\
	&\lesssim \frac{\lambda_2}{n^4}+\frac{\lambda_n}{n^3}+\sigma\sqrt{\log  n}+\frac{\sigma\sqrt{n\log n}}{n^2},
\end{align*}
where the last inequality follows from
\beq \label{cos.lb}
|\widehat\bv_-^\top\widehat\bv|=|1-\frac{1}{2}\|\widehat\bv_--\widehat\bv\|_2^2|\ge 1-c_1n^{-4},
\eeq
and
\beq
\sum_{i\ne 2}(\bu_i^\top\widehat\bv_-)^2\le \sum_{i\ne 2}(\bu_i^\top\bv_-)^2+\sum_{i\ne 2}[\bu_i^\top(\bv_--\widehat\bv_-)]^2\le \frac{c_2}{n^3},
\eeq
for sufficiently large $n$ and some constants $c_1,c_2>0$. Thus, by (\ref{l2.err}), it holds that
\beq \label{l2.lb}
\tilde\lambda_2\ge \lambda_2-|\tilde\lambda_2-\lambda_2|\gtrsim \frac{\delta}{n^2}-\sigma\sqrt{\log n}>0,
\eeq
for sufficiently large $n$.
As $\bB-\bE$ is a Laplacian matrix, $\{\tilde\lambda_i\}_{1\le i\le n}$ must contain an eigenvalue $0$, it follows that $\tilde\lambda_1=0$, so that $\min_{1\le i\le n}\tilde\lambda_i\ge 0$. This implies that
\begin{align}
&\widehat\bv^\top (\bB-\bE)\widehat\bv\ge \tilde\lambda_2( \tilde\bu_2^\top\widehat\bv)^2=\tilde\lambda_2( \widehat\bv_-^\top\widehat\bv)^2\nonumber\\
&\ge (\lambda_2-|\lambda_2-\tilde\lambda_2|)\cdot|\widehat\bv_-^\top\widehat\bv|\cdot(1-cn^{-4})\gtrsim \frac{\delta |\widehat\bv_-^\top\widehat\bv|}{n^2},
\end{align}
where the second last inequality follows from  the assumption $\|\bv-\widehat\bv\|_2+\|\bv-\widehat\bv_-\|_2\lesssim n^{-2}$, (\ref{cos.lb}) and (\ref{l2.lb}).

For the denominator of (\ref{cos}), we have
\begin{align*}
	&\| (\bB-\bE)\widehat\bv\|_2\le \|\bB\bv\|_2+\|\bB(\widehat\bv-\bv)\|_2+\|\bE\bv\|_2+\|\bE(\widehat\bv-\bv)\|_2\\
	&\le \lambda_2+\frac{\lambda_n}{n^2}+C\sigma\sqrt{n}+\frac{C\sigma\sqrt{n\log n}}{n^2}\lesssim \frac{\delta}{n^2}
\end{align*}
Thus, we have
\beq
c>\frac{\widehat\bv^\top (\bB-\bE)\widehat\bv}{\| (\bB-\bE)\widehat\bv\|_2}\ge  C|\widehat\bv_-^\top\widehat\bv|.
\eeq
Since $c$ can be chosen arbitrarily small, we conclude that $|\widehat\bv_-^\top\widehat\bv|\le c'$ for some sufficiently small constant $c'> cC^{-1}$. In other words, $\|\widehat\bv-\widehat\bv_-\|_2\ge \epsilon_2$ for some constant $\epsilon_2>0$. This contradicts the first statement that $\|\widehat\bv-\widehat\bv_-\|_2\le cn^{-2}$. This completes the proof of Step II.

\paragraph{Step III} Finally, we show each of the following events holds with  probability at least $1-n^{-c}$ for some constant $c>0$:
\[	E_1=\bigg\{\|\bE\bv\|_2\ge C_1\sigma\sqrt{n}\bigg\},\quad E_2 = \bigg\{\|\bE\|\le C_2\sigma\sqrt{n\log n}  \bigg\},
\]
\[
	E_3= \bigg\{|\bv^\top\bE\bv|\le C_3\sigma\sqrt{\log n}\bigg\}.
\]
We start with $E_3$.
For $E_3$, note that
\[
\bv^\top\bE\bv = 2\sum_{i>j}z_{ij}(v_i-v_j)^2,
\]
is a normal random variable with mean 0 and variance
\begin{align*}
	\E\bigg[ 2\sum_{i>j}z_{ij}(v_i-v_j)^2 \bigg]^2 &= 2\sum_{i>j}\E z_{ij}^2(v_i-v_j)^4\\
	&=2\sigma^2\sum_{i>j}(v_i-v_j)^4,
\end{align*}
where 
\begin{align} \label{bv}
&\bv=(v_1,...,v_n)^\top \nonumber\\
&= \frac{1}{\sqrt{n/2}}\bigg(\cos \bigg(\frac{\pi}{2n} \bigg),\cos \bigg(\frac{3\pi}{2n} \bigg),..., \cos \bigg(\frac{(2n-1)\pi}{2n} \bigg) \bigg)^\top.
\end{align}
Calculate that
\[
\sum_{i>j}(v_i-v_j)^4=\sum_{i=1}^n[(v_i-v_{i+1})^4+...+(v_i-v_{n})^4]
\le n^2\cdot \frac{c}{n^2}\le c,
\]
where the second last inequality follows from
\beq \label{diff.v}
\max_{i\ne j} |v_i-v_j|\le \frac{c}{\sqrt{n}}.
\eeq
We know that $	\bv^\top\bE\bv$ is a normal random variable with bounded variance. Thus, by the tail bound of normal random variable, we have $E_3$ holds with probability at least $1-n^{-c}$.

For $E_1$, we denote $\bE\bv =(\xi_1,...,\xi_n)^\top$, so that
\[
\xi_i=\sum_{j\ne i}z_{ij}(v_i-v_j),\qquad i=1,...,n.
\]
In particular, we know that $\xi_i$ is a normal random variable with mean 0 and variance
\[
\E\xi_i^2 = \sum_{k,j\ne i}\E  z_{ij}z_{ik}(v_i-v_j)(v_i-v_k)=\sigma^2\sum_{j\ne i}(v_i-v_j)^2.
\]
By (\ref{bv}), we  have
\[
\sum_{j\ne i}(v_i-v_j)^2\le n\cdot \frac{c}{n^2}\le c.
\]
If  $i<n/2$,  then
\[
\sum_{j\ne i}(v_i-v_j)^2\ge \sum_{j=i+n/4}^{i+n/2}(v_i-v_j)^2\ge \frac{n}{4}\cdot\frac{c_1}{n}\ge c_2,
\]
for some constants $c_1,c_2>0$. By a similar argument, the above inequality $\sum_{j\ne i}(v_i-v_j)^2\ge c_2$ also holds if $i\ge n/2$. Thus we obtain $\E\xi_i^2 \asymp  1$ for all $i=1,...,n$. Moreover, for any $i\ne j$, we have
\begin{align}
&\E \xi_i\xi_j=\sum_{\ell\ne j}\sum_{k\ne i}\E  z_{ik}z_{j\ell}(v_i-v_k)(v_j-v_\ell)\nonumber\\
&=\E z_{ij}^2(v_i-v_j)^2=\sigma^2(v_i-v_j)^2.
\end{align}
In particular,  by (\ref{diff.v}), we have $(v_i-v_j)^2\le c/n$. Therefore, $\bE\bv$ is a multivariate normal random vector with mean 0 and covariance matrix $\Sigma_n$, whose diagonal entries are of order $\sigma^2$, and the off-diagonal entries are bounded in absolute value by $O(\sigma^2/n)$. By the property  of multivariate normal distribution, there exists a  matrix $\Sigma_n^{1/2}$ such that, for $i.i.d.$ standard normal random variables $\{\zeta_i\}_{1\le i\le n}$, 
\[
\bE\bv =(\xi_1,...,\xi_n)^\top = \Sigma_n^{1/2}(\zeta_1,...,\zeta_n)^\top.
\] 
Hence, the quadratic form $\bv^\top\bE^2\bv$ can be written as
$
\bv^\top\bE^2\bv = \bzeta^\top \Sigma_n\bzeta,
$
where $\bzeta=(\zeta_1,...,\zeta_n)^\top.$ Note  that $\|\Sigma_n\|^2_F\lesssim \sigma^2 n$. By the decomposition $\Sigma_n=\Sigma_n^D+\Sigma_n^O$ of $\Sigma_n$ into a diagonal matrix $\Sigma_n^D$ and an off-diagonal matrix $\Sigma_n^O$, we also have
\[
\|\Sigma_n\|\le  \|\Sigma_n^D\|+\|\Sigma_n^O\|\le \max_{1\le i\le n}[\Sigma_n]_{ii}+\|\Sigma_n^O\|_F=O(\sigma^2).
\]
Now we could apply Hanson-Wright inequality \citep{rudelson2013hanson}, to obtain the  concentration inequality
\beq
P(|\bzeta^\top\Sigma_n\bzeta-\E\bzeta^\top\Sigma_n\bzeta|>t)\le 2\exp(-c\min\{t^2/(n\sigma^2), t/\sigma^4\}),
\eeq
for $t\ge 0.$
By choosing $t=C\sigma^2\sqrt{n\log n}$, we have
\beq
P(|\bzeta^\top\Sigma_n\bzeta-\E\bzeta^\top\Sigma_n\bzeta|>\sigma^2\sqrt{n\log n})\le n^{-c}.
\eeq
Finally, since
$
\E\bzeta^\top\Sigma_n\bzeta=\sum_{i=1}^n [\Sigma_n]_{ii}\asymp \sigma^2n,
$
we have
\[
P(\bzeta^\top\Sigma_n\bzeta\le 	\E\bzeta^\top\Sigma_n\bzeta-C\sigma^2\sqrt{n\log n})\le n^{-c},
\]
or
\[
P(\bv^\top\bE^2\bv\ge 	C\sigma^2n)\ge 1- n^{-c}.
\]
This proves $E_1$ holds with probability at least $1-n^{-c}$, by noting that $\bv^\top\bE^2\bv=\|\bE\bv\|_2^2$. 

For $E_2$, to obtain an upper bound for $\|\bE\|$, we note that
\[
\|\bE\|=\|\bL(\bZ)\|\le \|\bD(\bZ)\|+\|\bZ\|\le \max_{1\le i\le n} \sum_{j=1}^nZ_{ij}+\|\bZ\|.
\]
Since the Laplacian operator $\bL(\cdot)$ is invariant to the diagonal entries of the input, we assume without loss of generality that the diagonal entries of $\bZ$ is $i.i.d.$  Gaussian $N(0,2\sigma^2)$, so that $\bZ/\sigma$ is Gaussian Orthogonal Ensemble (GOE).  The upper bounds for $\max_{1\le i\le n} \sum_{j=1}^nZ_{ij}$ and $\|\bZ\|$ are obtained previously. Specifically, the well-celebrated Bai-Yin theorem \citep{bai1988necessary} implies that for any $\epsilon>0$,
\beq
\lim_{n\to\infty}P\bigg(\frac{\|\bZ\|}{\sigma\sqrt{n}}\le 2+\epsilon\bigg)=1,
\eeq
whereas  the standard tail bound for the Gaussian random variable and an union bound argument imply that for any $\epsilon>0$,
\beq
\lim_{n\to\infty}P\bigg(\max_{1\le i\le n} \sum_{j=1}^nZ_{ij}\le \sigma\sqrt{(2+\epsilon)(n+1)\log n}\bigg)=1.
\eeq
These results imply
\beq
\lim_{n\to\infty}P(\|\bE\|\le 8\sigma\sqrt{n\log n})=1.
\eeq

\subsection{Theoretical Guarantee for $\lambda$-Ridged Toeplitz Matrices: Proof of Proposition \ref{eff.thm.2}} \label{prop.2}

The proof of this proposition is separated into two parts, corresponding to the initialization (Step 1) and the iterative sorting part (Step 2) of the algorithm, respectively.

\paragraph{Part I} We start by showing that the initialization step successfully identifies the first or the last row of $\bTheta$. Set $\delta_n=n^{-c}$ for some $c>0$. Without loss of generality, we assume  $\Pi=\text{Id}_n$. In other words, we need to show that, 
\beq \label{o1}
P_{\bTheta,\Pi}(\widetilde\pi(1)\in\{1,n\})\ge 1-\delta_n/n,
\eeq
for any $(\bTheta,\Pi)\in \calR_n^R(\lambda)\times \calS_n$.
In fact, (\ref{o1}) is implied by
\beq 
P_{\bTheta,\Pi}(S_1<\min_{2\le i\le n-1}S_i)\ge 1-\delta_n/n.
\eeq
Then it suffices to show that, for some $\tau_n$,
\beq \label{o2}
P_{\bTheta,\Pi}(S_1<\tau_n)\ge 1-\frac{\delta_n}{2n},
\eeq
and
\beq \label{o3}
P_{\bTheta,\Pi}(\min_{2\le i\le n-1}S_i>\tau_n)\ge 1-\frac{\delta_n}{2n}.
\eeq
The rest of the proof is devoted to (\ref{o2}) and (\ref{o3}).

On the one hand, let $S_1=\sum_{i\ne 1}Y_{1i}=\sum_{1\le i\le  n-1}\theta_i+\sum_{i\ne 1}Z_{1i}$ where $Z_{1i}$'s are independent subgaussian random variables. By the concentration inequality for subgaussian random variables \citep{vershynin2010introduction}, we have
\beq \label{eq1}
P\bigg(S_1 -\sum_{i=1}^{n-1}\theta_i< C\sigma \sqrt{n}t\bigg)\ge 1-e^{-t^2}.
\eeq
On the other hand, for any $2\le i\le n-1$, we have $S_i=\sum_{j\ne i}Y_{ij}=\sum_{j=1}^{i-1}\theta_j+\sum_{j=1}^{n-i}\theta_j+\sum_{j\ne i}Z_{ij}$ where $Z_{ij}$'s are independent subgaussian random variables, so that
\beq\label{eq2}
P\bigg(S_i -\sum_{j=1}^{i-1}\theta_j-\sum_{j=1}^{n-i}\theta_j> -C\sigma \sqrt{n}t\bigg)\ge 1-e^{-t^2}.
\eeq
By setting $t=\sqrt{\log (4n/\delta_n)}$ in (\ref{eq1})  and $t=\sqrt{\log (4n^2/\delta_n)}$ in (\ref{eq1}), we have 
\beq
P\bigg(S_1 < \sum_{i=1}^{n-1}\theta_i+C\sigma \sqrt{n\log (4n/\delta_n)}\bigg)\ge 1-\frac{\delta_n}{4n},
\eeq
\beq\label{eq3}
P\bigg(S_i >\sum_{j=1}^{i-1}\theta_j+\sum_{j=1}^{n-i}\theta_j -C\sigma \sqrt{n\log (4n^2/\delta_n)}\bigg)\ge 1-\frac{\delta_n}{4n^2}.
\eeq 
Since for any $2\le i\le n-1$, we have
\[
\sum_{j=1}^{i-1}\theta_j+\sum_{j=1}^{n-i}\theta_j\ge\theta_1+\sum_{j=1}^{n-2}\theta_j,
\]
then we have
\[
P\bigg(S_i >\theta_1+\sum_{j=1}^{n-2}\theta_j -C\sigma \sqrt{n\log (4n^2/\delta_n)}\bigg)\ge 1-\frac{\delta_n}{4n^2}
\]
for all $2\le i\le n-1$.
A union bound applied to (\ref{eq3}) implies
\[
P\bigg(\min_{2\le i\le n} S_i >\theta_1+\sum_{j=1}^{n-2}\theta_j -C\sigma \sqrt{n\log (4n^2/\delta_n)}\bigg)\ge 1-\frac{\delta_n}{4n}.
\]
Now note that by assumption we have
\[
\theta_1+\sum_{j=1}^{n-2}\theta_j -C\sigma \sqrt{n\log (4n^2/\delta_n)}\ge \sum_{i=1}^{n-1}\theta_i+C\sigma \sqrt{n\log (4n/\delta_n)},
\]
as for $\lambda\ge 4C\sigma n$ we have
\[
\theta_1-\theta_{n-1} \ge \lambda\ge C\sigma \sqrt{n\log (4n/\delta_n)}+C\sigma \sqrt{n\log (4n^2/\delta_n)},
\]
for sufficiently large $n$. Therefore, we have shown (\ref{o2}) and (\ref{o3}) with $\tau_n=C\sigma \sqrt{n\log (4n/\delta_n)}+C\sigma \sqrt{n\log (4n^2/\delta_n)}$. This proves (\ref{o1}).

\paragraph{Part II} In this part, we show that, for each $1\le k\le n-2$,
\beq \label{ite1}
P_{\bTheta,\Pi}(\widetilde\pi({k+1})=k+1| \widetilde\pi(i)=i, \forall i\le k) \ge 1-\delta_n/n,
\eeq
\beq \label{ite10}
P_{\bTheta,\Pi}(\widetilde\pi({k+1})=n-k| \widetilde\pi(i)=n-i+1, \forall i\le k) \ge 1-\delta_n/n.
\eeq
This along with Part I implies that
\begin{align*}
	&P_{\bTheta,\Pi}(\widetilde\pi=id)\\
	 &=P(\widetilde\pi(1)=1) \prod_{k=1}^{n-2}P_{\bTheta,\Pi}(\widetilde\pi(k+1)=k+1| \widetilde\pi(i)=i, \forall i\le k)\\
	 & \ge P(\widetilde\pi(1)=1)(1-\delta_n/n)^{n-2},
\end{align*}
\begin{align*}
&	P_{\bTheta,\Pi}(\widetilde\pi=id^{-1}) \\
	&=P(\widetilde\pi(1)=n) \times\\
	&\quad\prod_{k=1}^{n-2}P_{\bTheta,\Pi}(\widetilde\pi(k+1)=n-k| \widetilde\pi(i)=n-i+1, \forall i\le k) \\
	&\ge P(\widetilde\pi(1)=n)(1-\delta_n/n)^{n-2},
\end{align*}
for $n\ge 1$. Combining the above two inequalities, by the result from Part I we have
\begin{align*}
	&P_{\bTheta,\Pi}(\widetilde\pi\in\{id,id^{-1}\}) \\
	&\ge [P(\widetilde\pi(1)=1)+P(\widetilde\pi(1)=n)](1-\delta_n/n)^{n-2}\\
	&\ge P(\widetilde\pi(1)\in\{1,n\})(1-\delta_n/n)^{n-2}\\
	&\ge (1-\delta_n/n)^{n-1}\ge 1-\delta_n,
\end{align*}
which in turn proves the proposition.

The rest of the proof is devoted to (\ref{ite1}), as the proof of (\ref{ite10})  follows by symmetry. Suppose $\widetilde\pi(i)=i, \forall i\le k$. Then the event $\widetilde\pi(k+1)=k+1$ is equivalent to
\[
\|\bY_{k,-k}-\bY_{k+1,-(k+1)}\|_1<\min_{j\in[n]\setminus[k+1] }\|\bY_{k,-k}-\bY_{j,-j}\|_1,
\] 
or
\[
\|\bY_{k,-k}-\bY_{k+1,-(k+1)}\|_1<\|\bY_{k,-k}-\bY_{j,-j}\|_1,
\]
for all $j\in[n]\setminus[k+1]$.
In the following, we show that, for any $j\in[n]\setminus[k+1]$, it holds that
\beq \label{ite2}
P(\|\bY_{k,-k}-\bY_{k+1,-(k+1)}\|_1<\|\bY_{k,-k}-\bY_{j,-j}\|_1) \ge 1-\frac{\delta_n}{2n^2}.
\eeq
Thus, by applying the union bound, we have
\begin{align} \label{ite3}
&P(\|\bY_{k,-k}-\bY_{k+1,-(k+1)}\|_1\le \min_{j\in[n]\setminus[k+1] }\|\bY_{k,-k}-\bY_{j,-j}\|_1)\nonumber \\
&\ge 1-\frac{\delta_n}{2n}.
\end{align}
To obtain (\ref{ite2}), we will show that there exists some $\tau_n$ such that
\beq\label{ite4}
P(\|\bY_{k,-k}-\bY_{k+1,-(k+1)}\|_1\le \tau_n) \ge 1-\frac{\delta_n}{4n^2},
\eeq
\beq\label{ite5}
P(\tau_n\le \|\bY_{k,-k}-\bY_{j,-j}\|_1) \ge 1-\frac{\delta_n}{4n^2}.
\eeq
On the one hand, if we denote $\bY_{k,-k}-\bY_{k+1,-(k+1)}=(\xi^{(k+1)}_1, \xi^{(k+1)}_2,..., \xi^{(k+1)}_{n-1})$, it follows that $\{\xi^{(k+1)}_i\}$ are independent subgaussian variables with variances bounded by $C\sigma^2$, and means satisfying
\beq 
\E\xi^{(k+1)}_i =\theta_{k-i}-\theta_{k-i+1},\qquad  \text{for $i<k$},
\eeq 
\beq 
\E\xi^{(k+1)}_i =\theta_{i-k+1}-\theta_{i-k},\qquad  \text{for $i>k$},
\eeq
and $\E \xi^{(k+1)}_k=0$. To see this, note that the $k$-th row of the signal matrix $\bTheta$ with the main diagonal removed is
\[
[\underbrace{\theta_{k-1}\quad \theta_{k-2} \quad ...\quad \theta_2\quad \theta_1}_{k-1}\quad \underbrace{\theta_1\quad \theta_2\quad...\quad \theta_{n-k}}_{n-k}]
\]
Then, we can write $\|\bY_{k,-k}-\bY_{k+1,-(k+1)}\|_1=\sum_{i=1}^{n-1}|\xi_i^{(k+1)}|=\sum_{i=1}^{n-1}|\E\xi_i^{(k+1)}+w_i|$ for some independent subgaussian random variables $w_i$ with mean 0 and variance bounded by $C\sigma^2$. By the concentration inequality for subgaussian random variables, we have
\begin{align*}
	P\bigg( \sum_{i=1}^{n-1}|\xi_i^{(k+1)}|\le\sum_{i=1}^{n-1}\E\xi_i^{(k+1)}+\sum_{i=1}^{n-1}|w_i|\le \tau_1\bigg)\ge 1-\frac{\delta_n}{2n^2},
\end{align*}
where
$	\tau_{1}=\sum_{i<k}|\theta_{k-i}-\theta_{k-i+1}|+\sum_{i>k}|\theta_{i-k+1}-\theta_{i-k}|
	+C\sigma n+C\sigma \sqrt{n\log (2n^2/\delta_n)}.
$
On the other hand, for any $j\in[n]\setminus[k+1]$, if we denote $\bY_{k,-k}-\bY_{j,-j}=(\xi^{(j)}_1, \xi^{(j)}_2,..., \xi^{(j)}_{n-1})$, it follows that $\{\xi^{(j)}_i\}$ are subgaussian variables with
\begin{align*}
	\E\xi^{(j)}_i &=\theta_{k-i}-\theta_{j-i},\qquad  \text{for $i<k$},\\
	\E\xi^{(j)}_i &=\theta_{i-k+1}-\theta_{j-i},\qquad  \text{for $k\le i<j$},\\
	\E\xi^{(j)}_i &=\theta_{i-k+1}-\theta_{i-j+1},\qquad  \text{for $ i\ge j$}.
\end{align*}
In addition, out of the $(n-1)$ elements in $\{\xi^{(j)}_i\}$, $(n-3)$ elements are mutually independent, whereas the rest two elements, corresponding to the entries $Y_{kj}$ and $Y_{jk}$, are correlated with covariance bounded by $C\sigma^2$, but independent from the others.
By the similar argument for the $(n-3)$ independent variables and by the sub-Gaussian property of the two dependent variables, we still have
\begin{align*}
	P\bigg( \sum_{i=1}^{n-1}|\xi_i^{(j)}|\ge \tau_2\bigg)\ge 1-\frac{\delta_n}{2n^2},
\end{align*}
where 
$
	\tau_2= \sum_{i<k}|\theta_{k-i}-\theta_{j-i}|
	+\sum_{k\le i<j}|\theta_{i-k+1}-\theta_{j-i}|+\sum_{i\ge j}|\theta_{i-k+1}-\theta_{i-j+1}|-Cn\sigma-C\sigma\sqrt{n\log (2n^2/\delta_n)}.
$
Now we claim that, if $\lambda\ge 4C\sigma n$, we have, $\tau_2\ge \tau_1$ for sufficiently large $n$, or
\begin{align}
	&    \quad\sum_{i<k}|\theta_{k-i}-\theta_{j-i}|
	+\sum_{k\le i<j}|\theta_{i-k+1}-\theta_{j-i}|\nonumber\\
	&\quad+\sum_{i\ge j}|\theta_{i-k+1}-\theta_{i-j+1}|\nonumber \\
	&\ge \sum_{i<k}|\theta_{k-i}-\theta_{k-i+1}|+\sum_{i>k}|\theta_{i-k+1}-\theta_{i-k}|\nonumber\\
	&\quad+2Cn\sigma+2C\sigma\sqrt{n\log (4n^2/\delta_n)}. \label{t2>t1}
\end{align}
To see this, it suffices to show that,
by  the definition of $\calT_n^R(\lambda)$,  it holds that
\begin{align}
	&  \sum_{i<k}|\theta_{k-i}-\theta_{j-i}|
	+\sum_{k\le i<j}|\theta_{i-k+1}-\theta_{j-i}|\nonumber\\
	&\quad+\sum_{i\ge j}|\theta_{i-k+1}-\theta_{i-j+1}|-\sum_{i<k}|\theta_{k-i}-\theta_{k-i+1}|\nonumber \\
	&\quad-\sum_{i>k}|\theta_{i-k+1}-\theta_{i-k}|\nonumber\\
	&\ge \sum_{i=1}^{\lceil n/2\rceil -1}|\theta_i-\theta_{i+1}|\ge \lambda.\label{TR.ineq}
\end{align}
Specifically, if we denote $\delta_k=\theta_{k}-\theta_{k+1}$ for $1\le k\le n-2$, then
\begin{align*}
&\sum_{i<k}|\theta_{k-i}-\theta_{k-i+1}|+\sum_{i>k}|\theta_{i-k+1}-\theta_{i-k}|\\
&= \sum_{i<k}\delta_{k-i}+\sum_{i>k}\delta_{i-k}=\sum_{i=1}^{k-1}\delta_i+\sum_{i=1}^{n-k-1}\delta_i,
\end{align*}
and
\begin{align}
	&  \sum_{i<k}|\theta_{k-i}-\theta_{j-i}|
	+\sum_{k\le i<j}|\theta_{i-k+1}-\theta_{j-i}|\nonumber\\
	&\quad+\sum_{i\ge j}|\theta_{i-k+1}-\theta_{i-j+1}|\nonumber \\
	&\ge \sum_{i<k}|\theta_{k-i}-\theta_{j-i}|+2\delta_1+\sum_{i\ge j}|\theta_{i-k+1}-\theta_{i-j+1}|\nonumber \\
	&=\sum_{i<k}(\delta_{k-i}+\delta_{k-i+1}+...+\delta_{j-i-1})\nonumber\\
	&\quad+\sum_{i\ge j}(\delta_{i-j+1}+\delta_{i-j+2}+...+\delta_{i-k})+2\delta_1\nonumber\\
	&=\sum_{i=1}^{k-1}(\delta_{i}+\delta_{i+1}+...+\delta_{i+j-k-1})\nonumber\\
	&\quad+\sum_{i=1}^{n-j+1}(\delta_{i}+\delta_{i+1}+...+\delta_{i+j-k-1})+2\delta_1.
\end{align}
Thus, we have
\begin{align*}
	&\sum_{i<k}|\theta_{k-i}-\theta_{j-i}|
	+\sum_{k\le i<j}|\theta_{i-k+1}-\theta_{j-i}|\\
	&\quad+\sum_{i\ge j}|\theta_{i-k+1}-\theta_{i-j+1}|-\sum_{i<k}|\theta_{k-i}-\theta_{k-i+1}|\\
	&\quad-\sum_{i>k}|\theta_{i-k+1}-\theta_{i-k}|\\
	&\ge \sum_{i=1}^{k-1}(\delta_{i}+\delta_{i+1}+...+\delta_{i+j-k-1})
	\\
	&\quad+\sum_{i=1}^{n-j+1}(\delta_{i}+\delta_{i+1}+...+\delta_{i+j-k-1})\\
	&\quad+2\delta_1-\sum_{i=1}^{k-1}\delta_i-\sum_{i=1}^{n-k-1}\delta_i\\
	&\ge \sum_{i=1}^{k}\delta_{i}+
	\sum_{i=1}^{n-k}\delta_{i}\ge \sum_{i=1}^{\lceil n/2\rceil -1}\delta_{i},
\end{align*}
where in the second last inequality we used the fact that
\beq
\sum_{i=1}^{n-j+1}(\delta_{i}+\delta_{i+1}+...+\delta_{i+j-k-1})-\sum_{i=1}^{n-k-1}\delta_i\ge  \sum_{i=2}^{n-k}\delta_{i}.
\eeq
To see this, note that
\begin{align*}
	&\sum_{i=1}^{n-j+1}(\delta_{i}+\delta_{i+1}+...+\delta_{i+j-k-1})-\sum_{i=1}^{n-k-1}\delta_i\\
	&=(\delta_{1}+\delta_{2}+...+\delta_{j-k})+(\delta_{2}+\delta_{3}+...+\delta_{j-k+1})+\\
	&\quad...+(\delta_{n-j+1}+...+\delta_{n-k})-\sum_{i=1}^{n-k-1}\delta_i.
\end{align*}
If $n-j+1<j-k$, then we have
\[
(\delta_{1}+\delta_{2}+...+\delta_{j-k})+(\delta_{n-j+1}+...+\delta_{n-k})-\sum_{i=1}^{n-k-1}\delta_i\ge \delta_{n-k},
\]
so that
\begin{align*}
	&(\delta_{1}+\delta_{2}+...+\delta_{j-k})+(\delta_{2}+\delta_{3}+...+\delta_{j-k+1})+\\
	&\quad...+(\delta_{n-j+1}+...+\delta_{n-k})-\sum_{i=1}^{n-k-1}\delta_i\\
	&\ge (\delta_{2}+\delta_{3}+...+\delta_{j-k+1})+\\
	&\quad...+(\delta_{n-j}+...+\delta_{n-k-1})+\delta_{n-k}\\
	&\ge \sum_{i=2}^{n-k} \delta_i.
\end{align*}
If $n-j+1\ge j-k$, then we can "extract" the first terms in each of the first $n-j$ sums, and all but the last term $\delta_{n-k}$  in the $(n-j+1)$-th sum, and get
\begin{align*}
	&(\delta_{1}+\delta_{2}+...+\delta_{j-k})+(\delta_{2}+\delta_{3}+...+\delta_{j-k+1})+\\
	&\quad...+(\delta_{n-j+1}+...+\delta_{n-k})-\sum_{i=1}^{n-k-1}\delta_i\\
	&\ge (\delta_{2}+\delta_{3}+...+\delta_{j-k})+(\delta_{3}+\delta_{4}+...+\delta_{j-k+1})+\\
	&\quad...+(\delta_{n-j+1}+...+\delta_{n-k-1})+\delta_{n-k}\ge \sum_{i=2}^{n-k} \delta_i.
\end{align*}
This proves (\ref{TR.ineq}) or (\ref{t2>t1}) under the condition that $\bTheta\in \calT_n^R(\lambda)$. Thus, we can take any $\tau_n\in[\tau_1,\tau_2]$, to get (\ref{ite4}) and (\ref{ite5}), which imply (\ref{ite2}) and (\ref{ite3}). 
This proves (\ref{ite1}) under the conditions of the proposition.

\subsection{Sharp $\lambda$-$\rho^*$ Correspondence: Proof of Proposition \ref{p.space}} \label{prop.3}

\paragraph{Sufficiency}
By the conditions of Proposition \ref{p.space}, for all $\bTheta\in\calT'_n$, we have $\rho(\bTheta,\mathcal{S}'_n)\ge C n\lambda^*$. Denote the diagonal values of $\bTheta$ as $\{\theta_0,...,\theta_{n-1}\}$. By the definition of $\rho(\bTheta,\mathcal{S}'_n)$, we have
\begin{align*}
	[\rho(\bTheta,\mathcal{S}'_n) ]^2&= \min_{\Pi_1,\Pi_2\in\mathcal{S}'_n}\|\Pi_1\bTheta\Pi_1^\top - \Pi_2\bTheta\Pi_2^\top\|_F^2\\
	&\le n(n-1)\cdot 4( \theta_1-\theta_{\lceil n/2\rceil })^2,
\end{align*}
where the last inequality follows from the fact that
$
	(\theta_j-\theta_k)^2\le (\theta_1-\theta_n)^2
$  for any $j\ne k$,
and the  "ridge" condition $
\theta_1-\theta_{\lceil n/2\rceil }\ge \theta_{\lceil n/2\rceil}-\theta_{n-1}$.
Thus, 	 $\rho(\bTheta,\mathcal{S}'_n)\ge C n\lambda^*$ implies 
$
\theta_1-\theta_{\lceil n/2\rceil }\ge \frac{C\lambda^*}{2},
$
that is, $\bTheta\in \calT_n^R(\lambda)$ for  $\lambda= C_2\lambda^*$ where $C_2=C/2$. Since this holds for any such $\bTheta\in \calT'_n$, we have proven the first part of Proposition \ref{p.space}.

\paragraph{Necessity}
To begin with, we first show that, there exists some $\calT'_n\subseteq \calT_n$ and $\mathcal{S}'_n\subseteq \mathcal{S}$ such that
\beq \label{cond.T'}
\rho^*(\calT'_n,\calS'_n)=C \sigma n\lambda,\qquad \calT'_n\subseteq \calT_n^R(\lambda),
\eeq
for some absolute constant $C>0$. To prove this,  without loss of generality, we assume $n=4m$ for some integer $m>0$, and consider the class of symmetric Toeplitz matrices 
\[
\calT'_{n} = \bigg\{\bTheta\in \calT_n:  
\begin{aligned}\theta_0=0, \theta_k=\alpha+(n-k)\beta, \\
\forall k\in\{1,2,...,n-1\}, \alpha\ge 0\end{aligned}\bigg\},
\]
where $\beta(n/2-1)=\lambda$. As a result, we can easily check that $\calT'_n\subseteq \calT_n^R(\lambda)$.
On the other hand, we consider the permutation set $\calS'_n\subseteq \calS_n$ including only two permutations $\pi_0=id$ and 
$
\pi_1= (n,..., 3n/4+1, n/4+1,...,3n/4, n/4,...,1).
$
In other words, $\pi_1$ exchanges the first $n/4$ elements with the last $n/4$ elements, and arranges them in the reversed order. Let $\Pi_0$ and $\Pi_1$ be the  permutation matrices associated with $\pi_0$ and $\pi_1$, respectively. In the following, we show that for any $\bTheta\in\calT'_n$, we have
\beq \label{cond.T'2}
\|\Pi_0\bTheta\Pi_0^\top-\Pi_1\bTheta\Pi_1^\top\|_F^2=C_0n^{2}\lambda^2,
\eeq
which implies (\ref{cond.T'}). To obtain (\ref{cond.T'2}), we denote $\bA=(a_{ij})=\Pi_0\bTheta\Pi_0^\top$ and $\bB=(b_{ij})=\Pi_1\bTheta\Pi_1^\top$, and calculate the differences $\|\bA_{.i}-\bB_{.i}\|_2^2$ carefully. 

Due to the invariance of the difference with respect to a translation of $\bTheta$, we assume without loss of generality that $\alpha=0$. For $i=1$,  we have
\begin{align*}
	\bA_{.1}&= \bTheta_{.1}=(0, (n-1)\beta,(n-2)\beta,...,2\beta,\beta)^\top,\\
	\bB_{.1}&=\Pi_1\bTheta_{.n}=\Pi_1 \cdot(\beta,2\beta,...,(n-1)\beta,0)^\top\\
	&=(\underbrace{0,(n-1)\beta,(n-2)\beta,..., (3n/4+1)\beta}_{n/4},\\
	&\qquad \underbrace{(n/4+1)\beta,...,(3n/4)\beta}_{n/2}, \underbrace{(n/4)\beta,...,2\beta,\beta}_{n/4})^\top.
\end{align*}
Thus,
\begin{align*}
&\|\bA_{.1}-\bB_{.1}\|_2^2=2\beta^2\sum_{k=0}^{n/4-1}(n/2-2k-1)^2\\
&=2\beta^2\sum_{k=1}^{n/4}(2k-1)^2= \frac{\beta^2 n(n/2-1)(n/2+1)}{6}.
\end{align*}
For $i=2$, we have
\begin{align*}
	\bA_{.2}&= \bTheta_{.2}=((n-1)\beta,0, (n-1)\beta,(n-2)\beta,...,3\beta,2\beta)^\top\\
	\bB_{.2}&=\Pi_1\bTheta_{.(n-1)}=\Pi_1 \cdot(2\beta,3\beta,...,(n-1)\beta,0,(n-1)\beta)^\top\\
	&=(\underbrace{(n-1)\beta,0,(n-1)\beta,..., (3n/4+2)\beta}_{n/4},\\
	&\quad \underbrace{(n/4+2)\beta,...,(3n/4+1)\beta}_{n/2}, \underbrace{(n/4+1)\beta,...,3\beta,2\beta}_{n/4})^\top.
\end{align*}
Thus,
\[
\|\bA_{.2}-\bB_{.2}\|_2^2=2\beta^2\sum_{k=1}^{n/4}(2k-1)^2= \frac{\beta^2 n(n/2-1)(n/2+1)}{6}.
\]
Similarly, one can show that,
\beq\label{n/4}
\|\bA_{.i}-\bB_{.i}\|_2^2=  \frac{\beta^2 n(n/2-1)(n/2+1)}{6}, 
\eeq
for  all $i\in\{1,2,...,n/4\}$.
Hence, by adding up the above differences, we have
\[
\|\bA-\bB\|_F^2\ge  \frac{\beta^2 n^2(n/2-1)(n/2+1)}{24}\ge \frac{\beta^2n^4}{24\times 8}
\]
for $n\ge 3$. On the other hand, we have
\begin{align*}
	&\|\bA-\bB\|_F^2\le 2\|\bA\|_F^2+2\|\bB\|_F^2=4\|\bTheta\|_F^2\\
	&\le 2\beta^2\sum_{k=1}^{n-1}k^3=\frac{\beta^2n^2(n-1)^2}{2}\le \beta^2n^4.
\end{align*}
Combining the above results, we have shown $\|\bA-\bB\|_F^2\asymp \beta^2n^4\asymp \lambda^2n^2$, or (\ref{cond.T'2}).

\subsection{An Information Lower Bound for Adaptive Sorting}\label{lb.sec}

In this part, we show that, for the adaptive sorting to achieve exact matrix reordering, a minimal signal strength condition that is strictly stronger than that required by the constrained LSE is in fact needed. The result is summarized as the following theorem.

\bet \label{eff.lb.thm}
Suppose the noise matrix has $i.i.d.$ entries up to symmetry generated from $N(0,\sigma^2)$. Then there exists some $\calT'_n\subset \calT_n$ and  $\mathcal{S}_n'\subseteq \mathcal{S}_n$ satisfying $\rho^*(\calT'_n,\mathcal{S}'_n)=C\sigma n^{3/2}$ for some absolute constant $C>0$, such that $\sup_{(\bTheta, \Pi)\in \calT'_n\times \mathcal{S}'_n}P_{\Theta,\Pi}(\widetilde\Pi\bTheta \widetilde\Pi^\top\ne \Pi\bTheta\Pi^\top)\ge 0.2.$
\eet

In particular, in light of the second part  of Proposition \ref{p.space}, it suffices to show the following proposition. 

\bep[Lower bound for $\lambda$-ridged Topelitz matrices] \label{eff.thm.3}
Suppose the noise matrix $\bZ$ has $i.i.d.$ entries up to symmetry generated from $N(0,\sigma^2)$. Then there exists some $\calT'_n\subseteq \calT_n^R(\lambda)$ where $\lambda= c\sigma \sqrt{n}$ for some absolute constant $c>0$, and $\calS'_n\subseteq\calS_n$, such that $\sup_{(\bTheta, \Pi)\in \calT'_n\times \mathcal{S}'_n}P_{\Theta,\Pi}(\widetilde\Pi\bTheta \widetilde\Pi^\top\ne \Pi\bTheta\Pi^\top)\ge 0.2.$ 
\eep

\begin{proof}
	To prove this proposition, we start with the set of symmetric  Toeplitz matrices
	\[
	\calT'_{n} = \bigg\{\bTheta\in \calT_n:  \begin{aligned}
	\theta_0=0, \theta_k=\alpha+(n-k)\beta, \\
	\forall k\in\{1,2,...,n-1\}, \alpha\ge 0\end{aligned}\bigg\},
	\]
	with $\beta=\lambda/(n-2)$, and the permutation set $\calS'_n\subseteq \calS_n$ including only two permutations $\pi_0=id$ and 
	$
	\pi_1= (n,..., 3n/4+1, n/4+1,...,3n/4, n/4,...,1),
	$
	constructed in the proof of the second statement of Proposition \ref{p.space}. In particular, we can easily check that $\calT'_n\subseteq \calT_n^R(\lambda)$, and that, for any $\bTheta\in\calT'_n$, we have
	\beq 
	\|\Pi_0\bTheta\Pi_0^\top-\Pi_1\bTheta\Pi_1^\top\|_F^2=C_0n^{2}\lambda^2,
	\eeq
	In what follows, we will show that, for $\lambda=c\sigma\sqrt{n}$, we have $P_{\Theta,\Pi}(\widetilde\Pi\bTheta \widetilde\Pi^\top\ne \Pi\bTheta\Pi^\top)\ge 0.2$ for any $(\bTheta,\Pi)\in\calT'_n\times \calS'_n$. In fact,  we only need to show that, in such a case, 
	\beq \label{ite1.2}
	P_{\bTheta,\Pi}(\widetilde\pi(1)\notin \{1,n\}) \ge 0.2.
	\eeq
	To show this, we show that
	\beq \label{ite1.3}
	P_{\bTheta,\Pi}(S_1>S_2)\ge 0.49,\qquad 	P_{\bTheta,\Pi}(S_n>S_{n-1})\ge 0.49,
	\eeq
	where $S_j=\sum_{i\ne j}Y_{ji}$. If we are able to show that event $\{S_1>S_2\}$ is independent of the event $\{S_n>S_{n-1}\}$, then
	we have
	\begin{align*}
		&P_{\bTheta,\Pi}(\widetilde\pi(1)\notin \{1,n\})\ge P_{\bTheta,\Pi}(\min\{S_1,S_n\}>\min_{2\le i\le n-1}S_i)\\
		&\ge	P_{\bTheta,\Pi}(S_1>S_2,S_n>S_{n-1})\\
		&\ge P_{\bTheta,\Pi}(S_1>S_2)\cdot P_{\bTheta,\Pi}(S_n>S_{n-1})\ge 0.2.
	\end{align*}
	The rest of the proof is devoted to (\ref{ite1.3}) and the independence between $\{S_1>S_2\}$ and $\{S_n>S_{n-1}\}$.
	
	On the one hand, note that $S_1-S_2=\sum_{i\ne 1}Y_{1i}-\sum_{i\ne 2}Y_{2i}=\sum_{i=1}^{n-1}\theta_i+\sum_{i\ne 1}Z_{1i}-\theta_1-\sum_{j=1}^{n-2}\theta_j-\sum_{j\ne 2}Z_{2j}$ where $Z_{ki}$'s are independent Gaussian random variables for $k=1,2$. By the tail bound of the Gaussian random variable $\sum_{i\ne 1}Z_{1i}-\sum_{j\ne 2}Z_{2j}\sim N(0, 2(n-1)\sigma^2)$, 
	we have
	\begin{align*}
		P(S_1-S_2\ge 0)&= P\bigg(\theta_{n-1}-\theta_1+\sum_{i\ne 1}Z_{1i}-\sum_{j\ne 2}Z_{2j}\ge 0\bigg)\\
		&=\Phi\bigg(\frac{\theta_{n-1}-\theta_1}{\sigma\sqrt{2(n-1)}}\bigg)\ge 0.49
	\end{align*}
	if $\lambda = c\sigma\sqrt{n}$ for some small constant $c>0$. Similarly, we can also show $P(S_n-S_{n-1}\ge 0)\ge 0.49$. This completes the proof of  (\ref{ite1.3}). 
	
	On the other hand, we write
	\begin{align*}
		&S_1-S_2\\
		&=\sum_{i=1}^{n-1}\theta_i-\theta_1-\sum_{j=1}^{n-2}\theta_j+\sum_{i\ne 1}Z_{1i}-\sum_{j\ne 2}Z_{2j}\\
		&=\sum_{i=1}^{n-1}\theta_i-\theta_1-\sum_{j=1}^{n-2}\theta_j+\sum_{i\notin \{1,2,n-1,n\}}Z_{1i}-\sum_{j\notin\{1,2,n-1,n\}}Z_{2j}\\
		&\quad+Z_{1,n-1}+Z_{1,n}-Z_{2,n-1}-Z_{2,n},
	\end{align*}
	and
	\begin{align*}
&		S_{n}-S_{n-1}\\
		&=\sum_{i=1}^{n-1}\theta_i-\theta_1-\sum_{j=1}^{n-2}\theta_j+\sum_{i\ne n}Z_{n,i}-\sum_{j\ne n-1}Z_{n-1,j}\\
		&=\sum_{i=1}^{n-1}\theta_i-\theta_1-\sum_{j=1}^{n-2}\theta_j+\sum_{i\notin \{1,2,n-1,n\}}Z_{n,i}-\sum_{j\notin \{1,2,n-1,n\}}Z_{n-1,j}\\
		&\quad+Z_{1,n-1}+Z_{2,n-1}-Z_{1,n}-Z_{2,n}.
	\end{align*}
	By the property of normal random variables, we have that $Z_{1,n-1}+Z_{1,n}$ is independent of $Z_{1,n-1}-Z_{1,n}$, and $Z_{2,n-1}+Z_{2,n}$ is independent of $Z_{2,n-1}-Z_{2,n}$. Therefore, noting that $\sum_{i\notin \{1,2,n-1,n\}}Z_{1i}-\sum_{j\notin\{1,2,n-1,n\}}Z_{2j}$ is independent of $\sum_{i\notin \{1,2,n-1,n\}}Z_{n,i}-\sum_{j\notin \{1,2,n-1,n\}}Z_{n-1,j}$, we conclude that the event $S_1-S_2\ge 0$ is independent of $S_n-S_{n-1}\ge 0$. 
	This completes the proof of the proposition.
\end{proof}

\section*{Acknowledgement}

We would like to thank the Editor, the Associate Editor, and three  anonymous referees for their helpful comments and suggestions on the previous version of the manuscript, which led to significant improvement of the paper. This work was partially done when R.M. was a PhD candidate in biostatistics at the University of Pennsylvania, supported by Professor Hongzhe Li, and a postdoctoral scholar in the Department of Statistics at Stanford University, hosted by Professor David Donoho.

\bibliographystyle{IEEEtran}
\bibliography{reference}

\vspace{1cm}

{\bf T. Tony Cai} received the Ph.D. degree from Cornell University, Ithaca, NY, USA, in 1996. He is currently the Daniel H. Silberberg Professor of statistics and data science at the Wharton School, University of Pennsylvania, Philadelphia, PA, USA. His research interests include statistical machine learning, high-dimensional statistics, large-scale inference, nonparametric function estimation, functional data analysis, and statistical decision theory. He is a fellow and President-Elect of the Institute of Mathematical Statistics. He was a recipient of the 2008 COPSS Presidents Award. He is a Past Editor of the Annals of Statistics. 

{\bf Rong Ma} is currently an Assistant Professor of biostatistics at Harvard T.H. Chan School of Public Health of Harvard University. He received his Ph.D. in biostatistics from the University of Pennsylvania, and was a postdoctoral scholar in statistics at Stanford University. His current research focuses on statistical inference for large random matrices, embedding theory, and manifold learning for biomedical research, especially single-cell genomics and multiomics. He was a recipient of the 2022 Lawrence D. Brown Ph.D. Student Award from the Institute of Mathematical Statistics. 

\end{document}